\documentclass[reqno]{amsart}
\usepackage[pdfencoding=auto, psdextra, hidelinks, pdfusetitle]{hyperref}
\usepackage{preamble/packages}
\usepackage{preamble/environments}
\usepackage{preamble/notation}
\usepackage[margin=2.9cm]{geometry}

\graphicspath{{./image/}}
\title{Nakajima quiver varieties in dimension four}
\author{Samuel Lewis and Pavel Shlykov}
\date{}

\begin{document}
\begin{abstract}
  This paper classifies all 4d Nakajima quiver varieties through a combinatorial approach. For each such variety, we describe the symplectic leaves and minimal degenerations between them. Using the resulting Hasse diagrams and secondary hyperplane arrangements, we fully classify the quiver varieties up to isomorphism, a step in the problem of classifying all 4d conical symplectic singularities and the \((2,2)\) case of quiver varieties. As an application, we answer in the negative a question posed by Bellamy, Craw, Rayan, Schedler, and Weiss regarding whether the \(G_4\) quotient singularity (or its projective crepant resolutions) can be realised as a quiver variety.
\end{abstract}
  \maketitle
  \setcounter{tocdepth}{1}
  \tableofcontents
    
  \section{Introduction}\label{sec:intro}
  Let \(\quiver=(\quiver[0],\quiver[1],s,t)\) be a finite, connected quiver of rank \(\rank\defn|\quiver[0]|\) and let \(\root\in\Int[0]<\rank>\) be a dimension vector. The vector space \(\Rep{\quiver,\root}\) of complex representations of \(\quiver\) with dimension vector \(\root\) is an affine variety, and its cotangent bundle admits a Hamiltonian action by the reductive gauge group \(G(\root)\defn\Prod{\GL[\root[i]]{\Comp}}[i\in\quiver[0]]\), with moment map $\mom$. The Nakajima quiver variety
\begin{equation*}
    \QV{\root}\defn\mom\inv(0)\git G(\root)
\end{equation*}
is then obtained by Hamiltonian reduction. This is a coarse moduli space parametrising \(S\)-equivalence classes of representations (with dimension vector \(\root\)) of the associated preprojective algebra \(\Uppi(\quiver)\), the semisimple representations corresponding to closed points.

Quiver varieties appear in many different parts of algebraic geometry and representation theory. For example, they give constructions of moduli spaces such as resolutions of Kleinian singularities, their Hilbert schemes of points, and Uhlenbeck and Gieseker instanton moduli spaces.

By work \cite{bellamy_symplectic_2021} of Bellamy and Schedler, \(\QV{\root}\) is a conical symplectic singularity, admitting projective symplectic (equivalently, crepant) resolutions in precise cases. These singularities are also an important topic in modern representation theory and algebraic geometry. 
They allow one to obtain many properties of deformation quantizations of a variety by studying the resolutions, enhance the search for compact symplectic varieties, and possess many fascinating properties, such as Białynicki-Birula decomposition or symplectic duality\footnote{We refer the reader to Kamnitzer's survey \cite{kamnitzer_symplectic_2022}.}, closely related to duality for $\on{3d}$ supersymmetric quantum field theories.

Quiver varieties make up a large part of these singularities, so it is by studying them we begin to classify all conical symplectic singularities in 4d. This is similar to the work of Nagaoka (see \cite[Thm. 1.4]{nagaoka_universal_2021}) on classifying affine 4d hypertoric varieties, another class of conical symplectic singularities.

In 2d, the quiver varieties are the well-studied Kleinian singularities. Here, \(\quiver\) has underlying graph \(\ade*\), for \(\ade\) a simply laced Dynkin diagram (\(\A\D\E\)), and \(\root\) is the associated minimal (\emph{isotropic}) imaginary root \(\mir(\ade)\) --- see \cref{ntn:mirhrr} for the list. Combinatorially, this dimension vector is special in the sense that for all \(i\in\quiver[0]\) we have
\begin{equation*}
    2\mir[i]-\Sum{\mir[j]}[i\text{---}j]=0.
\end{equation*}
It turns out that \(\dim\QV{\root}\) can be determined by generalising this idea to any dimension vector \(\root\).

To see this, let \(\graph=(\graph[0],\graph[1])\) be the underlying graph of \(\quiver\) and let \(i\) be a vertex with \(\ell_i\) loops. The \say{vertex} and \say{total} \emph{balances} of a dimension vector \(\root\in\Nat<\rank>\) are

\begin{equation*}
    (\root,\srt[i]) = 2(1-\ell_i)\root[i] - \Sum{\root[j]}[i\text{---}j], \qquad (\root,\root) = \Sum{\root[i](\root,\srt[i])}[i],
\end{equation*}

where \((\cdot,\cdot)\) is the bilinear Euler form associated with the generalised Cartan matrix of \(\graph\). We call a vertex \(i\) \emph{balanced} if \((\root,\srt[i]) = 0\) and call \(\root\) itself balanced if every vertex is balanced. In \cite{kac_infinite_1982}, Kac first studied the function
\begin{equation*}
    \p\colon\Int<\rank>\to\Int, \quad \root\mapsto1-\tfrac{1}{2}(\root,\root),
\end{equation*}
to count parameters in indecomposable representations of \(M\in\Rep{\quiver,\root}\). Following this, Crawley-Boevey in \cite[Thm. 1.1]{crawley-boevey_geometry_2001} shows that \(\dim\QV{\root}=2\p(\root)\), provided that \(\root\) is in a combinatorially defined (\cref{dfn:sig0}) subset \(\Sig\) of the positive root system \(\Root\) associated with \(\graph\). By \cite[Thm. 1.2]{crawley-boevey_geometry_2001}, \(\Sig\) consists of dimension vectors for \(0\)-stable (that is, simple) representations of the preprojective algebra \(\Uppi(\quiver)\). Geometrically, the functions that define the variety \(\QV{\root}\) will be \say{independent} when \(\root\in\Sig\).

The importance of \(\Sig\) is shown in \cite[Thm. 1.1, Ppn. 1.2]{crawley-boevey_decomposition_2002}, wherein every \(\root\in\Nat\Root\) admits a \emph{canonical decomposition} (see \eqref{eqn:canonical}) in terms of some \(\root*<k>\in\Sig\). Geometrically, \(\QV{\root}\) is then the product of the \(m_k^\text{th}\) symmetric power \(\Sym<m_k>(\QV{\root*<k>})\), where \(m_k\) is the coefficient of \(\root*<k>\) in the canonical decomposition. It follows that to classify all quiver varieties \(\variety\) with \(\dim\variety=2n\) one must find all the dimension vectors $\root\in\Sig$ with $\p(\root)\leq n$. In this paper we study the case \(n=2\) in order to classify 4d quiver varieties.

\subsection{Main results}\label{ssec:results}

We start by finding all dimension vectors \(\root\) such that

\begin{itemize}
    \item \(\p(\root)=2\), that is, total balance \(-2\),
    \item \((\root,\srt[i]) \leq0\), that is, nonpositive vertex balance, for all \(i\in\graph\).
\end{itemize}

The second condition (which turns out to be a consequence of \(\root\in\Sig\), see \cref{lem:fset} and \cite[Lem. 2.3]{crawley-boevey_decomposition_2002}) makes the problem finite, as there are only three qualitatively different cases:

\begin{itemize}
    \item[(\(\I\))] There are only two unbalanced vertices, both of weight \(1\) and balance \(-1\).
    \item[(\(\II\))] There is a unique unbalanced vertex, of weight \(1\) and balance \(-2\).
    \item[(\(\III\))] There is a unique unbalanced vertex, of weight \(2\) and balance \(-1\).
\end{itemize}

Across these three \say{Types} we give all options for an unbalanced vertex (see local structure, \cref{ssec:locstr}), and classify all ways they might be furnished with balanced vertices. Finding every dimension vector \(\root\) across Types \(\I\), \(\II\), and \(\III\), and removing those that lie outside \(\Sig\), we obtain our first result.

\begin{thm}\label{thm:dimvec}
    The following is a full list of dimension vectors \(\root\) such that \(\root\in\Sig\) and \(\p(\root)=2\). Here \(\ade\) is a simply laced Dynkin diagram and \(\hrr(\ade)\) the highest root in the associated (finite) root system, so that dotted boxes contain the corresponding minimal imaginary root. Unbalanced vertices are underlined, with all other vertices assumed balanced.
    \begin{longtable}{|C|C|C|}\caption{Dimension vectors \(\root\in\Sig\) with \(\p(\root)=2\)}\label{tbl:dimvec} \\ \hline
    \rowcolor{gray!50} \multicolumn{2}{|C|}{\root} & \rank, \height\footnotemark \\\hline\endfirsthead
    \caption{(continued)}\\\hline\rowcolor{gray!50} \multicolumn{2}{|C|}{\root} & \rank, \height \\\hline\endhead
    \hline\endfoot\footnotetext{We have included \(\height\) as it coincides with the dimension of the Coulomb branch for the associated quiver gauge theory.}
    \makecell{\I(\ade,n)\\\ade\neq\E[8]\\n\geq0} & \begin{tikzcd}[sep=1.5em, /tikz/execute at end picture={
    \node (large) [
    rectangle, draw=black,dotted, fit=(A1) (A2) (A3)] {};}]
            |[alias=A1]| \bal{1} \ar[r] \& |[alias=A2]| \hrr(\ade) \ar[r] \& |[alias=A3]| \bal{1} \\
            \& \cdots \ar[ul, bend left, "1"] \ar[ur, bend right, "n"']
        \end{tikzcd} & \makecell{\rank(\ade)+n\\\height(\ade)+n} \\
        \cellcolor{gray!20}\colorbox{gray!20}{\makecell{$\I(m)$\\$m\geq4$}} & \cellcolor{gray!20}\begin{tikzcd}[sep=1.5em]
        \&\& \bal{1} \ar[d] \&\& \bal{1} \ar[d] \\
        1 \ar[r] \& 2 \ar[r] \& 3 \ar[r, "1"] \& \cdots \ar[r, "m-4"] \& 3 \ar[r] \& 2 \ar[r] \& 1
        \end{tikzcd} & \cellcolor{gray!20}\colorbox{gray!20}{\makecell{$m+3$\\$3m-1$}} \\
        \I[a] & \begin{tikzcd}[sep=1.5em]
        \bal{1} \ar[dr] \&\& 2 \ar[d] \\
        \& 3 \ar[r] \& 4 \ar[r] \& 3 \ar[r] \& 2 \ar[r] \& 1 \\
        \bal{1} \ar[ur]
        \end{tikzcd} & \makecell{8\\17} \\
        \cellcolor{gray!20} \I[b] & \cellcolor{gray!20} \begin{tikzcd}[sep=1.5em]
        \bal{1} \ar[r] \& 3 \ar[dr] \\
        \&\& 5 \ar[r] \& 4 \ar[r] \& 3 \ar[r] \& 2 \ar[r] \& 1 \\
        \bal{1} \ar[r] \& 3 \ar[ur]
        \end{tikzcd} & \cellcolor{gray!20}\colorbox{gray!20}{\makecell{$9$\\$23$}} \\
        \I[c] & \begin{tikzcd}[sep=1.5em]
        \bal{1} \ar[dr] \&\&\&\& 3 \ar[d] \\
        \& 3 \ar[r] \& 4 \ar[r] \& 5 \ar[r] \& 6 \ar[r] \& 4 \ar[r] \& 2 \\
        \bal{1} \ar[ur]
        \end{tikzcd} & \makecell{9\\29} \\
           \cellcolor{gray!20} \II(\ade<1>,\ade<2>) & \cellcolor{gray!20} \begin{tikzcd}[sep=1.5em, /tikz/execute at end picture={
    \node (large) [
    rectangle, draw=black,dotted, fit=(A1) (A2)] {};
    \node (large) [
    rectangle, draw=black,dotted, fit=(A2) (A3)] {};}]
        |[alias=A1]|\hrr(\ade<1>) \ar[r] \& |[alias=A2]|\bal{1} \ar[r] \& |[alias=A3]|\hrr(\ade<2>)
        \end{tikzcd} & \mrow[gray!20]{\rank<(1)>+\rank<(2)>+1 \\ \mathsf{ht}^{(1)}+\mathsf{ht}^{(2)}+1} \\
            \mrow{\II(m)\\m\geq4} & \begin{tikzcd}[sep=1.5em]
            2 \ar[dr] \&\&\& \bal{1} \ar[d] \\
            \& 4 \ar[r, "1"] \& \cdots \ar[r, "m-4"] \& 4 \ar[r] \& 3 \ar[r] \& 2 \ar[r] \& 1 \\
            2 \ar[ur]
        \end{tikzcd} & \mrow{m+3\\4m-1} \\
            \cellcolor{gray!20}\II[a] & \cellcolor{gray!20}\begin{tikzcd}[sep=1.5em]
            \& 2 \ar[d] \&\& 3 \ar[d] \\
            \bal{1} \ar[r] \& 4 \ar[r] \& 5 \ar[r] \& 6 \ar[r] \& 4 \ar[r] \& 2 
        \end{tikzcd} & \mrow[gray!20]{8\\27} \\
    \II[b] & \begin{tikzcd}
        \&\&\& 5\ar[d] \\ \bal{1}\ar[r]\&4\ar[r]\&7\ar[r]\&10\ar[r]\&8\ar[r]\&6\ar[r]\&4\ar[r]\&2
    \end{tikzcd} & \mrow{9\\47} \\
            \mrow[gray!20]{\III(\ade)\\\ade\neq\A} & \cellcolor{gray!20}\begin{tikzcd}[sep=1.5em, /tikz/execute at end picture={
    \node (large) [
    rectangle, draw=black,dotted, fit=(A1) (A2)] {};}]
            |[alias=A1]| \hrr(\ade) \ar[r] \& |[alias=A2]| \bal{2} \ar[r] \& 1
        \end{tikzcd} & \mrow[gray!20]{\rank(\ade)+2 \\ \mathsf{ht}(\ade)+2} \\
        \mrow{\III(m, n)\\m\geq4\\n\in\set{6,7,8}} & \begin{tikzcd}[sep=1.5em, /tikz/execute at end picture={
    \node (large) [label = south:{\(\hrr(\E[n])\)}, 
    rectangle, fit=(A1) (A2)] {};}]
            1 \ar[dr] \\
            \& 2 \ar[r, "1"] \& \cdots \ar[r, "m-4"] \& |[alias=A1]| \bal{2} \ar[r] \& 3 \ar[r] \& |[alias=A2]| \cdots \\
            1 \ar[ur]
        \end{tikzcd} & \mrow{m+n-2\\2m+\set{5,11,24}} \\
        \rowcolor{gray!20}\III[a] & \cellcolor{gray!20}\begin{tikzcd}[sep=1.5em]
            \&\&\&\& \bal{2} \ar[d] \\
            1 \ar[r] \& 2 \ar[r] \& 3 \ar[r] \& 4 \ar[r] \& 5 \ar[r] \& 4 \ar[r] \& 3 \ar[r] \& 2 \ar[r] \& 1
        \end{tikzcd} & \mrow[gray!20]{10\\27} \\
        \III[b] & \begin{tikzcd}[sep=1.5em]
            \&\&4\ar[d] \\ \bal{2}\ar[r]\&5\ar[r]\&8\ar[r]\&7\ar[r]\&6\ar[r]\&5\ar[r]\&4\ar[r]\&3\ar[r]\&2\ar[r]\&1
        \end{tikzcd} & \mrow{11\\47} \\
    \end{longtable}
\end{thm}

\begin{rmk}
    Compare this to \cite[Thm. 5.1, 5.3]{schedler_symplectic_2022}, where the authors classified special \say{star-shaped} and \say{crab-shaped} \(\root\) with $\p(\root)=2$. All appear in our classification, mostly in the form of exceptions and minimal cases of infinite families, and some lie outside $\Sig$, their quiver varieties being symmetric powers of Kleinian singularities.
\end{rmk}

Bellamy and Schedler show in \cite{bellamy_symplectic_2021} that all possible subrefinements of the canonical decomposition classify symplectic leaves $\leaf$ within quiver varieties, and in some cases (importantly, our own) one can say something about the minimal degenerations between them. We classify them using \emph{isotropic decompositions}, which crucially when \(\p(\root)=2\) are all of the form
\begin{equation}\label{eqn:decomp}
    \root = \mir(\ade) + \hrr(\ade')
\end{equation}
for some \(\A\D\E\) \(\ade,\ade'\). Geometrically, the data $(\mir(\ade), \hrr(\ade'))$ tells us that the transversal slice to the corresponding 2d leaf \(\leaf\) inside \(\QV{\root}\) is the Kleinian singularity corresponding to \(\ade'\), and the transversal slice to \({0}\) inside \(\leaf\) (that is, the closure \(\cl{\leaf}\)) is the Kleinian singularity corresponding to \(\ade\). This, together with the number of projective crepant resolutions, allow us to distinguish the quiver varieties up to isomorphism.

\begin{thm}\label{thm:quivvar}
    Every dimension vector in \cref{thm:dimvec} yields a non-isomorphic conical symplectic singularity. In particular, \(\QV{\II(\ade<1>,\ade<2>})\) is a product of Kleinian singularities, and we have no examples of symmetric powers arising from dimension vectors in \(\Sig\).
\end{thm}

\begin{rmk}
    Note, that in our classification all the roots in $\Sig$ produce distinct quiver varieties. However, for higher dimensions this is not necessarily so. In \cite[Lem. 8.4]{bellamy_singularities_2025} the authors explain that for $\root<1>, \root<2> \in \Sig$ one can find $\root*$ in $\Sig$ such that $\QV{\root*} \simeq \QV{\root<1>} \times \QV{\root<2>}$ if both $\root<1>, \root<2>$ have a vertex of weight $1$. Namely, one needs to \say{glue} the roots by identifying this vertex in both of them, a generalisation of \(\II(\ade<1>,\ade<2>)\) beyond Kleinian singularities. In principle, one can do this in several different ways if one of the roots has more than one vertex of weight \(1\) and these dimension vectors will all correspond to the same variety --- one can see this already in the case $\p(\root)=3$.
\end{rmk}

We thus obtain a full classification of 4d Nakajima quiver varieties along with their Namikawa Weyl groups, the latter conveniently being the product of the finite Weyl groups of \(\A\D\E\) type that appear as highest roots in isotropic decompositions \eqref{eqn:decomp}.

\begin{thm}\label{thm:namikawa_weyl_groups_are_classified}
    A 4d Nakajima quiver variety has one of the following Namikawa Weyl groups:
    \rowcolors{1}{}{gray!20}
    \begin{longtable}{|C|C|C|C|}\caption{Namikawa Weyl groups for 4D quiver varieties}\label{tbl:quivvar} \\\hline
        \rowcolor{gray!50} \root & \textup{Parameter(s)} & \Weyl[\root] & |\Weyl[\root]| \\\hline\endfirsthead
        \caption{(continued)}\\\hline\rowcolor{gray!50} \root & \textup{Parameter(s)} & \Weyl[\root] & |\Weyl[\root]| \\\hline\endhead
        \hline\endfoot
        \I(\ell,m,n) & \ell,m,n\geq0 & \Weyl(\A[\ell-1])\times\Weyl(\A[m-1])\times\Weyl(\A[n-1]) & \ell! m!n! \\
        \I(\D[m], n) & m\geq4,\ n\geq0 & \Weyl(\A[m-1])\times\Weyl(\A[n-1]) & m!n! \\
        \I(\D[m], n)' & m\geq4,\ n\geq0 & \Weyl(\A[n-1])\times\Weyl(\D[m-1]) & 2^{m-2}(m-1)!n! \\
        \I(\E[6], n) & n\geq0 & \Weyl(\A[n-1])\times\Weyl(\D[5]) & 1920n! \\
        \I(\E[7], n) & n\geq0 & \Weyl(\A[n-1])\times\Weyl(\E[6]) & 51840n! \\
        \I(m) & m\geq4 & \Weyl(\A[m+1]) & (m+2)! \\
        \I[a] & & \Weyl(\D[6]) & 23040 \\
        \I[b] & & \Weyl(\D[7]) & 322560 \\
        \I[c] & & \Weyl(\E[7]) & 2903040 \\
        \II(\ade<1>,\ade<2>) & \ade<1>,\ade<2>\in\set{\A,\D,\E} & \Weyl(\ade<1>)\times\Weyl(\ade<2>) & |\Weyl(\ade<1>)||\Weyl(\ade<2>)| \\
        \II(m) & m\geq4 & \Weyl(\D[m+2]) & 2^{m+1}(m+2)! \\
        \II[a] & & \Weyl(\E[7]) & 2903040 \\
        \II[b] & & \Weyl(\E[8]) & 696729600 \\
        \III(\D[m], i) & 1\leq i\leq\ceil{\tfrac{m-3}{2}} & \Weyl(\A[1])\times\Weyl(\D[i+1])\times\Weyl(\D[m-i-1]) & 2^{m-1}m!\binom{m}{i+1} \\
        \III(\E[7]) & & \Weyl(\A[1])\times\Weyl(\A[7]) & 80640 \\
        \III(\E[8]) & & \Weyl(\A[1])\times\Weyl(\D[8]) & 10321920 \\
        \III(m, 6) & m\geq4 & \Weyl(\A[5])\times\Weyl(\D[m-2]) & 2^m90(m-2)! \\
        \III(m, 7) & m\geq4 & \Weyl(\D[6])\times\Weyl(\D[m-2]) & 2^m2880(m-2)! \\
        \III(m, 8) & m\geq4 & \Weyl(\E[7])\times\Weyl(\D[m-2]) & 2^m362880(m-2)! \\
        \III[a] & & \Weyl(\A[9]) & 3628800 \\
        \III[b] & & \Weyl(\D[10]) & 1857945600 \\
        \III(2\ade) & \ade\in\set{\A,\D,\E} & \Weyl(\ade)\times\Weyl(\A[1]) & 2|\Weyl(\ade)|
    \end{longtable}
\end{thm}

One of the main consequences of this work are the first steps in understanding the \((2,2)\) case of quiver varieties, the only (see \cite[Thm. 1.5]{bellamy_symplectic_2021}) quiver varieties with \emph{divisible} dimension vectors that admit projective symplectic resolutions. These are all of dimension \(10\) and take the form \(\QV{2\root}\), where \(\root\in\Sig\) satisfies \(\p(\root)=2\), so we have found all \((2,2)\) dimension vectors.

The classification is also connected with conical symplectic singularities arising in other areas, such as nilpotent orbit closures, hypertoric varieties, and finite group quotients. For the last area, we answer a question posed in \cite{bellamy_all_2024} regarding the exceptional complex reflection group \(G_4\). 

\begin{cor}\label{cor:g4}
    The \(G_4\) quotient singularity cannot be realised as a quiver variety, nor can any of its projective symplectic resolutions.
\end{cor}

In the proof of \cref{thm:quivvar} and \cref{cor:g4} we count the number of projective symplectic resolutions admitted by our quiver varieties, using the secondary and Coxeter arrangements associated with the Namikawa Weyl group. In this way we establish unexpected conjectures about these numbers for the cases which currently exceed available computational capabilities. In particular, we make the claim that every quiver variety with dimension vector of Type \(\II\) admits a unique projective symplectic resolution.

\subsection{Structure of the paper}

In \cref{sec:prelim} we give a background on quiver varieties, root combinatorics, and their interplay. Importantly, we summarise the 2d classification of quiver varieties, and provide balancing results that will aid us in the 4d classification.

These dimension vectors are given in \cref{sec:roots}, starting from the \say{local structure} of unbalanced vertices and working through all possible balanced terminations. We then rule out any such dimension vectors that lie outside of \(\Sig\), establishing \cref{thm:dimvec}.

In \cref{sec:symplectic} we go into more detail on the geometric structure admitted by quiver varieties, in order to distinguish $\QV{\root}$ for all $\root$ obtained in \cref{ssec:class}. In particular we compute all of their isotropic decompositions, summarised in \cref{tbl:decomp}. As a consequence of this work we may classify all Namikawa Weyl groups that arise for 4d quiver varieties, which is \cref{thm:namikawa_weyl_groups_are_classified}.

We complete the classification in \cref{sec:varieties}, proving \cref{thm:quivvar} using the invariants of minimial degenerations and projective crepant resolutions. This allows us to prove \cref{cor:g4} and state several conjectures for the number of such resolutions.

The symplectic leaf closures, and hence full poset of transversal slice singularities, are tabulated in \cref{app:hasse} to proivde another way of viewing the results of \cref{ssec:decompositions}.

\subsection{Notation}\label{ssec:notation}

For \(\root,\root*\in\Nat<\rank>\) we write \(\root\geq\root*\) if \(\root[i]\geq\root*[i]\) for all \(i\) and write \(\root\geq*\root*\) if further \(\root\neq\root*\). In general two dimension vectors will be supported on different graphs with their own root systems, but we abuse notation and write \(\Sig\) with respect to both. By \emph{weighted graph}, we mean a connected graph \(\graph\) with nonnegative integer weights \(\root[i]\) on the vertices. We draw such graphs with weights directly on each vertex and without explicit indexing, in this way using \(\root\) to infer the \(\graph\) that supports it. In similar vein, we write $\QV{\root}$ in place of $\QV{Q,\root}$ or $\QV{\graph,\root}$. In the displayed dimension vectors, vertices are assumed to be balanced (in the sense of \cref{dfn:bal}) unless underlined or overlined, which denotes negative or positive balance, respectively. This allows the adjacency of vertices to be inferred in the compact versions of dimension vectors used going forward. We use some graph theory notions like a \emph{branching vertex}, which has more than two neighbours, and a \emph{hanging vertex}, which only has one neighbour.

Our parameter conventions are \(m\geq0\) for \(\A[m]\) and \(m\geq4\) for \(\D[m]\) unless otherwise stated, where \(\A[0]\) is the empty graph and its extension \(\A*[0]\) consists of a single vertex with a loop. We now give an explanation of the general constructions shown in \cref{tbl:dimvec}, using the words \emph{affine} and \emph{extended} interchangeably with respect to \(\A\D\E\) Dynkin diagrams.

\begin{ntn}\label{ntn:chain}
    Let \(\ade\) be of type \(\A\D\E\), choose two vertices of weight \(1\) in its minimal imaginary root \(\mir(\ade)\), and let \(n\geq*0\). Then \(\I(\ade,n)\) is the dimension vector obtained by connecting these two vertices with a chain of \(n\) edges and weight \(1\) vertices, which gives both distinguished vertices a balance of \(-1\). There are some subtleties for specific \(\ade\):
    \begin{itemize}
        \item[(\(\A\))] We may view \(\I(\A[m],n)\) as a three-parameter family \(\I(\ell,m,n)\), consisting of two unbalanced vertices of weight \(1\) connected by three strings of balanced weight \(1\) vertices, of length \(\ell,m,n\geq1\). Note that \(\I(\A[m],0)\) coincides with some \(\II(\A, \A')\) described below, so we restrict to \(1\leq*\ell\leq m\leq n\).
        \item[(\(\D\))] There is choice of which pair of weight \(1\) vertices to chain. If they neighbour different branch vertices, then we form a chain of length \(m+n-2\) and denote the dimension vector by \(\I(\D[m],n)\). If they neighbour the same branch vertex, then we denote this by \(\I(\D[m],n)'\), which has a cycle of length \(n+2\).
        \item[(\(\E\))] There is a unique choice for \(\I(\E[6],n)\) and \(\I(\E[7],n)\), as in each case \(\mir(\ade)\) has only two vertices of weight \(1\). Since we require \(\mir(\ade)\) to have two vertices of weight \(1\), the dimension vector \(\I(\E[8],n)\) is not defined.
    \end{itemize}
    Similarly we may choose two \(\A\D\E\) graphs \(\ade<1>,\ade<2>\) and chain together the unique (up to symmetry) weight \(1\) vertices in \(\mir(\ade<1>)\) and \(\mir(\ade<2>)\). This is denoted \(\I(\ade<1>,\ade<2>,n)\). In both of the above constructions, if \(n=0\) then the distinguished weight \(1\) vertices are identified, giving Type \(\II\) dimension vectors \(\II(\ade)\) and \(\II(\ade<1>,\ade<2>)\), respectively.
    \begin{equation*}
        \setlength\arraycolsep{3em}
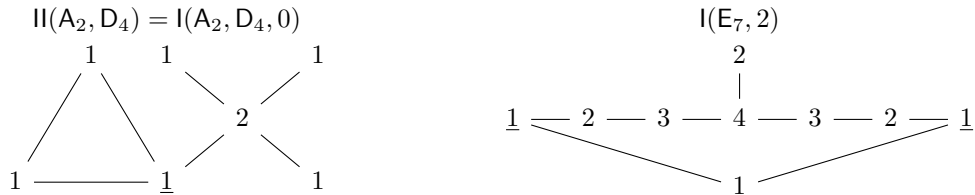
        \begin{array}{cc}
            \II(\A[2],\D[4])=\I(\A[2],\D[4],0) & \I(\E[7],2) \\
            \begin{tikzcd}[row sep=1em, column sep=1.5em]
                \& 1 \ar[ddl] \ar[ddr] \& 1 \&\& 1 \\ \&\&\& 2 \ar[ul] \ar[ur] \ar[dl] \ar[dr] \\ 1 \ar[rr] \&\& \bal{1} \&\& 1
            \end{tikzcd} & \begin{tikzcd}[row sep=1em, column sep=1.5em]
                \&\&\& 2 \ar[d] \\
                \bal{1} \ar[r] \& 2 \ar[r] \& 3 \ar[r] \& 4 \ar[r] \& 3 \ar[r] \& 2 \ar[r] \& \bal{1} \\ \&\&\& 1 \ar[ulll] \ar[urrr]
            \end{tikzcd}
        \end{array}
    \end{equation*}
\end{ntn}

\begin{ntn}\label{ntn:extend}
    Choose a vertex \(i\) of weight \(2\) in \(\mir(\ade)\) and unbalance it by attaching a weight \(1\) vertex.
    \begin{itemize}
        \item[(\(\A\))] There are no vertices of weight \(2\) in \(\mir(\A)\), so \(\III(\A[m])\) is not defined.
        \item[(\(\D\))] Up to symmetry, there are \(\ceil{\tfrac{m-3}{2}}\) unique choices of a weight \(2\) vertex, up to and including the centre of the chain.
        \item[(\(\E\))] For each \(\E[n]\) there is an option of weight \(2\) vertex that becomes a branch vertex when extended, we denote the resulting dimension vector \(\III(\E[n])'\). For \(\E[7]\) and \(\E[8]\) there is also a choice that does not create a new branch vertex, which we denote \(\III(\E[n])\).
    \end{itemize}
    By instead considering \(2\mir(\ade)\), which for any \(\ade\) has a unique weight \(2\) vertex that can be unbalanced, we define the dimension vector \(\III(2\ade)\).
    \begin{equation*}
        \setlength\arraycolsep{3em}
        \begin{array}{cc}
            \III(\D[6],2) & \III(2\A[3]) \\
            \begin{tikzcd}[row sep=1em, column sep=1.5em]
                1 \ar[dr] \&\& 1 \ar[d] \&\& 1 \ar[dl] \\ \& 2 \ar[r] \& \bal{2} \ar[r] \& 2 \\ 1 \ar[ur] \&\&\&\& 1 \ar[ul]
            \end{tikzcd} & \begin{tikzcd}[row sep=1em, column sep=1.5em]
                \& \bal{2} \ar[ddl] \ar[d] \ar[ddr] \\ \& 1 \\ 2 \ar[r] \& 2 \ar[r] \& 2
            \end{tikzcd}
        \end{array}
    \end{equation*}
\end{ntn}

\begin{rmk}
    Note that \(\III(\D[4])\defn\III(\D[4],1)\) is the \say{star} graph with five vertices and the corresponding variety is a special case (\(n=5\)) of the \emph{hyperpolygon space} (see \cite{bellamy_all_2024}).
\end{rmk}

\begin{ntn}\label{ntn:wedgetwo}
    Let \(m\geq4\) and \(n\in\set{6,7,8}\). The dimension vector \(\III(m,n)\) is constructed by taking \(\mir(\D[m])\) and identifying one of its branch vertices with the unbalanced vertex in \(\hrr(\E[n])\), and then removing the neighbours of the chosen branch vertex. Observe \(\III(4,n)=\III(\E[n])'\) for each \(n\), so in \cref{ntn:extend} we only need to consider \(\III(\E[7])\) and \(\III(\E[8])\) as separate dimension vectors.
    \begin{equation*}
        \setlength\arraycolsep{3em}
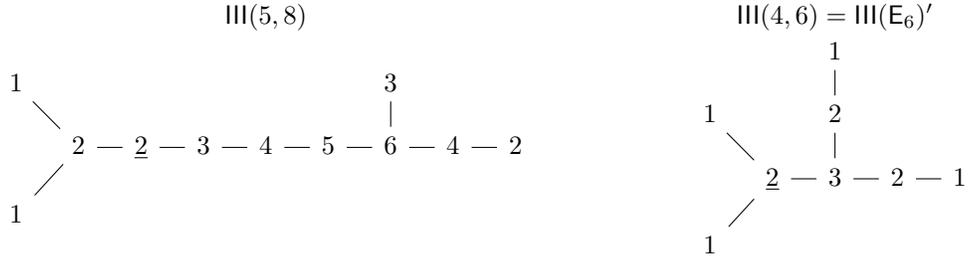
        \begin{array}{cc}
             \III(5,8) & \III(4,6)=\III(\E[6])' \\
             \begin{tikzcd}[sep=1em]
                 1 \ar[dr] \&\&\&\&\&\& 3 \ar[d] \\ \& 2 \ar[r] \& \bal{2} \ar[r] \& 3 \ar[r] \& 4 \ar[r] \& 5 \ar[r] \& 6 \ar[r] \& 4 \ar[r] \& 2 \\ 1 \ar[ur]
             \end{tikzcd} & \begin{tikzcd}[sep=1em]
                 \&\& 1 \ar[d] \\ 1 \ar[dr] \&\& 2 \ar[d] \\ \& \bal{2} \ar[r] \& 3 \ar[r] \& 2 \ar[r] \& 1 \\ 1 \ar[ur]
             \end{tikzcd}
        \end{array}
    \end{equation*}
\end{ntn}

\subsection*{Acknowledgements}

The first author was partially supported by a Research and Teaching Allowance from St John's College, University of Oxford. The second author was supported by the Engineering and Physical Sciences Research Council [grant number EP/W013053/1]. Both authors would like to thank Gwyn Bellamy for suggesting this problem and for many useful discussions thereafter. The first author would like to thank Nick Proudfoot for his advice on hypertoric varieties. The second author would like to thank Vasily Krylov for many helpful conversations. We also thank Antoine Bourget, Alastair Craw, Alastair King, and Lewis Topley for comments on the presentation of the material and suggesting related results. Finally, many thanks to Travis Schedler for pointing out \cref{lem:fnotsig} to us, streamlining the argument regarding \(\Sig\).

\subsection*{Open access} For the purpose of open access, the authors have applied a Creative Commons Attribution (CC:BY) licence to any Author Accepted Manuscript version arising from this submission.

  \section{Preliminaries}\label{sec:prelim}
  Following \cite{bellamy_symplectic_2021}, we give the required background on Nakajima quiver varieties, root system combinatorics, and the case of 2d quiver varieties. We then start the analysis of 4d quiver varieties by giving balancing results that will underpin the classification of dimension vectors.

\subsection{Quiver varieties} \label{ssec:quivvar_definitions}

Take a quiver \(\quiver\) as before, which has vertex set \(\quiver[0]\), arrow set \(\quiver[1]\), and arrow source and target maps \(s,t\colon\quiver[1]\to\quiver[0]\). A (complex, finite-dimensional) \emph{representation} \(M\) of \(Q\) is an assignment of vector spaces to vertices and linear maps to arrows. Then 
\begin{equation*}
    \Rep{\quiver,\root}\defn\Prod{\Hom{\Comp<\root[s(a)]>}{\Comp<\root[t(a)]>}}[a\in\quiver[1]]
\end{equation*} 
is the space of representations of \(\quiver\) with dimension vector \(\root\in\Nat<\rank>\), and viewed with the Zariski topology becomes an affine variety. The reductive group \(G(\root)\defn\prod_{i\in\graph}\GL[\root[i]]{\Comp}\) acts on \(\Rep{\quiver,\root}\) by conjugation, and this action is symplectic on the cotangent bundle \(T^\ast\Rep{\quiver,\root}\). Taking \(\double{\quiver}\) to be the doubled quiver of \(\quiver\), we may identify \(T^\ast\Rep{\quiver,\root}\) with \(\Rep{\double{\quiver},\root}\), and consider the moment map \(\mom\colon\Rep{\double{\quiver},\root}\to\Lie{G(\root)}\), identifying the Lie algebra with its dual via the trace pairing. The \emph{Nakajima quiver variety} associated with this data is then the affine quotient
\begin{equation*}
    \QV[\uplambda]{\root}\defn\mom\inv(\uplambda)\git G(\root).
\end{equation*}
We restrict to the case \(\uplambda = 0\), in which case this is a conical symplectic singularity. 

\begin{rmk}
    We note two generalisations of this construction. Firstly, one can also take the corresponding GIT quotient $\QV[\uplambda]{\root}<\stab> \defn \mom\inv(\uplambda)^{\stab} \git G(\root)$ with stability parameters explained in \cite{king_moduli_1994}. This variety also comes with a morphism to $\QV[\uplambda]{\root}$ given by variation of GIT. For classification purposes we study the affine varieties so for the most part these do not appear in the paper.
    Secondly, Nakajima's original construction \cite{nakajima_instantons_1994} is more general and involves also an additional vector $w$ of \emph{framing} weights. However, as observed by Crawley-Boevey in \cite[\S 1]{crawley-boevey_geometry_2001}, every framed quiver variety can be identified with an unframed one and here we consider the unframed case.
\end{rmk}

\subsection{Root combinatorics}\label{ssec:root_combinatorics}

Let \(\graph\) be the underlying graph of \(\quiver\), so that the data \((\graph,\root)\) may be considered as a (nonnegatively) weighted graph. In fact, \cite[Lem. 2.2]{crawley-boevey_noncommutative_1998} tells us that \(\Uppi(\quiver)\) (and hence \(\QV{\root}\)) depends only on \(\graph\), so this is sufficient input data. Write \(\GCM[\graph]\) for the associated \emph{generalised Cartan matrix}, which has diagonal entries \(2(1-\ell_i)\) and off-diagonal entries \(-\graph[ij]\), where \(\ell_i\) is the number of loops at vertex \(i\) and \(\graph[ij]\) is the number of edges between vertices \(i\) and \(j\). This induces the \emph{Euler form} \((\root,\root*)\defn\root\tran\GCM[\graph]\root*\) on \(\Int<\rank>\) which underpins our key combinatorial definitions.

\begin{dfn}\label{dfn:bal}
    Let \((\graph,\root)\) be a weighted graph and let \(i\in\graph\). The \emph{vertex balance} and \emph{total balance} are
    \begin{equation*}
        (\root,\srt[i]) = 2(1-\ell_i)\root[i] - \sum_{j\text{---}i}\root[j], \qquad (\root,\root) = \Sum{\root[i](\root[i],e_i)}[i\in\graph],
    \end{equation*}
    respectively. A vertex is \emph{balanced} if \((\root,\srt[i])=0\), and \emph{unbalanced} (positively or negatively) otherwise.
\end{dfn}

The following balancing condition is elementary, yet fundamental to our analysis.

\begin{lem}\label{lem:half}
  If \(i,j\in\graph\) share an edge, then \((\root,\srt[j])\leq*0\) whenever \(\root[j]\leq*\ceil{\tfrac{\root[i]}{2}}\). In particular, a balanced vertex necessarily has weight at least half of any neighbouring weight.
\end{lem}
\begin{proof}
  By assumption either \(\root[j]\leq\tfrac{\root[i]}{2} - 1\) (\(\root[i]\) even) or \(\root[j]\leq\tfrac{\root[i]}{2} - \tfrac{1}{2}\) (\(\root[i]\) odd). Then \((\root, \srt[j]) \leq 2\root[j] - \root[i]\), which is at most \(-2\) in the first case and at most \(-1\) in the second.
\end{proof}

Each loopfree (\(\ell_i=0\)) vertex \(i\) has a corresponding reflection \(s_i\root\defn\root-(\root,\srt[i])\srt[i]\), which satisfy \(s_i^2=1\), \(s_is_j=s_js_i\) when \(\graph[ij]=0\) and \(s_is_js_i=s_js_is_j\) when \(\graph[ij]=1\). The group \(\Weyl[\graph]\) is generated by these \(s_i\) and is a (generally infinite) Coxeter group.

\begin{dfn}\label{dfn:root}
    We say that \(\root\in\Int<\rank>\) is a \emph{root} if there exists \(w\in\Weyl[\graph]\) such that \(\root=w\root*\) for some \(\root*\) in either \(\set{e_i}[i\in\graph \text{ loopfree}]\) or \(\pm F\), where \(F\defn\set{\root*\in\Nat<\rank>}[\supp{\root*} \text{ connected}, (\root*,\srt[i])\leq0 \text{ for all } i\in\graph]\). In the first case we say that \(\root\) is \emph{real} and in the second case we say that \(\root\) is \emph{imaginary}.
\end{dfn}

Denote the set of roots by \(\Root*\subset\Int<\rank>\), and the positive roots by \(\Root\). See \cite{kac_infinite-dimensional_1990} for more details on root systems. Now, as \(\GCM[\graph]\) is a symmetric matrix with even diagonal entries, the total balance \((\root,\root)\) is even for all \(\root\), so it makes sense to define the function \(\p(\root)\defn1-\tfrac{1}{2}(\root,\root)\), following \cite{kac_infinite_1980}. It is well known that \(\root\in\Root\) is real (respectively imaginary) then \(\p(\root)=0\) (respectively \(\geq*0\)), however this is not a sufficient condition in general. In the latter case we call \(\root\) \emph{isotropic} if \(\p(\root)=1\), and \emph{anisotropic} otherwise. The \emph{height} of a root \(\root\) is \(\height\defn\Sum{\root[i]}[i]\).

\begin{ntn}\label{ntn:mirhrr}
    For \(\ade\) of type \(\A\D\E\), we write \(\hrr(\ade)\) for the highest root in its (necessarily finite) root system, and \(\mir(\ade)\) for the minimal imaginary root in the corresponding affine root system.
    \begin{longtable}{|C|C|C|}\caption{Isotropic and highest roots}\label{tbl:mirhrr}\\\hline
    \rowcolor{gray!50} \ade & \hrr(\ade) & \mir(\ade) \\\hline\endfirsthead
    \caption[]{(continued)}\\\hline\rowcolor{gray!50} \ade & \hrr(\ade) & \mir(\ade) \\\hline\endhead
    \hline\endfoot
    \mrow{\A[m]\\m\geq0} & \begin{tikzcd}[sep=.9em]
    \bal*{1} \ar[r] \& \cdots\cdots \ar[r] \& \bal*{1}
  \end{tikzcd} & \begin{tikzcd}[sep=.9em]
  \& 1 \ar[dl] \ar[dr] \\
    1 \ar[r] \& \cdots\cdots \ar[r] \& 1
  \end{tikzcd} \\
        \mrow[gray!20]{\D[m]\\m\geq4} & \cellcolor{gray!20}\begin{tikzcd}[sep=.9em]
    \&\&\&\& 1 \ar[dl] \\
    1 \ar[r] \& \bal*{2} \ar[r] \& \cdots\cdots \ar[r] \& 2 \\
    \&\&\&\& 1 \ar[ul]
  \end{tikzcd} & \cellcolor{gray!20}\begin{tikzcd}[sep=.9em]
    1 \ar[dr] \&\&\&\& 1 \ar[dl] \\
    \& 2 \ar[r] \& \cdots\cdots \ar[r] \& 2 \\
    1 \ar[ur] \&\&\&\& 1 \ar[ul]
  \end{tikzcd} \\
            \E[6] & \begin{tikzcd}[sep=.9em]
    \&\& \bal*{2} \ar[d] \\
    1 \ar[r] \& 2 \ar[r] \& 3 \ar[r] \& 2 \ar[r] \& 1
  \end{tikzcd} & \begin{tikzcd}[sep=.9em]
    \&\& 1 \ar[d] \\ \&\& 2 \ar[d] \\
    1 \ar[r] \& 2 \ar[r] \& 3 \ar[r] \& 2 \ar[r] \& 1
  \end{tikzcd} \\
            \cellcolor{gray!20}\E[7] & \cellcolor{gray!20}\begin{tikzcd}[sep=.9em]
    \&\& 2 \ar[d] \\
    \bal*{2} \ar[r] \& 3 \ar[r] \& 4 \ar[r] \& 3 \ar[r] \& 2 \ar[r] \& 1
  \end{tikzcd} & \cellcolor{gray!20}\begin{tikzcd}[sep=.9em]
    \&\&\& 2 \ar[d] \\
    1 \ar[r] \& 2 \ar[r] \& 3 \ar[r] \& 4 \ar[r] \& 3 \ar[r] \& 2 \ar[r] \& 1
  \end{tikzcd} \\
            \E[8] & \begin{tikzcd}[sep=.9em]
    \&\& 3 \ar[d] \\
    2 \ar[r] \& 4 \ar[r] \& 6 \ar[r] \& 5 \ar[r] \& 4 \ar[r] \& 3 \ar[r] \& \bal*{2}
  \end{tikzcd} & \begin{tikzcd}[sep=.9em]
    \&\& 3 \ar[d] \\
    2 \ar[r] \& 4 \ar[r] \& 6 \ar[r] \& 5 \ar[r] \& 4 \ar[r] \& 3 \ar[r] \& 2 \ar[r] \& 1
  \end{tikzcd} \\
    \end{longtable}
    The overlined vertices in \(\hrr(\ade)\) are positively unbalanced, and are \say{extended} with a vertex of weight \(1\) in each case to form \(\mir(\ade)\).
\end{ntn}

\begin{lem}\label{lem:pscale}
    If \(\root\in\Int[0]<\rank>\) and \(n\in\Int\), then \(\p(n\root) = n^2(\p(\root)-1)+1\). In particular, \(n\root\) is isotropic whenever \(\root\) is.
\end{lem}
\begin{proof}
    Applying bilinearity and the fact that \((\root,\root)=2-2\p(\root)\), we obtain
    \begin{equation*}
        \p(n\root)=1-\tfrac{1}{2}(n\root,n\root)=1-\tfrac{1}{2}n^2(\root,\root)=1-\tfrac{1}{2}n^2(2-2\p(\root))=n^2(\p(\root)-1)+1.
    \end{equation*}
    The second claim follows by taking \(\p(\root)=1\).
\end{proof}

We now define an  important subset of positive roots.

\begin{dfn}\label{dfn:sig0}
    The set \(\Sig\) consists of all \(\root\in\Root\) such that \(\p(\root)\geq*\Sum{\p(\root*<k>)}[k=1]<r>\) for every decomposition \(\root=\Sum{\root*<k>}[k=1]<r>\) with \(r\geq2\) and \(\root*<k>\in\Root\).
\end{dfn}

If \(\root\in\Sig\) is real then \(\root=\srt[i]\) is some \emph{simple root}. If  \(\root\in\Sig\) is isotropic, then it is the minimal imaginary root for some affine Dynkin diagram. See \cite[Exm. 1 and 2]{bellamy_symplectic_2021} for more details. Crawley-Boevey \cite[Ppn. 1.2]{crawley-boevey_decomposition_2002} also shows that for anisotropic root $\root* \in \Sig$, $n\root*$ is also in $\Sig$. 
By \cite[Thm. 1.2]{crawley-boevey_geometry_2001} being a root $\root \in \Sig$ is equivalent to the existence of a simple representation of the corresponding preprojective algebra (that is, a point in the quiver variety). 
The importance of $\Sig$ is explained in the following result:

\begin{thm}[Thm. 1.1, Ppn. 1.2, \cite{crawley-boevey_decomposition_2002}]\label{thm:canonical_decomposition_and_dim_of_quiver_variety}
    If \(\root*\in\Sig\) then \(\dim\QV{\root*}=2\p(\root*)\). Furthermore, for every \(\root\in\Nat\Root\) there exist \(m_1,\dots,m_k\in\Nat\) and \(\root*<1>,\dots,\root*<k>\in\Sig\) such that
    \begin{equation}\label{eqn:canonical}
        \root=\Sum{m_k\root*<k>}[k=1]<r>.
    \end{equation}
    This decomposition is canonical in the sense that every other decomposition of $\root$ in terms of $\Sig$ is a refinement of \eqref{eqn:canonical}. If the \(\root*<k>\) are distinct, and \(\Sym<m>\) denotes the $m^{\text{th}}$ symmetric power, then 
    \begin{equation}\label{eqn:qvdecomp}
        \QV{\root} \iso \Prod{\Sym<m_k>(\QV{\root*<k>})}[k=1]<r>.
    \end{equation}
    where . Moreover, $\QV{\root*<k>}$ is a point for real $\root*<k>$, and \(m_k=1\) for anisotropic \(\root*<k>\).
\end{thm}

Notice that if $\root \in \Sig$ then the corresponding canonical decomposition is $\root$ itself. Hence the 4d quiver varieties that we aim to classify will have dimension vectors \(\root\in\Sig\) satisfying \(\p(\root)=2\). All other possible quiver varieties of this dimension are products and symmetric powers of Kleinian singularities.

To motivate what follows, we recall the classification of isotropic roots, that is, the case \(\p(\root)=1\). This is equivalent to total balance \((\root,\root)=0\), so every vertex is balanced.

\begin{thm}[Lem. 1.9(d), \cite{kac_infinite_1980}]\label{thm:classification_of_dim2}
    All totally balanced positive roots are given by \(n\mir(\ade)\) for some simply laced Dynkin diagram \(\ade\) and some \(n\in\Int[1]\).
\end{thm}
Of these, only the case $n=1$ belongs to $\Sig$. They are connected and belong to the fundamental set of imaginary roots since every vertex has zero balance. Thus 2d quiver varieties take the form $\QV{\mir(\ade)}$, using the notation of \cref{tbl:mirhrr}. See \cite[Exm. 2]{bellamy_symplectic_2021} for more detail.

\subsection{Balancing lemmas}\label{ssec:balancing}
Our strategy for finding dimension vectors for 4d quiver varieties is based on the following statement.
\begin{lem}[Thm. 5.8, \cite{crawley-boevey_geometry_2001}]\label{lem:fset}
    Suppose \(\root\in\Sig\) is imaginary. Then \(\root\) has connected support and \((\root,\srt[i])\leq0\) for all \(i\). In other words, \(\root\) lies in the fundamental set of positive imaginary roots.
\end{lem}

This is a convenient necessary condition, as first we can classify dimension vectors \(\root\) with \(\p(\root)=2\) and nonpositive local balance at all vertices, and then check if any \(\root\not\in\Sig\).

We first investigate restrictions on the number of neighbours a balanced vertex may have.

\begin{lem}[Branching]\label{lem:branch}
    If \((\graph,\root)\) has a weighted subgraph of the form
    
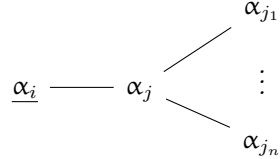
\begin{equation*}\label{eqn:branching}
        \begin{tikzcd}[row sep=0.5em]
            \&\& \root[j_1] \ar[dl] \\
            \bal{\root[i]} \ar[r] \& \root[j] \& \vdots \\
            \&\& \root[j_n] \ar[ul]
        \end{tikzcd}
    \end{equation*}
    for some \(n\in\Int[1]\), then \(\root[i]\leq2\root[j]-n\ceil{\tfrac{\root[j]}{2}}\) and \(n\leq3\). In particular \(\root[j]\geq\root[i]\) if \(n=2\) and \(\root[j]\geq2\root[i]\) (or \(\root[j]\geq2\root[i]+3\) for \(\root[j]\) odd) if \(n=3\).
\end{lem}
\begin{proof}
    Since \(\root[j]\) is balanced, by definition and \cref{lem:half} we have \(2\root[j]-\root[i]=\root[j_1]+\cdots+\root[j_n]\geq n\ceil{\tfrac{\root[j]}{2}}\) and rearranging gives the first claim. If \(n\geq4\), we see that \(\root[i]\leq2\root[j]-4\ceil{\tfrac{\root[j]}{2}}\leq2\root[j]-2\root[j]=0\), which is impossible as \(\root\geq*0\). Similarly, if \(n\geq2\) then \(\root[i]\leq2\root[j]-2\ceil{\tfrac{\root[j]}{2}}\leq2\root[j]-\root[j]=\root[j]\).
\end{proof}

The second part of \cref{lem:branch} implies that in decreasing chains, meaning a sequence of vertices with successively lower weights, branching (\(n\geq2\)) is impossible. Hence we can fully understand balanced decreasing chains that start with an unbalanced vertex.

\begin{lem}[Down-by-\(d\)]\label{lem:stepdown}
    Let \((\graph,\root)\) have a weighted subgraph of the form
    \begin{equation*}\label{eqn:decreasingchain}
        \begin{tikzcd}
            \bal{\root[i]} \ar[r] \& \root[j]\defn\root[i]-d \ar[r] \& \cdots\cdots
        \end{tikzcd}
    \end{equation*}
    for some \(0\leq*d\leq*\root[i]\). If this has a balanced continuation (meaning \(j\) has at least one other balanced neighbour), then \(\root[i]\geq3d\) (or \(\root[i]\geq3d+1\) for $d$ odd), terminating precisely when \(d\) divides \(\root[i]\).
\end{lem}
\begin{proof}
    By \cref{lem:branch} we know that \(j\) has at most one other neighbour, as \(\root[j]\leq*\root[i]\). Since this neighbour is also balanced, it has weight \(\root[i]-2d\geq\ceil{\tfrac{\root[i]-d}{2}}\) and the claimed inequalities follow. Assuming a balanced continuation exists, the next weights are \(\root[i]-2d, \root[i]-3d\) and so on. This eventually terminates with \(\root[i]-kd\) for some \(k\), but then \(2(\root[i]-kd)=\root[i]-(k-1)d\) and rearranging yields \(\root[i] = (k+1)d\).
\end{proof}

Increasing chains are more subtle, as they can grow unbounded and contain branches at multiple points. We consider the increasing steps relevant to our classification, and in each case determine when growing unbounded is the only possibility for the chain. As we are only interested in finite rank, we can then rule out such behaviour from the classification.

\begin{lem}[Up-by-\(d\)]\label{lem:stepup}
    Suppose \((\graph,\root)\) contains the increasing chain
    \begin{equation*}
        \begin{tikzcd}
            \bal{\root[i]} \ar[r] \& \root[j]\defn\root[i]+d \ar[r] \& \cdots
        \end{tikzcd}
    \end{equation*}
    for some \(d\geq0\). Then for vertex \(j\) to admit branching requires
    \begin{enumerate}
        \item[(\(d=0\))] \(\root[i]\) to be even,
        \item[(\(d=1\))] \(\root[i]\leq5\) if odd and \(\root[i]\leq2\) if even.
        \item[(\(d=2\))] \(\root[i]=1\) or \(\root[i]\leq10\) if even.
        \item[(\(d=3\))] \(\root[i]\leq15\) if odd and \(\root[i]\leq6\) if even.
    \end{enumerate}
\end{lem}
\begin{proof}
    Consider double branching of the form 
    \begin{equation}\label{eqn:branch}
        \begin{tikzcd}
            \&\& \ceil{\tfrac{\root[i]+d}{2}} + x \ar[dl] \\
            \bal{\root[i]} \ar[r] \& \root[i] + d \\
            \&\& \ceil{\tfrac{\root[i]+d}{2}} + y \ar[ul]
        \end{tikzcd}
    \end{equation}
    then the balancing condition at \(j\) implies \(x+y=d\) if \(\root[i]+d\) is even (equivalently \(\root[i]\) and \(d\) have the same parity) and \(x+y=d-1\) if \(\root[i]+d\) is odd. Triple branching is even more restrictive so we deal with that on a case-by-case basis.
    \begin{enumerate}
        \item[(\(d=0\))] By \cref{lem:branch} the only possible neighbours of $\root[j]=\root[i]$ are either a vertex of weight $\root[i]$ or two vertices of weight $\tfrac{\root[i]}{2}$. The latter is only possible if $\root[i]$ is even and is the only possible balanced termination of such a chain.
        
        \item[(\(d=1\))] By \cref{lem:branch}, a triple branching is only possible when \(\root[i]+1\geq2\root[i]\) (\(\root[i]\) odd) or \(\root[i]+1\geq2\root[i]+3\) (\(\root[i]\) even), so only with \(\root[i]=1\) could we induce a triple branch. Now consider \eqref{eqn:branch}, with \(x+y=0\) if \(\root[i]\) is even and \(x+y=1\) if \(\root[i]\) is odd. In the former case, both new vertices have weight \(\tfrac{\root[i]+2}{2}\) (\(x=0, y=0\)), which forces \(\root[i]\leq2\) by \cref{lem:stepdown} with difference \(\tfrac{\root[i]}{2}\). In the latter case we may only take \(x=1, y=0\) and apply \cref{lem:stepdown} with difference \(\tfrac{\root[i]-1}{2}\), giving \(\root[i]\leq5\).
        
        \item[(\(d=2\))] By \cref{lem:branch}, a triple branching is only possible when \(\root[i]+2\geq2\root[i]\) (\(\root[i]\) even) or \(\root[i]+2\geq\root[i]+3\) (\(\root[i]\) odd), so only \(\root[i]=2\) could induce a triple branch. Now consider \eqref{eqn:branch}, with \(x+y=2\) if \(\root[i]\) is even and \(x+y=1\) if \(\root[i]\) is odd. If in the even case we have \(x=2, y=0\) then the top vertex is already balanced, but applying \cref{lem:stepdown} (with difference \(\tfrac{\root[i]-2}{2}\)) to the lower branch gives \(\root[i]+2\geq\tfrac{3\root[i]-6}{2}\), that is \(\root[i]\leq10\). If \(x=y=1\) then both branches are the same and we apply \cref{lem:stepdown} with difference \(\tfrac{\root[i]}{2}\) to see \(\root[i]\leq4\). In the odd case the only option is \(x=1, y=0\), and \cref{lem:stepdown} with difference \(\tfrac{\root[i]+3}{2}\) forces \(\root[i]\leq1\).
        
        \item[(\(d=3\))] By \cref{lem:branch}, a triple branching is only possible when \(\root[i]+3\geq2\root[i]\) (\(\root[i]\) odd) or \(\root[i]+3\geq2\root[i]+3\) (\(\root[i]\) even), so only \(\root[i]=1\) or \(\root[i]=3\) could induce a triple branch. Now consider \eqref{eqn:branch} and suppose \(\root[i]\) is even. Then \(x+y=2\) and we first check \((x,y)=(2,0)\). If the lower branch is balanced then by \cref{lem:stepdown} we must have \(\root[i]+3\geq\tfrac{3\root[i]+6}{2}\), that is, \(\root[i]\leq0\). The other option is \((x,y)=(1,1)\), for which both branches are the same. Balancing and \cref{lem:stepdown} again enforces \(\root[i]+3\geq\tfrac{3\root[i]}{2}\), thus \(\root[i]\leq6\). If \(\root[i]\) is odd then \(\root[i]+3\) is even and \(x+y=3\). We first check \((x,y)=(3,0)\). The lower branch is already balanced, and for the upper branch we use \cref{lem:stepdown} to imply \(\root[i]+3\geq\tfrac{3\root[i]-9}{2}\), so \(\root[i]\leq15\). If instead \((x,y)=(2,1)\), then the lower branch implies \(\root[i]+3\geq\tfrac{3\root[i]+3}{2}\) by \cref{lem:stepdown}. Hence \(\root[i]\leq3\). \qedhere
     \end{enumerate}
\end{proof}

Next we look at possible terminations.

\begin{lem}\label{lem:1chain}
    A finite balanced chain of vertices of weight \(1\) has to either be some \(\mir(\A)\) or have unbalanced vertices at each end.
\end{lem}
\begin{proof}
    Consider two adjacent weight \(1\) vertices. If both are balanced then they can only have neighbours of weight $1$. Repeating this procedure, we either stop at an unbalanced $1$ or form a cycle.
\end{proof}

\begin{lem}\label{lem:12chain}
    A finite chain of successive weights starting from \(1\), where every weight after \(1\) is balanced, has to terminate either as some \(\mir(\D)\) or \(\mir(\E)\).
\end{lem}
\begin{proof}
    Let \(\root\) be a dimension vector with fixed unbalanced vertices that are connected in some way. Suppose that \(\root[i]=1\) has a neighbour \(\root[i+1]=2\) which needs to be balanced. If \(\root\) terminates without cycles that containing \(i+1\), then taking this balanced chain along with \(i\) gives a completely balanced dimension vector (if considered separately from \(\root\)). We know, however, by \cref{thm:classification_of_dim2} that these are given by multiples of \(\mir(\ade)\). However, no multiple of $\mir(\A)$ can contain vertices of weight $1$ and $2$.
\end{proof}

  \section{The dimension vectors}\label{sec:roots}
  The aim of this section is to prove \cref{thm:dimvec}, that is we classify all\(\root\in\Sig\) with \(\p(\root)=2\).

\subsection{Types of dimension vectors for 4d quiver varieties}\label{ssec:type}

We first use \cref{lem:fset} and find the \(\root\in F\) with total balance \(-2\). This restricts the behaviour of our dimension vectors to three cases.

\begin{lem}\label{lem:type}
    Let \(\root\in\Sig\) satisfy \(\p(\root)=2\). Then \(\root\) has at most two unbalanced vertices, namely:
    \begin{enumerate}
        \item[(\(\I\))] \(i\neq j\in\graph[0]\) such that \(\root[i]=\root[j]=1\) and \((\root, \srt[i]) = (\root, \srt[j]) = -1\),
        \item[(\(\II\))] \(i\in\graph[0]\) such that \(\root[i]=1\) and \((\root, \srt[i])=-2\),
        \item[(\(\III\))] \(i\in\graph[0]\) such that \(\root[i]=2\) and \((\root, \srt[i])=-1\).
    \end{enumerate}
\end{lem}
\begin{proof}
    As \(\root\in\Sig\), by \cref{lem:fset} we necessarily have \((\root,\srt[i])\leq0\) for all \(i\in\graph[0]\). Then \(\p(\root)=2\) is equivalent to
    \begin{equation*}
        -2=(\root,\root)=\sum_{i\in\graph[0]}\root[i](\root, \srt[i]).
    \end{equation*}
    As \(\root\geq*0\), it follows that each summand is a nonpositive integer. Hence there are at most two nonzero terms; either they are both \(-1\) (\(\I\)) or one is \(-2\) (\(\II\) or \(\III\)).
\end{proof}

Because \cref{lem:fset} is not sufficient, some of the dimension vectors obtained in this way may in fact lie outside \(\Sig\), a problem we deal with in \cref{ssec:sig}.

\subsection{Local unbalanced structure}\label{ssec:locstr}

Given \(\root\) taking one of the forms in \cref{lem:type}, we first determine its unbalanced vertices, which must then be furnished with balanced vertices to complete the dimension vector.

\begin{dfn}\label{dfn:locstr}
    The \emph{local structure} of a (possibly unbalanced) vertex \(i\) is the weighted subgraph \((\graph',\root')\subseteq(\graph,\root)\) consisting of \(i\) and neighbours \(i_1,\dots,i_n\), such that \((\root, \srt[i_j])|_{\graph'}\geq0\) for all \(j\). In other words, it is still possible for the neighbours of \(i\) to be balanced in \( (\graph,\root)\).
\end{dfn}

\begin{lem}\label{lem:locstr}
    Let \(\root\in\Sig\) satisfy \(\p(\root)=2\). Then any unbalanced vertex \(i\) has one of the following local structures; five of Type \(\I\), thirteen of Type \(\II\), and nine of Type \(\III\). \vspace{1.6em} 
    \begin{longtable}{|C|C|C|C|}\caption{Local structures for an unbalanced vertex when \(\p(\root)=2\)}\label{tbl:locstr}\\\hline\endfirsthead\caption{(continued)}\endhead
        \rowcolor{gray!20} \I<1> & \I<2> & \I<3> & \I<4> \\\hline
        \begin{tikzcd}
            \bal{1} \ar[r] \& 3
        \end{tikzcd} & \begin{tikzcd}
            1 \ar[r, shift left] \ar[r, shift right] \& \bal{1} \ar[r] \ar[loop, distance=2em, out=60, in=120, phantom] \ar[loop, distance=2em, out=240, in=300, phantom] \& 1
        \end{tikzcd} & \begin{tikzcd}
            2 \ar[r] \& \bal{1} \ar[r] \& 1
        \end{tikzcd} & \begin{tikzcd}
            \& 1 \ar[d] \\
            1 \ar[r] \& \bal{1} \ar[r] \& 1
        \end{tikzcd} \\\hline
        \rowcolor{gray!50} \cellcolor{gray!20} \I<5> & \II<1> & \II<2> & \II<3> \\\hline
        \begin{tikzcd}
            \bal{1} \ar[loop, distance=2em, out=150, in=210] \ar[r] \& 1
        \end{tikzcd} & \begin{tikzcd}
            \bal{1} \ar[r, shift left] \ar[r, shift right] \& 2
        \end{tikzcd} & \begin{tikzcd}
            \bal{1} \ar[r] \& 4
        \end{tikzcd} & \begin{tikzcd}
            1 \ar[r, shift left] \ar[r, shift right] \& \bal{1} \ar[r, shift left] \ar[r, shift right] \ar[loop, distance=2em, out=60, in=120, phantom] \ar[loop, distance=2em, out=240, in=300, phantom] \& 1
        \end{tikzcd} \\\hline
        \rowcolor{gray!50} \II<4> & \II<5> & \II<6> & \II<7> \\\hline
        \begin{tikzcd}
            1 \ar[r, shift left] \ar[r, shift right] \& \bal{1} \ar[r] \& 2
        \end{tikzcd} & \begin{tikzcd}
            2 \ar[r] \& \bal{1} \ar[r] \& 2
        \end{tikzcd} & \begin{tikzcd}
            1 \ar[r] \& \bal{1} \ar[r] \& 3
        \end{tikzcd} & \begin{tikzcd}
            \& 1 \ar[d, shift left] \ar[d, shift right] \\
            1 \ar[r] \& \bal{1} \ar[r] \& 1
        \end{tikzcd} \\\hline
        \rowcolor{gray!50} \II<8> & \II<9> & \II<10> & \II<11> \\\hline
        \begin{tikzcd}
            \& 2 \ar[d] \\
            1 \ar[r] \& \bal{1} \ar[r] \& 1
        \end{tikzcd} & \begin{tikzcd}
            1 \ar[dr] \&\& 1 \ar[dl] \\
            1 \ar[r] \& \bal{1} \ar[r] \& 1
        \end{tikzcd} & \begin{tikzcd}
            \bal{1} \ar[loop, distance=2em, out=150, in=210] \ar[r, shift left] \ar[r, shift right] \& 1
        \end{tikzcd} & \begin{tikzcd}
            \bal{1} \ar[loop, distance=2em, out=150, in=210] \ar[r] \& 2
        \end{tikzcd} \\\hline
        \rowcolor{gray!50} \II<12> & \II<13> & \cellcolor{gray!20} \III<1> & \cellcolor{gray!20} \III<2> \\\hline
        \begin{tikzcd}
            1 \ar[r] \& \bal{1} \ar[loop, distance=2em, out=60, in=120] \ar[loop, distance=2em, out=240, in=300, phantom] \ar[r] \& 1
        \end{tikzcd} & \begin{tikzcd}
            \bal{1} \ar[loop, distance=2em, out=150, in=210] \ar[loop, distance=2em, out=-30, in=30]
        \end{tikzcd} & \begin{tikzcd}
            \bal{2} \ar[r] \& 5
        \end{tikzcd} & \begin{tikzcd}
            1 \ar[r] \& \bal{2} \ar[r, shift left] \ar[r, shift right] \& 2
        \end{tikzcd}
        \\\hline
        \rowcolor{gray!20} \III<3> & \III<4> & \III<5> & \III<6> \\\hline
         \begin{tikzcd}
            2 \ar[r] \& \bal{2} \ar[loop, distance=2em, out=60, in=120, phantom] \ar[loop, distance=2em, out=240, in=300, phantom] \ar[r] \& 3
        \end{tikzcd} & \begin{tikzcd}
            1 \ar[r] \& \bal{2} \ar[r] \& 4
        \end{tikzcd} & \begin{tikzcd}
            \& 1 \ar[d] \\ 
            2 \ar[r] \& \bal{2} \ar[r] \& 2
        \end{tikzcd} & \begin{tikzcd}
            \& 3 \ar[d] \\
            1 \ar[r] \& \bal{2} \ar[r] \& 1
        \end{tikzcd} \\\hline
        \rowcolor{gray!20} \III<7> & \III<8> & \III<9> \\\cline{1-3}
        \begin{tikzcd}
            1 \ar[dr] \&\& 1 \ar[dl] \\
            2 \ar[r] \& \bal{2} \ar[r] \& 1
        \end{tikzcd} & \begin{tikzcd}
            1 \ar[dr] \& 1 \ar[d] \& 1 \ar[dl] \\
            1 \ar[r] \& \bal{2} \ar[r] \& 1
        \end{tikzcd} & \begin{tikzcd}
            \bal{2} \ar[loop, distance=2em, out=150, in=210] \ar[r] \& 1    
        \end{tikzcd} \\\cline{1-3}
    \end{longtable}
\end{lem}
\begin{proof}
    As \(\root\in\Sig\), we may apply \cref{lem:type}. Thus in each case the weights neighbouring \(i\) have sum
    \begin{enumerate}
        \item[(\(\I\))] \(3\), or \(1\) if \(\ell_i=1\);
        \item[(\(\II\))] \(4\), or \(2\) if \(\ell_i=1\), or \(0\) if \(\ell_i=2\);
        \item[(\(\III\))] \(5\), or \(1\) if \(\ell_i=1\).
    \end{enumerate}
    Recall that we view a loop at \(i\) as the vertex having two extra neighbours of weight \(\root[i]\). To make sure nothing is missed, we order (increasingly) the possibilities first by the number of loops, then the number of distinct neighbours, then the maximal neighbouring weight present, and then the minimal neighbouring weight present. If \(\root[i]=2\) then \cref{lem:half} forces at most one edge when connecting \(i\) to weight \(1\) vertices, so we can omit most cases of double edges and higher multiplicities. Similarly for \(\root[i]=1\) we have at most two edges shared with any neighbour. Hence we see that only the tabulated weighted subgraphs are possible.
\end{proof}

There is one dimension vector \(\root\) of Type \(\I\) that doesn't come from a local structure as in \cref{dfn:locstr}, as it consists only of unbalanced vertices:
\begin{equation*}
    \begin{tikzcd}
        \bal{1} \ar[r, shift left] \ar[r] \ar[r, shift right] \& \bal{1}
    \end{tikzcd}
    \quad=\I(1,1,1).
\end{equation*}

\subsection{The classification}\label{ssec:class}

We now proceed by type to complete the local structures in \cref{tbl:locstr} with balanced vertices, using the restrictions in \cref{ssec:balancing}.

\begin{ppn}\label{ppn:i}
    The dimension vectors of Type \(\I\) are \(\I(\ade<1>,\ade<2>,n)\), \(\I(\ade, n)\), \(\I(m)\), \(\I[a]\), \(\I[b]\), and \(\I[c]\).
\end{ppn}
\begin{proof}
    For \(M,N\in\set{1,2,3,4,5}\), we attempt to connect \(\I<M>\) with \(\I<N>\) using only balanced vertices. Observe that \(\I<1>\) can only pair up with itself, as all other local structures have an odd number of balanced vertices with weight \(1\), and so one would need to connect to the weight \(3\) vertex, contradicting \cref{lem:1chain}.
    \begin{itemize}
        \item[(1),(1)] Firstly if we look at the local structure of one of the unbalanced vertices. We have the following possible sets of neighbours to balance $3$:
        \begin{equation*}
            (1,2,2),\quad (1,4),\quad (2,3),\quad (5).
        \end{equation*}
        The first two determine the unbalanced vertices. The first yields 
        \begin{equation*}
            \begin{tikzcd}[scale cd=.8, cramped, sep=0]
                \&\&\bal{1} \\ 1\&2\&3\&2\&1 \\ \&\&\bal{1}
            \end{tikzcd}
            \quad=\I(4),
        \end{equation*}
        as the decreasing chains are forced by \cref{lem:stepdown}. Now take a vertex of weight \(4\). It requires total weight \(5\) to be balanced, we may balance with a \(2\) and a \(3\), giving. 
        \begin{equation*}
            \begin{tikzcd}[scale cd=.8, cramped, sep=0]
                \bal{1}\&\&2 \\ \&3\&4\&3\&2\&1 \\ \bal{1}
            \end{tikzcd}
            \quad=\I[a].
        \end{equation*}
        Otherwise continue the chain with a \(5\), which must then have neighbours at least \(3\) that sum to \(6\). A pair of \(3\)s contradicts \cref{lem:stepdown} as \(2\) does not divide \(5\). Continuing with a \(6\) is our last step by \cref{lem:stepup} (\(d=1\)), so the remaining terminalisation is
        \begin{equation*}
            \begin{tikzcd}[scale cd=.8, cramped, sep=0]
                \bal{1}\&\&\&\&3 \\ \&3\&4\&5\&6\&4\&2 \\ \bal{1}
            \end{tikzcd}
            \quad=\I[c].
        \end{equation*}
        Now suppose the weight \(3\) vertices are not identified, to each we must attach weights at least \(2\) that total \(5\). If on one copy of \(\I<1>\) we attach \(2\) and \(3\), the former weight starts a decreasing chain and the latter weight starts a constant balanced chain of \(3\)s. We can only terminate this chain by adding an unbalanced vertex to one of the weight \(3\) vertices, that is, by symmetrically attaching the other copy of \(\I<1>\) to obtain
        \begin{equation*}
            \begin{tikzcd}[scale cd=.8, cramped, sep=0]
                \&\&\bal{1}\&\&\bal{1} \\ 1\&2\&3\&\cdots\&3\&2\&1
            \end{tikzcd}
            \quad=\I(m)
        \end{equation*}
        with \(m\geq5\). Now continue one copy of \(\I<1>\) with a weight \(5\) vertex, so that neighbours must be at least \(3\) and total \(7\). There are two options to proceed. One is to take a weight $7$ vertex and proceed. This creates a growing chain of consecutive odd weights. By \cref{lem:stepup}(2) we know this cannot be terminated in a balanced way. Moreover, from the proof of \cref{lem:stepup}(2) one can see that there is no other way to finish it in an unbalanced way with local structure \(\I<1>\).
        We thus proceed with the only remaining option which is to take \(4\) and \(3\). The former is a forced terminating downwards chain and the latter necessarily being the other copy of \(\I<1>\):
        \begin{equation*}
            \begin{tikzcd}[scale cd=.8, cramped, sep=0]
                \bal{1}\&3 \\ \&\&5\&4\&3\&2\&1 \\ \bal{1}\&3
            \end{tikzcd}
            \quad=\I[b].
        \end{equation*}
        \item[(2),(2)] By \cref{lem:1chain}, the only possibility is (for some \(n\geq1\))
        \begin{equation*}
            \begin{tikzcd}[scale cd=.8, cramped, sep=0]
                1=\bal{1}\&\cdots\&\bal{1}=1
            \end{tikzcd}
            \quad=\I(\A[1],\A[1],n).
        \end{equation*}
        \item[(2),(3)] Similarly to the previous argument, vertices of weight \(1\) are connected in a chain of \(1\)s. On the other side, we apply \cref{lem:12chain} to see this must terminate as \(\mir(\ade)\) for some \(\ade\neq\A\), giving
        \begin{equation*}
            \begin{tikzcd}[scale cd=.8, cramped, sep=0]
                \hrr(\ade)\&\bal{1}\&\cdots\&\bal{1}=1
            \end{tikzcd}
            \quad=\I(\ade,\A[1],n).
        \end{equation*}
        \item[(2),(4)] Applying \cref{lem:1chain} again, the only possibility is to have a chain of \(1\)s connecting the unbalanced vertices and to connect another pair of \(1\)s in a cycle. This is
        \begin{equation*}
            \begin{tikzcd}[scale cd=.8, cramped, sep=0]
                \&1 \\\\ 1\&\cdots\&\bal{1}\&\cdots\&\bal{1}=1
            \end{tikzcd}
            \quad=\I(\A[m],\A[1],n).
        \end{equation*}
        \item[(2),(5)] This is another use of \cref{lem:1chain}, so for some \(n\geq1\) we get
        \begin{equation*}
            \begin{tikzcd}[scale cd=.8, cramped, sep=0]
                \bal{1}\ar[loop above, looseness=4]\&\cdots\&\bal{1}=1
            \end{tikzcd}
            \quad=\I(\A[0],\A[1],n).
        \end{equation*}
        \item[(3),(3)] \cref{lem:12chain} tells us that the only balanced termination of
            \begin{tikzcd}
                \bal{1} \ar[r] \& 2 \ar[r] \& \cdots
            \end{tikzcd} is with some \(\mir(\D)\) or \(\mir(\E)\). Suppose we start with one of our unbalanced vertices. However, notice that one of the weight \(1\) vertices could be unbalanced. If that did not happen this means that each of the unbalanced vertices has its own \(\ade\) root, that is
        \begin{equation*}
            \begin{tikzcd}[scale cd=.8, cramped, sep=0]\hrr(\ade<1>)\&\bal{1}\&\cdots\&\bal{1}\&\hrr(\ade<2>)
            \end{tikzcd}
            \quad=\I(\ade<1>,\ade<2>,n),
        \end{equation*}
        for some \(n\geq1\) and \(\ade<1>,\ade<2>\) not of type \(\A\). If that did happen, however, we get one of the following:
        \begin{equation*}
            \begin{tikzcd}[scale cd=.8, cramped, sep=0]
                1\&\&\&\&1 \\ \&2\&\cdots\&2 \\ \bal{1}\ar[rrrr, bend right, "\cdots" description]\&\&\&\&\bal{1} 
            \end{tikzcd}, \qquad\qquad
            \begin{tikzcd}[scale cd=.8, cramped, sep=0]
                1\&\&\&\&\bal{1}\ar[dd, bend left = 5em, "\cdots" description, sloped] \\ \&2\&\cdots\&2 \\ 1\&\&\&\&\bal{1}
            \end{tikzcd}, \qquad\qquad
            \begin{tikzcd}[scale cd=.8, cramped, sep=0]
                \&\&1 \\ \&\&2 \\ \bal{1} \ar[rrrr, bend right, "\cdots" description] \&2\&3\&2\&\bal{1}
            \end{tikzcd}, \qquad\qquad
            \begin{tikzcd}[scale cd=.8, cramped, sep=0] 
                \&\&\&2 \\ \bal{1}\ar[rrrrrr, bend right, "\cdots" description] \&2\&3\&4\&3\&2\&\bal{1}
            \end{tikzcd}
            \quad=\I(\ade,n).
        \end{equation*}
        \item[(3),(4)] Combining Lemmas \ref{lem:1chain} and \ref{lem:12chain}, we see the only possibility is
        \begin{equation*}
            \begin{tikzcd}[scale cd=.8, cramped, sep=0]
                \&\&\&\& 1 \\ \hrr(\ade)\&\bal{1}\&\cdots\&\bal{1}\&\cdots\&1
            \end{tikzcd}
            \quad=\I(\ade,\A[m],n),
        \end{equation*}
        for some \(n\geq1\) and \(\ade\neq\A\).
        \item[(3),(5)] As in the above argument, this combination forces
        \begin{equation*}
            \begin{tikzcd}[scale cd=.8, cramped, sep=0]
                \hrr(\ade)\&\bal{1}\&\cdots\&\bal{1}\ar[loop above, looseness=4]
            \end{tikzcd}
            \quad=\I(\ade,\A[0],n)
        \end{equation*}
        with the same restrictions on \(n\) and \(\ade\).
        \item[(4),(4)] \cref{lem:1chain} tells us that we have to connect chains of \(1\)s, and there are two ways to do this. If we pair up the hanging \(1\)s between the copies of \(\I<4>\), we obtain \(\I(\ell,m.n)\), otherwise we have
        \begin{equation*}
            \begin{tikzcd}[scale cd=.8, cramped, sep=0]
                \&1\&\&\&\&1 \\ 1\&\cdots\&\bal{1}\&\cdots\&\bal{1}\&\cdots\&1
            \end{tikzcd}
            \quad=\I(\A[\ell],\A[m],n)
        \end{equation*}
        for some \(\ell,m\geq2\) and \(n\geq1\).
        \item[(4),(5)] This is the previous case but with \(m=0\), that is
        \begin{equation*}
            \begin{tikzcd}[scale cd=.8, cramped, sep=0]
                \&1 \\ 1\&\cdots\&\bal{1}\&\cdots\&\bal{1}\ar[loop above, looseness=4]
            \end{tikzcd}
            \quad=\I(\A[\ell],\A[0],n).
        \end{equation*}
        \item[(5),(5)] This is the limiting case of the previous two arguments, and we get
        \begin{equation*}
            \begin{tikzcd}[scale cd=.8, cramped, sep=0]\bal{1}\ar[loop above, looseness=4]\&\cdots\&\bal{1}\ar[loop above, looseness=4]
            \end{tikzcd}
            \quad=\I(\A[0],\A[0],n). \qedhere
        \end{equation*}
    \end{itemize}
\end{proof}

Having finished the most complicated type, we move on to those that have at most one unbalanced vertex.

\begin{ppn}\label{ppn:ii}
    The dimension vectors of Type \(\II\) are \(\II(\ade<1>,\ade<2>)\), \(\II(\ade)\) for \(\ade\neq\E[8]\), \(\II(m)\), \(\II[a]\), and \(\II[b]\).
\end{ppn}
\begin{proof}
    As before, we proceed by completing local structure in \cref{tbl:locstr} with balanced vertices.
    \begin{enumerate}
        \item  By \cref{lem:stepup}, the hanging weight \(2\) can either be completed by a pair of \(1\)s or a \(2\). The former terminates the process and the latter repeats this same choice. Hence for \(m\geq4\) this is
        \begin{equation*}
            \begin{tikzcd}[scale cd=.8, cramped, sep=0]
                \&\&\&1 \\ \bal{1}=2\&\cdots\&2 \\ \&\&\&1
            \end{tikzcd}
            \quad=\II(\D[m])'.
        \end{equation*}
        \item The neighbours of \(4\) must be each at least \(2\) and must total \(2\times4-1=7\). The only options that do not lead to a contradiction are \((3,4)\), \((2,2,3)\), \((2,5)\), and \((7)\). Choosing \((2,5)\) forces
        \begin{equation*}
            \begin{tikzcd}[scale cd=.8, cramped, sep=0]
                \&2\&\&3 \\ \bal{1}\&4\&5\&6\&4\&2
            \end{tikzcd}
            \quad=\II[a].
        \end{equation*}
        by the proof of \cref{lem:12chain}. Choosing $(7)$ means that the corresponding vertex needs to have neighbours that add up to $10$. There are $3$ possibilities: \((10)\), \((4,6)\), and \((5,5)\). By \cref{lem:stepdown} the last two are impossible to terminate in a balanced way. Proceeding, for \((10)\) the only possible neighbours are \((5,8)\), \((6,7)\), and \((13)\). The first one is 
        \begin{equation*}
            \begin{tikzcd}[scale cd=.8, cramped, sep=0]
                \&\&\&5 \\ \bal{1}\&4\&7\&10\&8\&6\&4\&2
            \end{tikzcd}
            \quad=\II[b].
        \end{equation*}
        The other two options are impossible to terminate by \cref{lem:stepdown} and \cref{lem:stepup}(3), respectively. Finally, choosing \((3,4)\) (with \((2,2,3)\) as a special case) has a forced termination from the vertex of weight \(3\) and, by \cref{lem:stepup} the chain of \(4\)s terminates as
        \begin{equation*}
            \begin{tikzcd}[scale cd=.8, cramped, sep=0]
                2\&\&\&\bal{1} \\ \&4\&\cdots\&4\&3\&2\&1 \\ 2
            \end{tikzcd}
            \quad=\II(m).
        \end{equation*}
        \item This is already complete, and is \(\II(\A[1],\A[1])=\I(\A[1],\A[1],0)\).
        \item By \cref{lem:12chain} we obtain
        \begin{equation*}
            \begin{tikzcd}[scale cd=.8, cramped, sep=0]
                \ade\&1=\bal{1}
            \end{tikzcd}
            \quad=\II(\ade,\A[1]).
        \end{equation*}
        \item Completing the weight \(2\) vertices separately gives rise to \begin{equation*}
            \begin{tikzcd}[scale cd=.8, cramped, sep=0]
                \ade<1>\&\bal{1}\&\ade<2>
            \end{tikzcd}
            \quad=\II(\ade<1>,\ade<2>).
        \end{equation*}
        for \(\ade\neq\A\) by \cref{lem:12chain}. We could also terminate by linking the weight \(2\) vertices in some way, but this is where the results in \cref{ssec:balancing} cannot help us. We proceed from one of the ends and go case by case. The possible neighbours to $2$ (that don't force a dimension vector without cycles) are \((1,2)\) and \((3)\). The first one gives rise to a chain of $2$ which can only be connected to another $2$. However, as the dimension vector has cycles it has to connect to the initial \(2\) from \(\I<5>\), resulting in
        \begin{equation*}
            \begin{tikzcd}[scale cd=.8, cramped, sep=0]
                1\&\&\&\&1 \\ \&2\&\cdots\&2 \\ \&\&\bal{1}
            \end{tikzcd}
            \quad=\II(\D[m]).
        \end{equation*}
        Now suppose the next weight is \(3\), which can have neighbours \((2,2)\) or \((4)\). The first choice results in a dimension vector without cycles unless one of these $2$s was in \(\I<5>\), that is,
        \begin{equation*}
            \begin{tikzcd}[scale cd=.8, cramped, sep=0]
                \&1 \\ \&2 \\ 2\&3\&2 \\\\ \&\bal{1} \ar[uul, bend left] \ar[uur, bend right]
            \end{tikzcd}
            \quad=\II(\E[6]).
        \end{equation*}
        Now suppose the next weight is \(4\). The possible neighbours are $(3,2)$ or $(5)$. The first terminates the root in a down-by-one chain without cycles unless $2$ was in \(\I<5>\). A choice of $(5)$ forces $6$ and the only possible branching by \cref{lem:stepup}(1) is $\E[8]$, which does not have a vertex of weight $1$ to be unbalanced. Thus the only remaining option is 
        \begin{equation*}
            \begin{tikzcd}[scale cd=.8, cramped, sep=0]
                \&\&2 \\ 2\&3\&4\&3\&2 \\\\\\ \&\& \bal{1}\ar[uuull, bend left] \ar[uuurr, bend right]
            \end{tikzcd}
            \quad=\II(\E[7]).
        \end{equation*}
        \item Since the hanging \(1\) does not have another to connect with, this local structure cannot be completed in a balanced way.
        \item By \cref{lem:1chain}, pairing the hanging \(1\)s with an arbitrary chain completes the diagram, for \(m\geq2\), as
        \begin{equation*}
            \begin{tikzcd}[scale cd=.8, cramped, sep=0]
                \&1 \\ 1=\bal{1}\&\cdots\&1
            \end{tikzcd}
            \quad=\II(\A[1],\A[m]).
        \end{equation*}
        \item Combining Lemmas \ref{lem:1chain} and \ref{lem:12chain} we obtain
        \begin{equation*}
            \begin{tikzcd}[scale cd=.8, cramped, sep=0]
                \&\&1 \\ \hrr(\ade)\&\bal{1}\&\cdots\&1
            \end{tikzcd}
            \quad=\II(\ade,\A[m])
        \end{equation*}
        with \(\ade\neq\A\) and \(m\geq2\).
        \item As before, the only option is to pair up the \(1\)s with arbitrary chains, that is, for \(\ell,m\geq2\),
        \begin{equation*}
            \begin{tikzcd}[scale cd=.8, cramped, sep=0]
                \&1\&\&1 \\ 1\&\cdots\&\bal{1}\&\cdots\&1
            \end{tikzcd}
            \quad=\II(\A[\ell],\A[m]).
        \end{equation*}
        \item This is already complete, and is \(\II(\A[0],\A[1])\).
        \item By \cref{lem:12chain}, the hanging weight \(2\) can only be completed through some \(\mir(\D)\) or \(\mir(\E)\):
        \begin{equation*}
            \begin{tikzcd}[scale cd=.8, cramped, sep=0]
                \hrr(\ade)\&\bal{1}\ar[loop above, looseness=4]
            \end{tikzcd}
            \quad=\II(\ade,\A[0]).
        \end{equation*}
        \item Again, by \cref{lem:1chain} the only way to complete the hanging weight \(1\) vertices is to connect them with an arbitrary weight \(1\) chain (for \(m\geq2\)).
        \begin{equation*}
            \begin{tikzcd}[scale cd=.8, cramped, sep=0]
                \&1 \\ \bal{1}\ar[loop above, looseness=4]\&\cdots\&1
            \end{tikzcd}
            \quad=\II(\A[0],\A[m]).
        \end{equation*}
        \item This is already complete, and is \(\II(\A[0],\A[0])\). \qedhere
	\end{enumerate}
\end{proof}

Finally we consider which dimension vectors arise from terminating the last nine local structures.

\begin{ppn}\label{ppn:iii}
    The dimension vectors of Type \(\III\) are \(\III(\ade)\), \(\III(2\ade)\), \(\III(m,n)\), \(\III[a]\), and \(\III[b]\).    
\end{ppn}
\begin{proof}
    We proceed in entirely the same way as the previous two results.
  \begin{enumerate}
    \item The weight \(5\) vertex must have neighbours each at least \(3\) that total \(8\), so the options are two \((4,4)\) or \((8)\). The former has a forced descending termination and the latter can only (by \cref{lem:stepdown}) continue by attaching \(4\) and \(7\). These are
    \begin{equation*}
        \begin{tikzcd}[scale cd=.8, cramped, sep=0]
            \&\&\&\&\bal{2} \\ 1\&2\&3\&4\&5\&4\&3\&2\&1
        \end{tikzcd}\quad=\III[a], \qquad \begin{tikzcd}[scale cd=.8, cramped, sep=0]
            \&\&4 \\ \bal{2}\&5\&8\&7\&6\&5\&4\&3\&2\&1            
        \end{tikzcd}
        \quad=\III[b].
    \end{equation*}
    \item This diagram is already complete, and is \(\III(2\A[1])\).
    \item We must complete the vertex of weight \(2\) and the vertex of weight \(3\). By \cref{lem:stepup}(0), the former only terminates with balanced vertices through a chain of $2$ similar to \(\D[m]\). By the proof of \cref{lem:12chain}, the latter terminates with a choice of \(\mir(\E[n])\), yielding
    \begin{equation*}
        \begin{tikzcd}[scale cd=.8, cramped, sep=0]
            1 \\ \&2\&\cdots\&\bal{2}\&3\&\cdots \\ 1           
        \end{tikzcd}
        \quad=\III(m,n).
    \end{equation*}
    \item We need to complete the vertex of weight $4$, with possible neighbours totalling \(6\). Thus the possibilities are \((3,3)\), \((2,4)\), or \((6)\). By \cref{lem:stepdown}, the first option forces
    \begin{equation*}
        \begin{tikzcd}[scale cd=.8, cramped, sep=0]
            \&\&\&1 \\ \&\&\&\bal{2} \\ 1\&2\&3\&4\&3\&2\&1            
        \end{tikzcd}
        \quad=\III(\E[7]).
    \end{equation*}
    Further, by \cref{lem:stepup}, $(2,4)$ can only end as a chain of $4$s that splits at the end, that is,
    \begin{equation*}
        \begin{tikzcd}[scale cd=.8, cramped, sep=0]
            1\&\bal{2}\&\&\&\&2 \\ \&\&4\&\cdots\&4 \\ \&2\&\&\&\&2
        \end{tikzcd}
        \qquad=\III(2\D[m]).
    \end{equation*}
    Next, suppose next vertex is $6$. Its possible neighbours are $(4,4), (3,5)$ or $(8)$. The former yields
    \begin{equation*}
        \begin{tikzcd}[scale cd=.8, cramped, sep=0]
            \&\&\&2 \\ \&\&\&4 \\ 1\&\bal{2}\&4\&6\&4\&2
        \end{tikzcd}
        \quad=\III(2\E[6]).
    \end{equation*}
    Applying \cref{lem:stepdown} to $(3,5)$, gives rise to
    \begin{equation*}
        \begin{tikzcd}[scale cd=.8, cramped, sep=0]
            \&\&\&\&\&3 \\ 1\&2\&3\&4\&5\&6\&4\&\bal{2}\&1
        \end{tikzcd}
        \quad=\III(\E[8]).
    \end{equation*}
    Lastly, $8$ needs to have neighbours of weight $10$. By \cref{lem:stepdown} the only options are $(4,6)$ (as \(3\) does not divide \(8\)), which terminates in
    \begin{equation*}
        \begin{tikzcd}[scale cd=.8, cramped, sep=0]
            \&\&\&\& 4 \\ 1\&\bal{2}\&4\&6\&8\&6\&4\&2
        \end{tikzcd}
        \quad=\III(2\E[7])
    \end{equation*}
    and $10$, which, by \cref{lem:stepdown} continues in $12$. Here, by \cref{lem:stepup} we have to branch as further going up is impossible. The only branching is $(6,8)$ which gives
    \begin{equation*}
        \begin{tikzcd}[scale cd=.8, cramped, sep=0]
            \&\&\&\&\&\&6 \\ 1\&\bal{2}\&4\&6\&8\&10\&12\&8\&4
        \end{tikzcd}
        \quad=\III(2\D[m]).
    \end{equation*}
    \item Both of the remaining vertices with weight \(2\) must be made balanced. By \cref{lem:stepup} we may terminate with a pair of weight \(1\) vertices or we may also connected them with an arbitrary chain of weight \(2\) vertices. The first option results in
    \begin{equation*}
        \begin{tikzcd}[scale cd=.8, cramped, sep=0]
            1\&\&\&1\&\&\&1 \\ \&2\&\cdots\&\bal{2}\&\cdots\&2 \\ 1\&\&\&\&\&\&1
        \end{tikzcd}
        \quad=\III(\D[m],i)
    \end{equation*}
    for \(i\neq1\). The second option, for \(m\geq2\), is
    \begin{equation*}
        \begin{tikzcd}[scale cd=.8, cramped, sep=0]
            \&2 \\ 2\&\cdots\&\bal{2}\&1
        \end{tikzcd}
        \quad=\III(2\A[m]).
    \end{equation*}
    \item By \cref{lem:12chain} the only possible balanced terminations of the weight \(3\) vertex chain are
    \begin{equation*}
        \begin{tikzcd}[scale cd=.8, cramped, sep=0]
            \&\&1 \\ 1\&\&2 \\ \&\bal{2}\&3\&2\&1 \\ 1           
        \end{tikzcd}\quad=\III(4,6), \qquad 
        \begin{tikzcd}[scale cd=.8, cramped, sep=0]
            1\&\&\&2 \\ \&\bal{2}\&3\&4\&3\&2\&1 \\ 1
        \end{tikzcd}\quad=\III(4,7), \qquad 
        \begin{tikzcd}[scale cd=.8, cramped, sep=0]
            1\&\&\&\&\&3 \\ \&\bal{2}\&3\&4\&5\&6\&4\&2 \\ 1
        \end{tikzcd}\quad=\III(4,8).
    \end{equation*}
    \item The only vertex to be made balanced is the remaining \(2\). By \cref{lem:stepup} we obtain an arbitrary chain of weight \(2\) which terminates by splitting into two vertices of weight $1$:
    \begin{equation*}
        \begin{tikzcd}[scale cd=.8, cramped, sep=0]
            1\&\&\&\&1 \\ 1\&\bal{2}\&\cdots\&2 \\ 1\&\&\&\&1
        \end{tikzcd}
        \quad=\III(\D[m],1).
    \end{equation*}
    \item This diagram is already complete, and is \(\III(\D[4])\).
    \item This diagram is already complete, and is \(\III(2\A[0])\). \qedhere
  \end{enumerate}
\end{proof}

This concludes the list of dimension vectors \(\root\) with nonpositive local balance and total balance \(-2\), comprising \cref{tbl:dimvec} along with \(\I(\ade<1>,\ade<2>,n)\) with positive parameter and \(\III(2\ade)\).

\subsection{Filtering for \(\Sig\)}\label{ssec:sig}

We now finish the proof of \cref{thm:dimvec} by removing the \(\root\not\in\Sig\) that satisfy \(\p(\root)=2\), as they will not give fundamentally (in the sense of the canonical decomposition) new quiver varieties of this dimension. There are two dimension vectors that we can immediately exclude.

\begin{lem}\label{lem:sigfails}
    The dimension vectors \(\I(\ade<1>,\ade<2>,n)\) and \(\III(2\ade)\) lie outside \(\Sig\).
\end{lem}
\begin{proof}
    Simply observe that \(\I(\ade<1>,\ade<2>,n)=\mir(\ade<1>)+\mir(\ade<2>)+\hrr(\A[n-1])\) (as \(n\geq1\), this decomposition exists) and \(\III(2\ade)=2\mir(\ade)+\srt[i]\) for some vertex \(i\). Applying \(\p\) to every term shows that neither dimension vector are in \(\Sig\).
\end{proof}

Apart from the above families, we claim that all of the dimension vectors \(\root\in\Root\) with \(\p(\root)=2\) lie in \(\Sig\), giving those seen in \cref{tbl:dimvec}. We can check this using the following result.

\begin{lem}\label{lem:fnotsig}
    If \(\root\in F\setminus\Sig\) then one of the following occurs:
    \begin{enumerate}
        \item There is some \(n\geq2\) such that \(\root=n\mir(\ade)\).
        \item There are two connected components of \(\root\) that are connected to each other by \begin{tikzcd}
            1\ar[r]\&1.
        \end{tikzcd}
        \item There are two connected components of \(\root\) that are connected to each other by \begin{tikzcd}
            1\ar[r]\&n,
        \end{tikzcd}
        where \(n\geq2\) is the extending weight of \(n\mir(\ade)\).
    \end{enumerate}
\end{lem}
\begin{proof}
    This is \cite[Thm. 8.1]{crawley-boevey_geometry_2001} with \(\uplambda=0\) and in characteristic zero.
\end{proof}

The dimension vectors with \(\p(\root)=2\) of the form \cref{lem:fnotsig}(2) are precisely \(\I(\ade<1>,\ade<2>,n)\) for \(n\geq*0\). For (3) we have \(\III(2\ade)\), as expected from \cref{lem:sigfails}. This gives enough to finish the classification of dimension vectors.

\begin{proof}[Proof of \cref{thm:dimvec}]
    In \cref{ssec:class} we find all \(\root\in\Root\) satisfying \(\p(\root)=2\) and \((\root,\srt[i])\leq0\) for every \(i\in\graph[0]\). Now use \cref{lem:fnotsig}. Since (1) is \(\p(\root)=1\) by \cref{thm:classification_of_dim2}, none of our dimension vectors apply. For (2), the only dimension vector with a weighted subgraph of the form \(\hrr(\A[2])\) is \(\I(\ade,n)\) for \(n\geq*0\), but the two ends of this edge do not lie in disconnected components. Similarly, none of the remaining dimension vectors take form (3). Hence only \cref{tbl:dimvec} remains in \(\Sig\).
\end{proof}

  \section{Symplectic stratification and minimal degenerations}\label{sec:symplectic}
  Since different dimension vectors could correspond to the same quiver variety, we would now like to have some geometric invariants to distinguish them. In this section we explain how to compute the minimal degenerations between the symplectic leaves of quiver varieties, which serves as a useful invariant. We then summarise these degenerations in the 4d case using \cref{tbl:decomp}. This also allows us to compute the Namikawa Weyl groups of the quiver varieties, establishing \cref{thm:namikawa_weyl_groups_are_classified}.

\subsection{Symplectic singularities}

We begin by recalling the definition of conical symplectic singularity, which first appeared in \cite{beauville_symplectic_2000}.

\begin{dfn}\label{dfn:sympsing}
    We say that an affine algebraic variety \(\variety\) over \(\Comp\) has \emph{symplectic singularities} if: 
    \begin{itemize}
        \item \(\variety\) is a normal Poisson variety.
        \item There exists a smooth, dense open subset $U \subset X$ (with complement of codimension \(2\)) on which the Poisson structure comes from a symplectic form $\upomega$.
        \item There exists a resolution of singularities, that is a birational and projective morphism $Y\ra X$, such that the pullback of $\upomega$ to $Y$ has no poles. 
    \end{itemize}
    If, in addition to the above conditions, one has a $\Comp<\times>$-action on \(\variety\) which acts on $\upomega$ with some positive weight and contracts \(\variety\) to the unique fixed point, \(X\) is called a \emph{conical symplectic singularity}.
\end{dfn}

Symplectic singularities possess additional structure, which we explain in the following subsection.

\subsection{Symplectic leaves}\label{ssec:sympleaf}

Being Hamiltonian reductions, quiver varieties come equipped with a Poisson structure. Symplectic leaves are maximal connected smooth Poisson submanifolds where the Poisson bracket is non-degenerate. For a better exposition on symplectic leaves, we refer the reader to \cite{martino_stratifications_2008}. By \cite[Thm. 1.2]{bellamy_symplectic_2021}, quiver varieties are conical symplectic singularities, so we can use the following theorem (also look at \cite[\S 3.v]{nakajima_quiver_1998}):

\begin{thm}\cite[Thm. 2.3]{kaledin_symplectic_2006}
    A normal symplectic singularity has finitely many symplectic leaves, all of which are algebraic.
\end{thm}
By \cite[Cor. 1.4]{crawley-boevey_geometry_2001}, quiver varieties are irreducible. In particular, this means that there is only one open symplectic leaf. As quiver varieties are \emph{conical} symplectic singularities, there is at most one 0d leaf. We thus need to study the intermediate leaves, whenever they exist.

The symplectic leaves of quiver varieties can be described using the following combinatorial datum.

\subsection{Representation types} 

For each point \(x\in\QV{\root}\) there is a corresponding closed orbit $\mathcal{O}_x \subset \mom\inv(0)$. By \cite[Prop. 3.1]{king_moduli_1994}, points $y\in\mathcal{O}_x$ are semisimple representations of the preprojective algebra \(\Uppi(\quiver)\). Thus we have a decomposition
\begin{equation}\label{eqn:rep_type_origins}
    y=\root*[1]<n_1>\oplus\cdots\oplus\root*[k]<n_k>
\end{equation}
with $\root*\in\Sig$. Moreover, \cref{eqn:rep_type_origins} does not depend on the point $y \in \mathcal{O}_x$.
This motivates the following combinatorial notion.

\begin{dfn}\label{dfn:rtype}
  We say that \(\rtype = (n_1,\root*<1>;\cdots;n_k,\root*<k>)\) is a \emph{representation type} for \(\root\) if \(\root*<i>\in\Sig\) for all \(i\) and \(\root=\Sum{n_i\root*<i>}[i=1]<k>\), $n_i \in \BN$.
\end{dfn}


This means that each representation type $\rtype$ may be associated with the set $\QV{\root}_{\rtype}$ of points where the corresponding representations decompose precisely as $\rtype$. This stratification of $\QV{\root}$ has the following geometric meaning:

\begin{thm}\cite[Thm. 1.9]{bellamy_symplectic_2021}\label{thm:rep_type_are_leaves}
    The symplectic leaves of $\QV{\root}$ are given precisely by the representation-type strata $\QV{\root}_{\rtype}$.
\end{thm}

Whenever $\root$ itself belongs to $\Sig$, the open symplectic leaf corresponds to the trivial representation type $\rtype=(1,\root)$. At the other end, the 0d leaf has the representation type $\rtype=(\dots;\root[i], \srt[i];\dots)$, a sum of simple roots with coefficients given by the dimension vector. Combinatorially, one can think that symplectic leaves are thus subrefinements of the canonical decomposition from \cref{thm:canonical_decomposition_and_dim_of_quiver_variety}.

\subsection{Minimal degenerations between leaves and transverse slices}\label{ssec:mindegen}

Symplectic leaves $\leaf[\rtype]$ form a poset with respect to the inclusion order on the closures $\cl{\leaf}[\rtype]$. The neighbours with respect to this order provide a useful geometric invariant.

\begin{dfn}\label{dfn:mindeg}
    A \emph{minimal degeneration} is a pair of leaf closures $\cl{\leaf}[1] \subset \cl{\leaf}[2]$ such that there is no leaf closure $\cl{\leaf}[3]$ with $\cl{\leaf}[1] \subset \cl{\leaf}[3] \subset \cl{\leaf}[2]$. 
\end{dfn}
\begin{rmk}
    In light of \cref{thm:canonical_decomposition_and_dim_of_quiver_variety} and \cref{thm:rep_type_are_leaves}, minimal degenerations for quiver varieties are $1$-step subrefinements of the canonical decomposition.
\end{rmk}

Given a minimal degeneration $\cl{\leaf[1]} \subset \cl{\leaf}[2]$, one can take an étale local transverse slice $S_{(12)}$ to $\leaf[1]$ inside $\cl{\leaf}[2]$. Following \cite[\S 8]{bellamy_singularities_2025}, we view a transverse slice $(S, o)$ to a leaf $\leaf$ with \(\dim\leaf=2k\) as a pointed Poisson variety together with a local isomorphism 
\begin{equation} \label{eqn:qv_is_slice_slice_times_leaf_locally}
(\QV{\root}, p) \iso \left (\Aff<2k> \times S, (0, o) \right )
\end{equation}
of Poisson varieties for each $p \in \leaf$, where $\Aff<2k>$ is equipped with the usual symplectic structure. Since the map in \eqref{eqn:qv_is_slice_slice_times_leaf_locally} is Poisson, ${o}$ will be a symplectic leaf in $S$ and hence the leaf in $\Aff<2k> \times S$ containing $(0, o)$
is $\Aff<2k> \times {o}$. The slice is étale local by its construction in \cite[\S 4]{crawley-boevey_normality_2003}.

By minimality, $S_{(12)}$ will then be an isolated singularity. In \cite{kaledin_symplectic_2006}, Kaledin shows that the stratification by symplectic leaves is the same as stratification by singular loci, so they a give \say{local analytic picture}. Thus, in order to distinguish the varieties arising from \cref{tbl:dimvec} we will list their symplectic leaves and the minimal degenerations between them. We may do this with Hasse diagrams, wherein nodes correspond to leaves and edges to minimal degenerations (via the corresponding isolated singularities).

\begin{exm}\label{exm:hanany}
    For \(k\in\Int[0]\) consider the dimension vector
    \begin{equation*}
        \begin{tikzcd}
            \& \& \& \& k+3 \ar[d] \\
            1 \ar[r] \& 2 \ar[r] \& \cdots\cdots \ar[r] \& 2k+5 \ar[r] \& 2k+6 \ar[r] \& k+4 \ar[r] \& \bal{2}.
        \end{tikzcd}
    \end{equation*}
    The unbalanced vertex has balance \(-k\) and weight $2$, so that the associated quiver variety has dimension \(2k+2\). If \(k=0\) then we see that this is just the Kleinian singularity of type \(\E[8]\) and taking \(k=1\) recovers the 4d exceptional case \(\III[b]\). In the physics literature, dimension vectors in this family (along with several others we recover\footnote{For example: (2.2) --- \(\III[b]\), (3.1) --- \(\III(\E[8])\), (4.3) --- \(\III[a]\), (5.1) --- \(\I[b]\), (6.3) --- \(\III(\E[7])\), (7.1) --- \(\I[a]\), (8.1) --- \(\I(m)\)} in 4d) were studied in \cite{ferlito_3d_2018}. The Hasse diagram of this family is
    \begin{equation*}
        \begin{tikzcd}
            \leaf[2k+2] \ar[r, "{\D[2k+8]}"] \&  \leaf[2k] \ar[r, "{\D[2(k-1)+8]}"] \& \cdots\cdots \ar[r, "{\D[10]}"] \& \leaf[2] \ar[r, "{\E[8]}"] \& \leaf[0],
        \end{tikzcd} 
    \end{equation*}
    and one property of interest for these diagrams is \emph{inversion}, meaning that knowing the Hasse diagram for the Coulomb branch allows one to recover the Hasse diagram for the Higgs branch. For more information about such a conjecture in broader sense, see \cite[\S 10.4]{braden_quantizations_2016}. 
\end{exm}

There is work \cite{bellamy_minimal_nodate} in progress by Bellamy and Schedler to classify the possible isolated singularities that occur as minimal degenerations between symplectic leaf closures in quiver varieties. These include minimal nilpotent orbit closures and Kleinian singularities, both of which appear for 4d quiver varieties. In particular, we encounter only three special cases of degenerations between leaves:

\begin{itemize}
    \item \(\dim4\) --- A degeneration between $\leaf[0]$ and $\leaf[4]$.
    \item \(\dim2\) --- A degeneration between $\leaf[0]$ and some 2d leaf \(\leaf[2i]\).
    \item \(\codim2\) --- A degeneration between some 2d leaf $\leaf[2i]$ and $\leaf[4]$.
\end{itemize}

As we will show in \cref{ssec:decompositions}, the first case is uncommon and includes only two varieties. For the second case, the transverse slice to the point $\leaf[0]$ in $\leaf[2i]$ will be its closure $\cl{\leaf[2i]}$. One therefore needs to describe the closures of all 2d leaves. For the third case we need to describe transverse slices to the symplectic leaves of codimension (coincidentally, dimension) \(2\), and this can also be done combinatorially.

\begin{dfn}\cite[Dfn. 1.18]{bellamy_symplectic_2021} \label{dfn:isotropic_decomposition}
    A root decomposition $\rtype: \root= \root*<1> + \ldots + \root*<m>+ n_1 \root**<1>+ \ldots + n_k \root**<k>  $ is called \emph{isotropic} if 
    \begin{enumerate}
        \item $\root*<i>, \root**<j> \in \Sig $,
        \item $\p(\root*<i>)=1$,
        \item $\root**<j>$ are real and pairwise distinct,
        \item the \emph{undoubled quiver} $\quiver[\rtype]''$ is extended \(\A\D\E\), where $\overline{\quiver[\rtype]''}$ is the quiver with $m+k$ loop-free vertices and for \(i\neq j\), $-(\root<i>, \root<j>)$ arrows \(i\rightarrow j\), where $\root<i>, \root<j>  \in \set{\root*<1>,\dots,\root*<m>,\root**<1>,\dots,\root**<k>}$,
        \item the dimension vector $(1, \dots, 1, n_1, \dots , n_k)$ on $\quiver[\rtype]''$ is the corresponding minimal imaginary root.
    \end{enumerate}
    The quiver in (4) is known as the \emph{ext-quiver}, due to the number of arrows coinciding with the dimension of \(\operatorname{Ext}^1\) spaces between vertex simple modules over \(\Uppi(\quiver[\rtype]'')\).
\end{dfn}

\subsubsection{Codimension two degenerations}

\begin{thm} \cite[Thm. 1.20]{bellamy_symplectic_2021} \label{thm:isotropic_rep_type_are_codim2_leaves}
    Let $\root \in \Sig$ be imaginary. Then the codimension two strata of  $\QV{\root}$ are in bijection with the isotropic decompositions of $\root$. Moreover, the slice du Val singularity corresponds to the ext-quiver $Q''$.
\end{thm}

Thus to understand such degenerations, we need to find all isotropic decompositions of our dimension vectors. This is carried out in \cref{ssec:decompositions}.

\subsubsection{Dimension two degenerations} \label{sssec:normalisations_of_leaves}

Recall from \cref{thm:rep_type_are_leaves} that symplectic leaves of quiver varieties correspond to representation-type datum.

In \cite[Thm. 2.4]{bellamy_singularities_2025} the authors, via representation types, construct normalisations of closures of symplectic leaves using products of quotients of \say{smaller} quiver varieties. Next they ask whether the leaves themselves are normal. It turns out that the answer is \say{no} in general. However, we \emph{do} have normality for a special type of leaf.


\begin{ppn}\cite[Ppn. 2.8]{bellamy_singularities_2025} \label{ppn:ADE_leaf_closure}
    If a representation type \(\rtype = (1,\mir;n_2,\root**<2>;\dots;n_r,\root**<r>)\) has \(\p(\mir)=1\) and \(\p(\root**<k>)=0\), then \(\cl{\leaf}[\rtype]\) is normal and isomorphic to the Kleinian singularity associated with \(\mir\).
\end{ppn}

The representation types above are usually not isotropic. However, they are in the case \(\on{p}(\root) = 2\). As every non-trivial representation type of \(\root\) is of the form \(\rtype = (1,\mir;n_2,\root**<2>;\cdots;n_k,\root**<k>)\) with every \(\root**<i>\) a real root (as otherwise \(\root\not\in\Sig\)), we may apply \cref{ppn:ADE_leaf_closure}. We thus proceed to classify all the possible representation types for the dimension vectors in \cref{thm:dimvec}. 

\subsection{Isotropic decompositions}\label{ssec:decompositions}

By the previous discussion we know that if \(\p(\root)=2\) and \(\mir\leq*\root\in\Sig\) is isotropic, then \(\root-\mir\) cannot admit another isotropic subroot, so we exhaustively list all such subroots \(\mir\) for each \(\root\) in our classification. We first exclude subroots of the form $\mir(\A)$ by looking for cycles, then of the form $\mir(\D)$ by looking for pairs of branch vertices, and lastly of the form $\mir(\E)$. We demonstrate this routine in an example.

\begin{exm}\label{exm:isodecomp}
    Consider the Type \(\I\) dimension vector
    \begin{equation*}
        \root = \I(\E[7], n) = \begin{tikzcd}[scale cd=.8, cramped, sep=0]
        \&\&\&2 \\ \bal{1}\ar[drr, bend right, "1"']\&2\&3\&4\&3\&2\&\bal{1}\ar[dll, bend left, "n"] \\ \&\&1\&\cdots\&1
        \end{tikzcd}
    \end{equation*}
    As $\root$ contains a cycle, it must be a subroot of the form \(\mir(\A[n])\) for some $n$, explicitly
    \begin{equation*}
        \begin{tikzcd}[scale cd=.8, cramped, sep=0]
        \&\&\&2 \\ \bal{1}\ar[drr, bend right, "1"']\&2\&3\&4\&3\&2\&\bal{1}\ar[dll, bend left, "n"] \\ \&\&1\&\cdots\&1
        \end{tikzcd} \quad=\quad \ubrace*{\begin{tikzcd}[scale cd=.8, cramped, sep=0]
        \&\&\&0 \\ 1\ar[drr, bend right]\&1\&1\&1\&1\&1\&1\ar[dll, bend left] \\ \&\&1\&\cdots\&1
        \end{tikzcd}}{\mir(\A[n+5])} \quad+\quad \ubrace*{\begin{tikzcd}[scale cd=.8, cramped, sep=0]
        \&\&\&2 \\ 0\ar[drr, bend right]\&1\&2\&3\&2\&1\&0\ar[dll, bend left] \\ \&\&0\&\cdots\&0
        \end{tikzcd}}{\hrr(\E[6])}
    \end{equation*}
    Next, $\root$ has only one branch vertex, so there are no subroots of the form $\mir(\D)$. Lastly, one can see that $\mir(\E[6])\not\leq*\root$ (branches too short) and $\mir(\E[8])$ (branch vertex of too low weight), but does contain  $\mir(\E[7])$ provided that \(n\geq*1\). This gives the following decomposition:
    \begin{equation*}
        \begin{tikzcd}[scale cd=.8, cramped, sep=0]
            \&\&\&2 \\ \bal{1}\ar[drr, bend right, "1"']\&2\&3\&4\&3\&2\&\bal{1}\ar[dll, bend left, "n"] \\ \&\&1\&\cdots\&1
        \end{tikzcd} \quad=\quad \ubrace*{\begin{tikzcd}[scale cd=.8, cramped, sep=0]
            \&\&\&2 \\ 1\ar[drr, bend right]\&2\&3\&4\&3\&2\&1\ar[dll, bend left] \\ \&\&0\&\cdots\&0
        \end{tikzcd}}{\mir(\E[7])} \quad+\quad \ubrace*{\begin{tikzcd}[scale cd=.8, cramped, sep=0]
            \&\&\&0 \\ 0\ar[drr, bend right]\&0\&0\&0\&0\&0\&0\ar[dll, bend left] \\ \&\&1\&\cdots\&1
        \end{tikzcd}}{\hrr(\A[n-1])}
    \end{equation*}
\end{exm}

We proceed in this way for the remaining dimension vectors of \cref{thm:dimvec}, proceeding by type.

\subsubsection{Type \(\I\)}\label{sssec:typei}

We start with \(\I(\ade, n)\), which by construction contains \(\mir(\ade)\) (provided that \(n\geq*1\)) and some \(\mir(\A)\). If \(\ade=\A[m]\), then we have:

\begin{equation*}
    \ubrace*{\begin{tikzcd}[scale cd=.8, cramped, sep=0]
        \&1\&\cdots\&1 \\ \bal{1}\&1\&\cdots\&1\&\bal{1} \\ \&1\&\cdots\&1
    \end{tikzcd}}{\I(\ell,m,n)} \quad=\quad \ubrace*{\begin{tikzcd}[scale cd=.8, cramped, sep=0]
        \&1\&\cdots\&1 \\ 1\&1\&\cdots\&1\&1 \\ \&0\&\cdots\&0
    \end{tikzcd}}{\mir(\A[\ell+m-1])} \quad+\quad \ubrace*{\begin{tikzcd}[scale cd=.8, cramped, sep=0]
        \&0\&\cdots\&0 \\ 0\&0\&\cdots\&0\&0 \\ \&1\&\cdots\&1
    \end{tikzcd}}{\hrr(\A[n-1])}
\end{equation*}

By symmetry there are two other decompositions with isotropic subroots \(\mir(\A[\ell+n-1])\) (whenever \(m\geq*1\)) and \(\mir(\A[m+n-1])\) (whenever \(\ell\geq*1\)). For \(\ade=\D[m]\), we see that:

\begin{equation*}
    \ubrace*{\begin{tikzcd}[scale cd=.8, cramped, sep=0]
        1\&\&\&\&1 \\ \&2\&\cdots\&2 \\ \bal{1}\ar[dr, bend right, "1"']\&\&\&\&\bal{1}\ar[dl, bend left, "n"] \\ \&1\&\cdots\&1
    \end{tikzcd}}{\I(\D[m],n)} \quad=\quad \ubrace*{\begin{tikzcd}[scale cd=.8, cramped, sep=0]
        0\&\&\&\&0 \\ \&1\&\cdots\&1 \\ 1\ar[dr, bend right]\&\&\&\&1\ar[dl, bend left] \\ \&1\&\cdots\&1
    \end{tikzcd}}{\mir(\A[m+n-3])} \quad+\quad \ubrace*{\begin{tikzcd}[scale cd=.8, cramped, sep=0]
        1\&\&\&\&1 \\ \&1\&\cdots\&1 \\ 0\ar[dr, bend right]\&\&\&\&0\ar[dl, bend left] \\ \&0\&\cdots\&0
    \end{tikzcd}}{\hrr(\A[m-1])} \quad\underbrace{=}_{n\geq*1}\quad \ubrace*{\begin{tikzcd}[scale cd=.8, cramped, sep=0]
        1\&\&\&\&1 \\ \&2\&\cdots\&2 \\ 1\ar[dr, bend right]\&\&\&\&1\ar[dl, bend left] \\ \&0\&\cdots\&0
    \end{tikzcd}}{\mir(\D[m])} \quad+\quad \ubrace*{\begin{tikzcd}[scale cd=.8, cramped, sep=0]
        0\&\&\&\&0 \\ \&0\&\cdots\&0 \\ 0\ar[dr, bend right]\&\&\&\&0\ar[dl, bend left] \\ \&1\&\cdots\&1
    \end{tikzcd}}{\hrr(\A[n-1])}
\end{equation*}

\makeatletter\DeclareRobustCommand\vdots{\vbox{\baselineskip4\p@\lineskiplimit\z@\hbox{.}\hbox{.}\hbox{.}}}\makeatother

The situation is similar for \(\I(\D[m],n)'\), but now the cycle length depends on $n$ only:

\begin{equation*}
    \ubrace*{\begin{tikzcd}[scale cd=.8, cramped, sep=0]
        1\&\&\&\&\bal{1}\&1 \\ \&2\&\cdots\&2\&\&\vdots \\ 1\&\&\&\&\bal{1}\&1
    \end{tikzcd}}{\I(\D[m],n)'}
    \quad=\quad \ubrace*{\begin{tikzcd}[scale cd=.8, cramped, sep=0]
        0\&\&\&\&1\&1 \\ \&0\&\cdots\&1\&\&\vdots \\ 0\&\&\&\&1\&1
    \end{tikzcd}}{\mir(\A[n+1])} \quad+\quad \ubrace*{\begin{tikzcd}[scale cd=.8, cramped, sep=0]
        1\&\&\&\&0\&0 \\ \&2\&\cdots\&1\&\&\vdots \\ 1\&\&\&\&0\&0
    \end{tikzcd}}{\hrr(\D[m-1])}
    \quad\underbrace{=}_{n\geq*1}\quad \ubrace*{\begin{tikzcd}[scale cd=.8, cramped, sep=0]
        1\&\&\&\&1\&0 \\ \&2\&\cdots\&2\&\&\vdots \\ 1\&\&\&\&1\&0
    \end{tikzcd}}{\mir(\D[m])} \quad+\quad \ubrace*{\begin{tikzcd}[scale cd=.8, cramped, sep=0]
        0\&\&\&\&0\&1 \\ \&0\&\cdots\&0\&\&\vdots \\ 0\&\&\&\&0\&1
    \end{tikzcd}}{\hrr(\A[n-1])}
\end{equation*}

This leaves \(\ade=\E[6]\) (the case \(\E[7]\) was covered in \cref{exm:isodecomp}):

\begin{equation*}
    \ubrace*{\begin{tikzcd}[scale cd=.8, cramped, sep=0]
        \&\&1 \\ \&\&2 \\ \bal{1}\ar[dr, bend right, "1"']\&2\&3\&2\&\bal{1}\ar[dl, bend left, "n"] \\ \&1\&\cdots\&1
    \end{tikzcd}}{\I(\E[6],n)}
    \quad=\quad \ubrace*{\begin{tikzcd}[scale cd=.8, cramped, sep=0]
        \&\&0 \\ \&\&0 \\ 1\ar[dr, bend right]\&1\&1\&1\&1\ar[dl, bend left] \\ \&1\&\cdots\&1
    \end{tikzcd}}{\mir(\A[n+3])} \quad+\quad \ubrace*{\begin{tikzcd}[scale cd=.8, cramped, sep=0]
        \&\&1 \\ \&\&2 \\ 0\ar[dr, bend right]\&1\&2\&1\&0\ar[dl, bend left] \\ \&0\&\cdots\&0
    \end{tikzcd}}{\hrr(\D[5])}
    \quad\underbrace{=}_{n\geq*1}\quad \ubrace*{\begin{tikzcd}[scale cd=.8, cramped, sep=0]
        \&\&1 \\ \&\&2 \\ 1\ar[dr, bend right]\&2\&3\&2\&1\ar[dl, bend left] \\ \&0\&\cdots\&0
    \end{tikzcd}}{\mir(\E[6])} \quad+\quad \ubrace*{\begin{tikzcd}[scale cd=.8, cramped, sep=0]
        \&\&0 \\ \&\&0 \\ 0\ar[dr, bend right]\&0\&0\&0\&0\ar[dl, bend left] \\ \&1\&\cdots\&1
    \end{tikzcd}}{\hrr(\A[n-1])}
\end{equation*}

We have thus completed the analysis of \(\I(\ade,n)\), and move onto the infinite family \(\I(m)\).

\begin{equation*}
    \ubrace*{\begin{tikzcd}[scale cd=.8, cramped, sep=0]
        \&\&\bal{1}\&\&\bal{1} \\ 1\&2\&3\&\cdots\&3\&2\&1
    \end{tikzcd}}{\I(m)} \quad=\quad \ubrace*{\begin{tikzcd}[scale cd=.8, cramped, sep=0]
        \&\&1\&\&1 \\ 0\&1\&2\&\cdots\&2\&1\&0
    \end{tikzcd}}{\mir(\D[m])} \quad+\quad \ubrace*{\begin{tikzcd}[scale cd=.8, cramped, sep=0]
        \&\&0\&\&0 \\ 1\&1\&1\&\cdots\&1\&1\&1
    \end{tikzcd}}{\hrr(\A[m+1])}
\end{equation*}

The first exception of Type \(\I\) is \(\I[a]\), admitting the sole isotropic subroot of type \(\D[5]\).

\begin{equation*}
    \ubrace*{\begin{tikzcd}[scale cd=.8, cramped, sep=0]
        \bal{1}\&\&2 \\ \&3\&4\&3\&2\&1 \\ \bal{1}
    \end{tikzcd}}{\I[a]} \quad=\quad \ubrace*{\begin{tikzcd}[scale cd=.8, cramped, sep=0]
        1\&\&1 \\ \&2\&2\&1\&0\&0 \\ 1
    \end{tikzcd}}{\mir(\D[5])} \quad+\quad \ubrace*{\begin{tikzcd}[scale cd=.8, cramped, sep=0]
        0\&\&1 \\ \&1\&2\&2\&2\&1 \\ 0
    \end{tikzcd}}{\hrr(\D[6])}
\end{equation*}

Next we consider the exception \(\I[b]\) which contains only a type \(\E[6]\) isotropic subroot.

\begin{equation*}
    \ubrace*{\begin{tikzcd}[scale cd=.8, cramped, sep=0]
        \bal{1}\&3 \\ \&\&5\&4\&3\&2\&1 \\ \bal{1}\&3
    \end{tikzcd}}{\I[b]} \quad=\quad \ubrace*{\begin{tikzcd}[scale cd=.8, cramped, sep=0]
        1\&2 \\ \&\&3\&2\&1\&0\&0 \\ 1\&2
    \end{tikzcd}}{\mir(\E[6])} \quad+\quad \ubrace*{\begin{tikzcd}[scale cd=.8, cramped, sep=0]
        0\&1 \\ \&\&2\&2\&2\&2\&1 \\ 0\&1
    \end{tikzcd}}{\hrr(\D[7])}
\end{equation*}

The final exceptional case is \(\I[c]\). We have only a \(\D[7]\) isotropic subroot.

\begin{equation*}
    \ubrace*{\begin{tikzcd}[scale cd=.8, cramped, sep=0]
        \bal{1}\&\&\&\&3 \\ \&3\&4\&5\&6\&4\&2 \\ \bal{1}
    \end{tikzcd}}{\I[c]} \quad=\quad \ubrace*{\begin{tikzcd}[scale cd=.8, cramped, sep=0]
        1\&\&\&\&1 \\ \&2\&2\&2\&2\&1\&0 \\ 1
    \end{tikzcd}}{\mir(\D[7])} \quad+\quad \ubrace*{\begin{tikzcd}[scale cd=.8, cramped, sep=0]
        0\&\&\&\&2 \\ \&1\&2\&3\&4\&3\&2 \\ 0
    \end{tikzcd}}{\hrr(\E[7])}
\end{equation*}

\subsubsection{Type \(\II\)}

We start with \(\II(\ade<1>,\ade<2>)\). The immediate isotropic subroots in this construction are that of \(\ade<1>\) and \(\ade<2>\). However notice that \(\II(\ade<1>, \ade<2>)=\mir(\ade<2>) + \hrr(\ade<1>)\) and dually if we subtract \(\mir(\ade<2>)\). To see that there are no other isotropic subroots, notice that every minimal imaginary root only has a weight one vertex at the end of a chain. The construction $\II(\ade)$ is covered by the $n=0$ case of $\I(\ade, n)$, so we may proceed to the infinite family \(\II(m)\). The only isotropic decomposition is

\begin{equation*}
    \ubrace*{\begin{tikzcd}[scale cd=.8, cramped, sep=0]
        2\&\&\&\bal{1} \\ \&4\&\cdots\&4\&3\&2\&1 \\ 2
    \end{tikzcd}}{\II(m)} \quad=\quad \ubrace*{\begin{tikzcd}[scale cd=.8, cramped, sep=0]
        1\&\&\&1 \\ \&2\&\cdots\&2\&1\&0\&0 \\ 1
    \end{tikzcd}}{\mir(\D[m])} \quad+\quad \ubrace*{\begin{tikzcd}[scale cd=.8, cramped, sep=0]
        1\&\&\&0 \\ \&2\&\cdots\&2\&2\&2\&1 \\ 1
    \end{tikzcd}}{\hrr(\D[m+2])}
\end{equation*}

The two exceptional cases of Type \(\II\) admit an isotropic subroot of type \(\D[6]\) and \(\E[7]\) respectively.

\begin{align*}
    \ubrace*{\begin{tikzcd}[scale cd=.8, cramped, sep=0]
        \&2\&\&3 \\ \bal{1}\&4\&5\&6\&4\&2
    \end{tikzcd}}{\II[a]} \quad&=\quad \ubrace*{\begin{tikzcd}[scale cd=.8, cramped, sep=0]
        \&1\&\&1 \\ 1\&2\&2\&2\&1\&0
    \end{tikzcd}}{\mir(\D[6])} \quad+\quad \ubrace*{\begin{tikzcd}[scale cd=.8, cramped, sep=0]
        \&1\&\&2 \\ 0\&2\&3\&4\&3\&2
    \end{tikzcd}}{\hrr(\E[7])}\\
       \ubrace*{\begin{tikzcd}[scale cd=.8, cramped, sep=0]
        \&\&\&5 \\ \bal{1}\&4\&7\&10\&8\&6\&4\&2
    \end{tikzcd}}{\II[b]} \quad&=\quad \ubrace*{\begin{tikzcd}[scale cd=.8, cramped, sep=0]
        \&\&\&2 \\ 1\&2\&3\&4\&3\&2\&1\&0
    \end{tikzcd}}{\mir(\E[7])} \quad+\quad \ubrace*{\begin{tikzcd}[scale cd=.8, cramped, sep=0]
        \&\&\&3 \\ 0\&2\&4\&6\&5\&4\&3\&2
    \end{tikzcd}}{\hrr(\E[8])} 
\end{align*}

\subsubsection{Type \(\III\)}

By construction, subtracting \(\mir(\ade)\) from $\III(\ade)$ results in \(\srt[i]\), where \(i\) is the vertex of weight \(2\) which was unbalanced.

\begin{align*}
    \ubrace*{\begin{tikzcd}[scale cd=.8, cramped, sep=0]
        1\&\&\&\&\&1\&\&\&\&\&\&1 \\ \&2\&2\&\cdots\&2\&\bal{2}\&2\&\cdots\&2\&2 \\ 1\&\&\&\&\&\&\&\&\&\&\&1
    \end{tikzcd}}{\III(\D[m],i)} \quad&=\quad \ubrace*{\begin{tikzcd}[scale cd=.8, cramped, sep=0]
        1\&\&\&\&\&0\&\&\&\&\&\&1 \\ \&2\&2\&\cdots\&2\&2\&2\&\cdots\&2\&2 \\ 1\&\&\&\&\&\&\&\&\&\&\&1
    \end{tikzcd}}{\mir(\D[m])} \quad+\quad \ubrace*{\begin{tikzcd}[scale cd=.8, cramped, sep=0]
        0\&\&\&\&\&1\&\&\&\&\&\&0 \\ \&0\&0\&\cdots\&0\&0\&0\&\cdots\&0\&0 \\ 0\&\&\&\&\&\&\&\&\&\&\&0
    \end{tikzcd}}{\hrr(\A[1])} \\
    \quad&=\quad \ubrace*{\begin{tikzcd}[scale cd=.8, cramped, sep=0]
        0\&\&\&\&\&1\&\&\&\&\&\&1 \\ \&0\&0\&\cdots\&1\&2\&2\&\cdots\&2\&2 \\ 0\&\&\&\&\&\&\&\&\&\&\&1
    \end{tikzcd}}{\mir(\D[m-i+1])} \quad+\quad \ubrace*{\begin{tikzcd}[scale cd=.8, cramped, sep=0]
        1\&\&\&\&\&0\&\&\&\&\&\&0 \\ \&2\&2\&\cdots\&2\&2\&1\&\cdots\&0\&0 \\ 1\&\&\&\&\&\&\&\&\&\&\&0
    \end{tikzcd}}{\hrr(\D[i+1])} \\
    \quad&=\quad \ubrace*{\begin{tikzcd}[scale cd=.8, cramped, sep=0]
        1\&\&\&\&\&1\&\&\&\&\&\&0 \\ \&2\&2\&\cdots\&2\&2\&1\&\cdots\&0\&0 \\ 1\&\&\&\&\&\&\&\&\&\&\&0
    \end{tikzcd}}{\mir(\D[i+3])} \quad+\quad \ubrace*{\begin{tikzcd}[scale cd=.8, cramped, sep=0]
        0\&\&\&\&\&0\&\&\&\&\&\&1 \\ \&0\&0\&\cdots\&0\&0\&1\&\cdots\&2\&2 \\ 0\&\&\&\&\&\&\&\&\&\&\&1
    \end{tikzcd}}{\hrr(\D[m-i-1])}
\end{align*}

Notice that \(\III(\D[m],1)\) admits \(\mir(\D[m])\) in three different ways (five if \(m=4\)), so the last two decompositions will \say{split}. However, this is consistent with the fact \(\D[2]\iso\A[1]\times\A[1]\). If \(\ade=\E\), we have

\begin{align*}
    \ubrace*{\begin{tikzcd}[scale cd=.8, cramped, sep=0]
        \&\&\&1 \\ \&\&\&\bal{2} \\ 1\&2\&3\&4\&3\&2\&1
    \end{tikzcd}}{\III(\E[7])} \quad&=\quad \ubrace*{\begin{tikzcd}[scale cd=.8, cramped, sep=0]
        \&\&\&0 \\ \&\&\&2 \\ 1\&2\&3\&4\&3\&2\&1
    \end{tikzcd}}{\mir(\E[7])} \quad+\quad \ubrace*{\begin{tikzcd}[scale cd=.8, cramped, sep=0]
        \&\&\&1 \\ \&\&\&0 \\ 0\&0\&0\&0\&0\&0\&0
    \end{tikzcd}}{\hrr(\A[1])} \\
    \quad&=\quad \ubrace*{\begin{tikzcd}[scale cd=.8, cramped, sep=0]
        \&\&\&1 \\ \&\&\&2 \\ 0\&1\&2\&3\&2\&1\&0
    \end{tikzcd}}{\mir(\E[6])} \quad+\quad \ubrace*{\begin{tikzcd}[scale cd=.8, cramped, sep=0]
        \&\&\&0 \\ \&\&\&0 \\ 1\&1\&1\&1\&1\&1\&1
    \end{tikzcd}}{\hrr(\A[7])}
\end{align*}

\begin{align*}
    \ubrace*{\begin{tikzcd}[scale cd=.8, cramped, sep=0]
        \&\&\&\&\&3 \\ 1\&2\&3\&4\&5\&6\&4\&\bal{2}\&1
    \end{tikzcd}}{\III(\E[8])} \quad&=\quad \ubrace*{\begin{tikzcd}[scale cd=.8, cramped, sep=0]
        \&\&\&\&\&3 \\ 1\&2\&3\&4\&5\&6\&4\&2\&0
    \end{tikzcd}}{\mir(\E[8])} \quad+\quad \ubrace*{\begin{tikzcd}[scale cd=.8, cramped, sep=0]
        \&\&\&\&\&0 \\ 0\&0\&0\&0\&0\&0\&0\&0\&1
    \end{tikzcd}}{\hrr(\A[1])}\\
    \quad&=\quad \ubrace*{\begin{tikzcd}[scale cd=.8, cramped, sep=0]
        \&\&\&\&\&2 \\ 0\&0\&1\&2\&3\&4\&3\&2\&1
    \end{tikzcd}}{\mir(\E[7])} \quad+\quad \ubrace*{\begin{tikzcd}[scale cd=.8, cramped, sep=0]
        \&\&\&\&\&1 \\ 1\&2\&2\&2\&2\&2\&1\&0\&0
    \end{tikzcd}}{\hrr(\D[8])}
\end{align*}

For the infinite family \(\III(m,n)\), (which we recall encompasses \(\III(\E[6])'=\III(\E[6]), \III(\E[7])'\) and \(\III(\E[8])'\)) we display an extra vertex in the chains for clarity.

\begin{align*}
    \ubrace*{\begin{tikzcd}[scale cd=.8, cramped, sep=0]
        \&\&\&\&\&1 \\ 1\&\&\&\&\&2 \\ \&2\&\cdots\&2\&\bal{2}\&3\&2\&1 \\ 1
    \end{tikzcd}}{\III(m,6)} \quad&=\quad \ubrace*{\begin{tikzcd}[scale cd=.8, cramped, sep=0]
        \&\&\&\&\&0 \\ 1\&\&\&\&\&1 \\ \&2\&\cdots\&2\&2\&2\&1\&0 \\ 1
    \end{tikzcd}}{\mir(\D[m+1])} \quad+\quad \ubrace*{\begin{tikzcd}[scale cd=.8, cramped, sep=0]
        \&\&\&\&\&1 \\ 0\&\&\&\&\&1 \\ \&0\&\cdots\&0\&0\&1\&1\&1 \\ 0
    \end{tikzcd}}{\hrr(\A[5])} \\
    \quad&=\quad \ubrace*{\begin{tikzcd}[scale cd=.8, cramped, sep=0]
        \&\&\&\&\&1 \\ 0\&\&\&\&\&2 \\ \&0\&\cdots\&1\&2\&3\&2\&1 \\ 0
    \end{tikzcd}}{\mir(\E[6])} \quad+\quad \ubrace*{\begin{tikzcd}[scale cd=.8, cramped, sep=0]
        \&\&\&\&\&0 \\ 1\&\&\&\&\&0 \\ \&2\&\cdots\&1\&0\&0\&0\&0 \\ 1
    \end{tikzcd}}{\hrr(\D[m-2])}
\end{align*}
\begin{align*}
    \ubrace*{\begin{tikzcd}[scale cd=.8, cramped, sep=0]
        1\&\&\&\&\&\&2 \\ \&2\&\cdots\&2\&\bal{2}\&3\&4\&3\&2\&1 \\ 1
    \end{tikzcd}}{\III(m,7)} \quad&=\quad \ubrace*{\begin{tikzcd}[scale cd=.8, cramped, sep=0]
        1\&\&\&\&\&\&1 \\ \&2\&\cdots\&2\&2\&2\&2\&1\&0\&0 \\ 1
    \end{tikzcd}}{\mir(\D[m+2])} \quad+\quad \ubrace*{\begin{tikzcd}[scale cd=.8, cramped, sep=0]
        0\&\&\&\&\&\&1 \\ \&0\&\cdots\&0\&0\&1\&2\&2\&2\&1 \\ 0
    \end{tikzcd}}{\hrr(\D[6])} \\\\
    \quad&=\quad \ubrace*{\begin{tikzcd}[scale cd=.8, cramped, sep=0]
        0\&\&\&\&\&\&2 \\ \&0\&\cdots\&1\&2\&3\&4\&3\&2\&1 \\ 0
    \end{tikzcd}}{\mir(\E[7])} \quad+\quad \ubrace*{\begin{tikzcd}[scale cd=.8, cramped, sep=0]
        1\&\&\&\&\&\&0 \\ \&2\&\cdots\&1\&0\&0\&0\&0\&0\&0 \\ 1
    \end{tikzcd}}{\hrr(\D[m-2])}
\end{align*}
\begin{align*}
    \ubrace*{\begin{tikzcd}[scale cd=.8, cramped, sep=0]
        1\&\&\&\&\&\&\&\&3 \\ \&2\&\cdots\&2\&\bal{2}\&3\&4\&5\&6\&4\&2 \\ 1
    \end{tikzcd}}{\III(m,8)} \quad&=\quad \ubrace*{\begin{tikzcd}[scale cd=.8, cramped, sep=0]
        1\&\&\&\&\&\&\&\&1 \\ \&2\&\cdots\&2\&2\&2\&2\&2\&2\&1\&0 \\ 1
    \end{tikzcd}}{\mir(\D[m+4])} \quad+\quad \ubrace*{\begin{tikzcd}[scale cd=.8, cramped, sep=0]
        0\&\&\&\&\&\&\&\&2 \\ \&0\&\cdots\&0\&0\&1\&2\&3\&4\&3\&2 \\ 0
    \end{tikzcd}}{\hrr(\E[7])} \\
    \quad&=\quad \ubrace*{\begin{tikzcd}[scale cd=.8, cramped, sep=0]
        0\&\&\&\&\&\&\&\&3 \\ \&0\&\cdots\&1\&2\&3\&4\&5\&6\&4\&2 \\ 0
    \end{tikzcd}}{\mir(\E[8])} \quad+\quad \ubrace*{\begin{tikzcd}[scale cd=.8, cramped, sep=0]
        1\&\&\&\&\&\&\&\&0 \\ \&2\&\cdots\&1\&0\&0\&0\&0\&0\&0\&0 \\ 1
    \end{tikzcd}}{\hrr(\D[m-2])}
\end{align*}

If \(m=4\), observe that the second decomposition occurs twice, as there are two ways to subtract \(\mir(\E[n])\), analogously to \(\III(\D[m],i)\). Finally we have the Type \(\III\) exceptions:

\begin{align*}
    \ubrace*{\begin{tikzcd}[scale cd=.8, cramped, sep=0]
        \&\&\&\&\bal{2} \\ 1\&2\&3\&4\&5\&4\&3\&2\&1
    \end{tikzcd}}{\III[a]} \quad&=\quad \ubrace*{\begin{tikzcd}[scale cd=.8, cramped, sep=0]
        \&\&\&\&2 \\ 0\&1\&2\&3\&4\&3\&2\&1\&0
    \end{tikzcd}}{\mir(\E[7])} \quad+\quad \ubrace*{\begin{tikzcd}[scale cd=.8, cramped, sep=0]
        \&\&\&\&0 \\ 1\&1\&1\&1\&1\&1\&1\&1\&1
    \end{tikzcd}}{\hrr(\A[9])} \\
    \ubrace*{\begin{tikzcd}[scale cd=.8, cramped, sep=0]
        \&\&4 \\ \bal{2}\&5\&8\&7\&6\&5\&4\&3\&2\&1
    \end{tikzcd}}{\III[b]} \quad&=\quad \ubrace*{\begin{tikzcd}[scale cd=.8, cramped, sep=0]
        \&\&3 \\ 2\&4\&6\&5\&4\&3\&2\&1\&0\&0
    \end{tikzcd}}{\mir(\E[8])} \quad+\quad \ubrace*{\begin{tikzcd}[scale cd=.8, cramped, sep=0]
        \&\&1 \\ 0\&1\&2\&2\&2\&2\&2\&2\&2\&1
    \end{tikzcd}}{\hrr(\D[10])}
\end{align*}

\subsubsection{Summary}

Noting that \(\dim2\) degenerations correspond to \(\mir\) and \(\codim2\) degenerations correspond to \(\hrr\), we tabulate the decomposition data for later.

\begin{longtable}{|C|C|C|C|C|}\caption{Isotropic decompositions for \(\root\in\Sig\) satisfying \(p(\root)=2\)}\label{tbl:decomp}\\\hline
    \rowcolor{gray!50} \root & \text{Restriction(s)} & \rank & \mir & \hrr \\\hline \endfirsthead
    \caption[]{(continued)}\\\hline\rowcolor{gray!50}
    \root & \text{Restriction(s)} & \rank & \mir & \hrr \\\hline \endhead
    \hline\endfoot
    \multirow{4}{*}{$\I(\ell,m,n)$} & 1\leq*\ell\leq m\leq n & \multirow{4}{*}{$\ell+m+n-1$} & \A[m+n-1];\,\A[\ell+n-1];\,\A[\ell+m-1] & \A[\ell-1];\,\A[m-1];\,\A[n-1] \\
    & 1=\ell\leq* m\leq n & & \A[m-1];\,\A[n-1] & \A[n];\,\A[m] \\
    & 1=\ell=m\leq*n & & \A[1] & \A[n-1] \\
    & 1=\ell=m=n & & & \\
    \rowcolor{gray!20} & n\geq2 & & \D[m];\,\A[m+n-3] & \A[n-1];\,\A[m-1] \\
    \rowcolor{gray!20} \multirow{-2}{*}{$\I(\D[m],n)$} & n\in\set{0,1} & \multirow{-2}{*}{$m+n$} & \A[m+n-3] & \A[m-1] \\
    \multirow{2}{*}{$\I(\D[m],n)'$} & n\geq2 & \multirow{2}{*}{$m+n$} & \D[m];\,\A[n+1] & \A[n-1];\,\D[m-1] \\
    & n\in\set{0,1} & & \A[n+1] & \D[m-1] \\
    \rowcolor{gray!20} & n\geq2 & & \E[6];\,\A[n+3] & \A[n-1];\,\D[5] \\
    \rowcolor{gray!20} \multirow{-2}{*}{$\I(\E[6],n)$} & n\in\set{0,1} & \multirow{-2}{*}{$n+6$} & \A[n+3] & \D[5] \\
    \multirow{2}{*}{$\I(\E[7],n)$} & n\geq2 & \multirow{2}{*}{$n+7$} & \E[7];\,\A[n+5] & \A[n-1];\,\E[6] \\
    & n\in\set{0,1} & & \A[n+5] & \E[6] \\
    \rowcolor{gray!20} \I(m) & m\geq4 & m+3 & \D[m] & \A[m+1] \\
    \I[a] & & 8 & \D[5] & \D[6] \\
    \rowcolor{gray!20} \I[b] & & 9 & \E[6] & \D[7] \\
    \I[c] & & 9 & \D[7] & \E[7] \\
    \rowcolor{gray!20} & \ade<1>\neq\A[0]\neq\ade<2> & & \ade<1>;\,\ade<2> & \ade<2>;\,\ade<1> \\
    \rowcolor{gray!20} & \ade<1>=\A[0]\neq\ade<2> & & \A[0] & \ade<2> \\
    \rowcolor{gray!20} \multirow{-3}{*}{$\II(\ade<1>, \ade<2>)$} & \ade<1>=\A[0]=\ade<2> & \multirow{-3}{*}{$\rank<(1)>+\rank<(2)>+1$} & & \\
    \II(m) & m\geq4 & m+3 & \D[m] & \D[m+2] \\
    \rowcolor{gray!20} \II[a] & & 8 & \D[6] & \E[7] \\
    \II[b] & & 9 & \E[7] & \E[8] \\
    \rowcolor{gray!20} & 2\leq i\leq\ceil{\tfrac{m-3}{2}} & & \D[m];\,\D[m-i+1];\,\D[i+3] & \A[1];\,\D[i+1];\,\D[m-i-1] \\
    \rowcolor{gray!20} & 1=i\leq\ceil{\tfrac{m-3}{2}} & & \D[m];\,\D[m];\,\D[m];\,\D[4] & \A[1];\,\A[1];\,\A[1];\,\D[m-2]\\
    \rowcolor{gray!20}\multirow{-3}{*}{$\III(\D[m],i)$} & m=4 & \multirow{-3}{*}{$m+2$} & \D[4];\,\D[4];\,\D[4];\,\D[4];\,\D[4] & \A[1];\,\A[1];\,\A[1];\,\A[1];\,\A[1]\\
    \III(\E[7]) & & 9 & \E[6];\,\E[7] & \A[7];\,\A[1] \\
    \rowcolor{gray!20} \III(\E[8]) & & 10 & \E[7];\,\E[8] & \D[8];\,\A[1] \\
    & m\geq5 & & \D[m+1];\,\E[6] & \A[5];\,\D[m-2] \\
    \multirow{-2}{*}{$\III(m,6)$} & m=4 & \multirow{-2}{*}{$m+4$} & \D[5];\,\E[6];\,\E[6] & \A[5];\,\A[1];\,\A[1] \\
    \rowcolor{gray!20} & m\geq5 & & \D[m+2];\,\E[7] & \D[6];\,\D[m-2] \\
    \rowcolor{gray!20}\multirow{-2}{*}{$\III(m,7)$} & m=4 & \multirow{-2}{*}{$m+5$} & \D[6];\,\E[7];\,\E[7] & \D[6];\,\A[1];\,\A[1] \\
    & m\geq5 & & \D[m+4];\,\E[8] & \E[7];\,\D[m-2] \\
    \multirow{-2}{*}{$\III(m,8)$} & m=4 & \multirow{-2}{*}{$m+6$} & \D[8];\,\E[8];\,\E[8] & \E[7];\,\A[1];\,\A[1] \\
    \rowcolor{gray!20} \III[a] & & 10 & \E[7] & \A[9] \\
    \III[b] & & 11 & \E[8] & \D[10]
\end{longtable}

\subsection{Namikawa's Weyl group}

In \cite{namikawa_poisson_2010}, Namikawa defines a finite reflection group \(\Weyl\), associated to a conical symplectic singularity \(\variety\).
Suppose $ Y \rightarrow X$ is a $\BQ$-factorial terminalisation of \(\variety\). The group \(\Weyl\) then acts as a reflection group on $H^2(Y,\Comp) $, which depends on \(\variety\) only.

Let $\leaf$ be a codimension $2$ symplectic leaf in \(\variety\) and let $x \in \leaf$. Then the formal neighbourhood of \(x\) in \(\variety\) is isomorphic to $\Comp<2n-2>\times \Comp<2>/ \ade$, where $\dim X=2n$ and $\ade \subset SL_2(\Comp)$ is a finite subgroup. Denote by $S_{\leaf}$ the corresponding Kleinian singularity, so by the McKay correspondence we can associate with $\ade$ a Weyl group $\Weyl[\leaf]$ of type $\A\D\E$ and consider the root space $\lie{\hat{h}}_{\leaf} \iso H^2(\widetilde{S}_{\leaf}, \Comp) $. The fundamental group $\uppi_1(\leaf)$ acts on $\Weyl[\leaf]$ and $\lie{\hat{h}}_{\leaf}$  by automorphisms of the Dynkin diagram. Denote by $\Weyl[\leaf]'$ the centralizer of $\uppi_1(\leaf)$ in $\Weyl[\leaf]$. Taking the product over all leaves we obtain Namikawa's Weyl group
\begin{equation}\label{eqn:nwgroup}
    \Weyl \defn \Prod{\Weyl[\leaf]'}[\leaf] 
\end{equation}
\begin{thm}\cite[Lem. 2.8]{losev_deformations_2022}
    There is a vector space isomorphism
    \begin{equation*}
        H^2(Y,\Comp) \iso H^2(X^{\on{reg}}, \Comp) \oplus \bigoplus_{\leaf} \lie{h}_{\leaf},  
    \end{equation*}
    where $\lie{h}_{\leaf} \defn \left (\lie{\hat{h}}_{\leaf} \right ) ^{\uppi_1(\leaf)}$.
\end{thm}

We would like to compute \eqref{eqn:nwgroup} when \(X\defn\QV{\root}\) and \(\p(\root)=2\), in which case \(Y\defn\QV{\root}<\stab>\) will be projective symplectic resolution for generic \(\stab\) (this is \cite[Cor. 1.3]{bellamy_birational_2023}). To do this, we need to find all the codimension $2$ leaves $\leaf$ of a quiver variety and compute the centraliser $\Weyl[\leaf]'$. For the first step, the  \cref{thm:isotropic_rep_type_are_codim2_leaves} explains that all codimension $2$ leaves of a quiver variety are given by the isotropic decompositions. The following theorem shows that the second step is not required when \(\p(\root)=2\).

\begin{thm}
    For each symplectic leaf $\leaf$ of codimension two, corresponding to a representation type $\rtype$ with a single imaginary root, we have $\Weyl[\leaf]'= \Weyl[\leaf]$.
\end{thm}
\begin{proof}
    Firstly, $H^2(\widetilde{S}_{\leaf}, \Comp) \iso \lie{\hat{h}}_{\leaf}$ is the space $\Comp<Q''>$, where \(\quiver''\) is the ext-quiver constructed in \cref{dfn:isotropic_decomposition}. Following \cite[Thm. 3.2, Thm. 3.3, and Ppn. 3.4]{wu_namikawa-weyl_2023} (see also \cite{mcgerty_kirwan_2018}, \cite{losev_isomorphisms_2012}), there is a surjection
    \begin{equation*}
        \upkappa\colon\Int<Q_0> \rightarrow H^2(\QV{\root}<\stab>, \Comp).
    \end{equation*}
    In particular, its restriction to each $\lie{h}_{\leaf} \subseteq \lie{\hat{h}}_{\leaf}$ can be obtained using \cite[Ppn. 3.4]{wu_namikawa-weyl_2023} (see also Thm. 1.6 and Rmk. 1.7 in loc. cit.). Explicitly, the restriction of the projection is given by 
    \begin{equation*}
        \upkappa_{\leaf}: \upchi \longmapsto (\upchi \cdot \root*<2>, \ldots, \upchi \cdot \root*<m>, \upchi \cdot \root**<1> , \ldots, \upchi \cdot \root**<k>) \in \lie{h}_{\leaf}.
    \end{equation*}
    Notice that whenever the representation type $\rtype$ only has a single imaginary root $\root*<1>$, the image of this map is given by $(\upchi \cdot \root**^{1} , \ldots, \upchi \cdot \root**^{k})$, which is a full-rank lattice in $\lie{\hat{h}}_{\leaf}$, as all the remaining roots are simple. Hence $\lie{\hat{h}}_{\leaf}=\lie{h}_{\leaf}$.
\end{proof}

\begin{cor}\label{cor:nwgroup}
    For $\p(\root)=2$ the \emph{Namikawa Weyl} group of \(\QV{\root}\) is the product of the Weyl groups corresponding to the Kleinian slice singularities across the isotropic decompositions of $\root$.
\end{cor}

Applying this to \cref{tbl:decomp}, we obtain a proof of \cref{thm:namikawa_weyl_groups_are_classified}.

\begin{rmk}
    Unfortunately, this simplification does not hold in higher dimensions. Indeed, \cite[Exm. 3.6, 3.7]{wu_namikawa-weyl_2023} gives counterexamples in 8d and 6d, respectively.
\end{rmk}

  \section{Arrangements and resolutions}\label{sec:varieties}

Now that we have computed all minimal degenerations between symplectic leaves, we may use this to compare the quiver varieties associated with \cref{tbl:dimvec} among themselves, and with products and symmetric powers of Kleinian singularities. To finish our classification we then compare the number of projective symplectic resolutions for the remaining cases, using a collection of hyperplanes known as the secondary arrangement.

\subsection{Products and symmetric powers of Kleinian singularities}

Recall \cref{thm:canonical_decomposition_and_dim_of_quiver_variety} tells us that other than those induced by \cref{tbl:dimvec}, the only possible quiver varieties in 4d are given by the products of lower-dimensional quiver varieties with symmetric powers of Kleinian singularities. In our case, the only such options are \(\Sym<2>(\Comp<2>/\ade)\) or \((\Comp<2>/\ade<1>)\times(\Comp<2>/\ade<2>)\). However, there could still be some \(\root\in\Sig\) such that \(\QV{\root}\) is isomorphic to such a variety. We start by investigating this possibility.

In \cite[Lem. 8.4]{bellamy_singularities_2025} the authors show that the quiver variety corresponding to $\II(\ade<1>, \ade<2>)$ is the product of Kleinian singularities of type $\ade<1>$ and $\ade<2>$. We are thus left with symmetric powers to study. By \cite[Thm. 3.4]{crawley-boevey_decomposition_2002}, $\Sym<n>(\QV{\root**})$ may be realised as a quiver variety with dimension vector \(n\root**\).

\begin{lem}\label{lem:rtypesympow}
    The representation types for $2\root**$ are $2\root**$, $ \root** +\root**$, $\root** +\sum \root**[i]\srt[i]$ and $\sum 2\root**[i]\srt[i]$.
\end{lem}
\begin{proof}
    This follows by simply observing the minimal imaginary subroots of \(2\root**\).
\end{proof}

In particular, as the quiver representations classify points in $\Sym<2>(\Comp<2>/\Gamma)$, these representation types label the symplectic leaves, as per \cref{thm:rep_type_are_leaves}. Explicitly, $\rtype = (1,\root**;1,\root**)$ are the generic points of $\leaf[4]$, as this is the canonical decomposition,  $\rtype=(2,\root**)$ is the diagonal leaf, $\rtype= (1,\root**; \dots;\root**[i],\srt[i];\dots)$ are points of the form $(x,0) \iso (0,x)$ and $\rtype= (\dots;2\root**[i],\srt[i];\dots)$ is \(\leaf[0]\).

\begin{thm}\label{thm:hassesympow}
    Hasse diagrams for $\Sym<2>(\Comp^2/\ade)$ have the following form (labeled by $\rtype$): 
    \begin{equation}\label{eqn:hassesympow}
        \begin{tikzcd}
        \&  \root** +\root** \ar[dl, "\ade"'] \ar[dr, "\mathsf{A}_1"] \\
        \root** +\sum \root**[i]\srt[i] \&\& 2\root** \\
        \&  \sum 2\root**[i]\srt[i] \ar[ul, "\ade"] \ar[ur, "\ade"']
        \end{tikzcd}
    \end{equation}
\end{thm}
\begin{proof}
    The degenerations between the 0d and 2d leaves are given by the closures of the corresponding leaves, which are both $\Comp<2>/ \Gamma$. For the upper degenerations we invoke a more general technique. 
    \cite[Thm. 8.1]{bellamy_singularities_2025} explains that the quiver variety corresponding to the ext-quiver $Q''$ from \cref{dfn:isotropic_decomposition} describes étale local slice to a point in the corresponding symplectic leaf.    
    The ext-quiver for the diagonal leaf is
    \begin{equation*}
        \begin{tikzcd}
            2 \arrow[loop left],
        \end{tikzcd}
    \end{equation*}
    so locally the variety is $\Sym<2>(\Comp<2>)\iso\Comp<2>\times \Comp<2>/\Cyc[2]$. The smooth part is the leaf direction and the étale slice is given by $\Comp<2> /\Cyc[2]$, meaning that the corresponding minimal degeneration is $\A[1]$.

    Next, we study the etale local slice to the leaf corresponding to $\root** +\sum \srt[i]$. The ext-quiver  
    \begin{equation*}
        \begin{tikzcd}
            1 \arrow[loop left] \&  \mir(\ade)
        \end{tikzcd}
    \end{equation*}
    is disconnected and consists of two parts, so locally we have the structure of $\Comp<2> \times \Comp<2> / \ade$. As before, the minimal degeneration corresponds to $\ade$.
\end{proof}

\begin{rmk}
    In the spirit of \cref{cor:nwgroup} one can note that by, for example, \cite{bellamy_birational_2020}, the Namikawa-Weyl group of $\Sym<2>(\Comp^2/\ade)$ is also computed by taking the product of groups, corresponding to the two codimension $2$ slice singularities. 
\end{rmk}

\begin{lem}\label{lem:symmetric_square_of_A1_is_not_same_as_product}
    The varieties $\Sym<2>(\Comp<2>/\Cyc[2]) \iso \Comp<4>/(S_2 \ltimes\Cyc[2]<2>)$ and $(\Comp<2>/\Cyc[2])\times(\Comp<2>/\Cyc[2])$ are non-isomorphic.
\end{lem}
\begin{proof}
    Write \(\ade[1]\defn\Cyc[2]<2>\), \(\ade[2]\defn\Sym[2]\ltimes\Cyc[2]<2>\), and consider the set \(\Sing[i]\subset\Comp<4>\) of points with a non-trivial stabiliser under the action of \(\ade[i]\). Both \(\ade[1]\) and \(\ade[2]\) act symplectically, so $\Sing[i]$ will have codimension \(2\) and real codimension \(4\). Hence by \cite[Thm. 2.3]{godbillon_elements_1971}, $\uppi_1(\Comp<4>\setminus\Sing[i]) \iso \uppi_1(\Comp<4>)$ is trivial. Using similar reasoning, the fundamental groups of the factors are isomorphic to \(\uppi_1((\Comp<4>/\ade[i])\setminus\Sing[i]\), but as the action of these groups is not free, the corresponding fundamental groups are \(\Cyc[2]<2> \neq \Sym[2] \ltimes\Cyc[2]<2>\).
\end{proof}

We can now rule out new symmetric powers from the classification of 4d quiver varieties.

\begin{thm} \label{thm:no_sym_products_in_table}
    There is no $\root \in \Sig$ with $\p(\root)=2$ such that \(\QV{\root}\iso\Sym<2>(\Comp<2>/\ade)\) for some \(\A\D\E\) \(\ade\).
\end{thm}
\begin{proof}
    As explained in \cref{ssec:mindegen}, \cref{tbl:quivvar} contains all the information about minimal degenerations of \(\QV{\root}\) in this case. A simple check shows that the only Hasse diagrams that coincide between these and ones of the form \eqref{eqn:hassesympow} are for $\Sym<2>(\Comp<2>/\Cyc[2]$ and $(\Comp<2>/\Cyc[2])\times(\Comp<2>/\Cyc[2])$. These two are different by \cref{lem:symmetric_square_of_A1_is_not_same_as_product}.
\end{proof}

\subsection{Distinguishing the quiver varieties}\label{ssec:compare}

We begin working towards the proof of \cref{thm:quivvar}.

\begin{thm} 
    \label{thm:most_quiver_varieties_are_different}
    All of the quiver varieties \(\QV{\root}\) with \(\root\) in \cref{tbl:quivvar} are non-isomorphic, with the possible exception of $\I(1,1,4)$ and $\II(\D[4])$.
\end{thm}
\begin{proof}
    We proceed as in the proof of \cref{thm:no_sym_products_in_table}. It is clear that, within a given family, the decompositions are all different. We proceed by the number of decompositions (equivalently, the number of intermediate leaves), enumerated below.
    \begin{enumerate}
        \setcounter{enumi}{-1}
        \item Observe that only \(\I(1,1,1)\) and \(\II(\A[0],\A[0])\) have no nontrivial decompositions.
        \item The cases where \(\mir\in\set{\D,\E}\) are all distinct, except from \(\II(5)\) and \(\I[a]\) which both have \(\mir=\D[5]\). However the former has \(\hrr=\D[7]\) and the latter has \(\hrr=\D[6]\), so these decompositions do not coincide. This leaves the cases where \(\mir=\A\), within which only \(\II(\A[0],\ade<2>)\) has \(\mir=\A[0]\) and only \(\I(\E[7], n)\) with \(n\in\set{0,1}\), has \(\hrr=\E[6]\). We see that \(\I(1,1,4)\) and \(\II(\D[4])\) are distinct dimension vectors with the same decomposition \((\mir,\hrr) = (\A[1],\A[3])\). The data also coincides for \(\I(\D[4],n)\) and \(\I(\D[4],n)'\) more generally, but this is because the dimension vectors themselves are equal.
        \item We observe that the construction \(\II(\ade<1>,\ade<2>)\) is unique in having \say{symmetric} decompositions, meaning that \(\mir<1>=\hrr<2>\) and \(\mir<2>=\hrr<1>\). The families \(\III(m,n)\) (\(m\geq5\)) can be separated from the rest by noting that they have one decomposition with \(\mir=\D\) and another with \(\mir=\E\), all of which are different to each other as well. A similar argument holds for \(\I(\E[6],n)\), \(\I(\E[7],n)\) and \(\I(\D[m],n)\), \(\I(\D[m],n)'\), leaving \(\I(1,m,n)\) that has unique decompositions.
        \item There are five such dimension vectors: \(\I(\ell,m,n)\) (\(1\leq*\ell\leq m\leq n\)) has only \(\A\) in its decompositions, \(\III(\D[m],i)\) (\(m\geq4,i\geq2\)) has \(\A\) and \(\D\) in its decompositions, and the three cases \(\III(m,n)\) (\(n\in\set{6,7,8}\)) all have a different \(\E[n]\) in their decompositions.
        \item The dimension vector \(\III(\D[m],1)\) (\(m\geq4\)) is unique in having four decompositions.
        \item Similarly to the previous case, only \(\III(\D[4],1)\) has five decompositions. \qedhere
    \end{enumerate}
    Since we know that \(\QV{\I(1,1,1)}=\CO_{min}\) is singular and \(\QV{\II(\A[0],\A[0])}=\Comp<4>\) is smooth, these quiver varieties are non-isomorphic. For $\I(1,1,4)$ and $\II(\D[4])$, which share the Hasse diagram below, we will need an additional invariant to distinguish their quiver varieties geometrically.
    \begin{equation}
        \setlength\arraycolsep{2.5em}
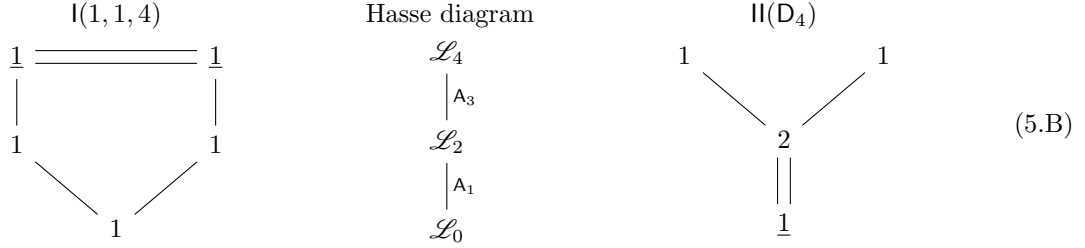
        \begin{array}{ccc}
            \I(1,1,4) & \text{Hasse diagram} & \II(\D[4]) \\
            \begin{tikzcd}
            \bal{1} \ar[rr, shift left] \ar[rr, shift right] \&\& \bal{1} \\
            1 \ar[u] \&\& 1 \ar[u] \\
            \& 1 \ar[ul] \ar[ur]
            \end{tikzcd} & \begin{tikzcd}
            \leaf[4] \ar[d, "{\A[3]}"] \\
            \leaf[2] \ar[d, "{\A[1]}"] \\
            \leaf[0]
            \end{tikzcd} & \begin{tikzcd}
            1 \ar[dr] \&\& 1 \ar[dl] \\
            \& 2 \ar[d, shift left] \ar[d, shift right] \\
            \& \bal{1}
            \end{tikzcd}
        \end{array}
    \end{equation}
\end{proof}

\subsection{Secondary arrangements and projective resolutions} \label{ssec:secondary_arrangements}

In this subsection we explain how to compute the number of projective symplectic resolutions for the quiver varieties from \cref{tbl:dimvec}. We then proceed to apply this, as a further invariant, to \(\QV{\II(\D[4])}\) and \(\QV{\I(1,1,4)}\).


For \(\root**\in\Root\) consider the hyperplane \(\root**^{\perp} \defn \set{\stab\in\on{Hom}(\BZ^{Q_0},\BQ)}[\stab(\root**)=0]\). This is a subset of \(\Stab\defn\root^\perp\), the space of GIT stability parameters.

\begin{dfn}[\cite{bellamy_birational_2023}, Dfn. 4.13]\label{dfn:secarr}
    The \emph{secondary arrangement} associated with \(\QV{\root}\) is given by
    \begin{equation*}
        \Hpln[\root] \defn \set{\root*^{\perp} \cap \root^{\perp}}[\root = \root* + (\root - \root*) \, \text{is a decomposition into two roots in } \Root] \subseteq \Stab.
    \end{equation*}
\end{dfn}

\begin{rmk}
    In general, the secondary arrangement is the hyperplane arrangement defining the Mori fan, which parametrises all crepant resolutions. In the case of symplectic singularities, symplectic and crepant resolutions are the same (\cite[Ppn. 3.2]{kaledin_crepant_2003}, see also \cite{fu_symplectic_2003}). Moreover, all our dimension vectors have at least one vertex of weight $1$. Hence their quiver varieties admit projective symplectic resolutions, as shown in \cite{bellamy_symplectic_2021}, so we will henceforth speak about these. Note that in general quiver varieties also admit proper, non-projective crepant resolutions. For a detailed discussion, see \cite{kaplan_nonprojective_2025}.
\end{rmk}

Denote by \(N(\root)\) the number of projective symplectic (equivalently, crepant) resolutions of \(\QV{\root}\). In \cite[\S 4.5]{bellamy_birational_2023} (see also Cor. 4.7 and Ppn. 4.12 loc. cit.) the authors explain that \(\Hpln[\root]\) is precisely the GIT fan. Moreover, in each of the GIT regions presented there is a chamber for the action of the Namikawa Weyl group $\Weyl[\root]$. Therefore the number of projective symplectic resolutions can be obtained from the secondary arrangement by enumerating the set \(\Cham\Hpln[\root]\) of chambers (connected components of the complement) and factoring out by the (transitive) action of $\Weyl[\root]$. In other words, if \(\root\) is indivisible there are
\begin{equation}\label{eqn:countres}
    N(\root)= \frac{|\Cham\Hpln[\root]|}{|\Weyl[\root]|}
\end{equation}
projective symplectic resolutions of \(\QV{\root}\).

\begin{rmk}\label{rmk:coxarr}
    Notice the secondary arrangement contains the Coxeter arrangement \(\Hpln*[\root]\) associated with the Namikawa Weyl group, which has hyperplanes corresponding to real roots. If this inclusion is strict, then \(\QV{\root}\) has at least two projective symplectic resolutions.
\end{rmk}

We give an explicit example which will be useful later.

\begin{exm}\label{exm:secarr}
    We compute the secondary arrangement for \(\I(1,1,3)\) and thus determine the number of projective symplectic resolutions admitted by the associated quiver variety. Using the indexing
    
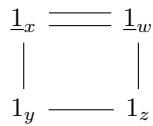
\begin{equation*}
        \begin{tikzcd}
            \bal{1}_x \ar[r, shift left] \ar[r, shift right] \ar[d] \& \bal{1}_w \ar[d] \\
            1_y \ar[r] \& 1_z
        \end{tikzcd}
    \end{equation*}
    we have \(\Hpln[\root]\subset\Stab=\set{w+x+y+z=0}\). Recall from \cref{ssec:root_combinatorics} that \(\Root\) is defined with respect to \(\Weyl[\graph]=\langle s_w,s_x,s_y,s_z\rangle\) with involutive generators, \(s_w\) and \(s_x\) free, and all other pairs having a braid relation of length three. The isotropic decomposition \(\root=\mir(\A[1])+\hrr(\A[2])\) gives the hyperplane \(\set{w+x=0}\), (equivalently, \(\set{y+z=0}\) in \(\Stab\)). We also have the hyperplane \(\set{w+z=0}\) arising from the decomposition \(\root=\hrr(\A[2])+\hrr(\A[2])\) into real roots. Similarly there are two decompositions of the form \(\root=\hrr(\A[3])+\hrr(\A[1])\), giving hyperplanes \(\set{w=0}\) and \(\set{x=0}\). Finally observe that the subroot \(\root*\defn(\begin{tikzcd}[cramped]
            1 \ar[r] \& 1 \ar[r, shift left] \ar[r, shift right] \& 1
        \end{tikzcd})\not\in\Sig\)
    is imaginary, being the image of the fundamental imaginary subroot \(\mir(\A[1])\) under the reflection \(s_y\) or \(s_z\). This gives hyperplanes \(\set{y=0}\) and \(\set{z=0}\). Hence \(\Hpln[\root]\) consists of \(6\) hyperplanes, cutting \(\Real<3>\) into \(24\) regions\footnote{An interactive version of \(\Hpln[\I(1,1,3)]\) is available: \href{https://www.desmos.com/3d/m5yp6auo7r}{https://www.desmos.com/3d/m5yp6auo7r}.}. The number of projective symplectic resolutions is then
    \begin{equation*}
       N(\root)=\frac{|\Cham\Hpln[\root]|}{|\Weyl_{\root}|} = \frac{24}{6} = 4,
    \end{equation*}
    using \cref{tbl:quivvar} for the denominator.
\end{exm}

This is the last invariant needed to distinguish all 4d quiver varieties.

\begin{cor} \label{cor:two_last_different_varieties}
    As affine schemes, \(\QV{\II(\D[4])}\) and \(\QV{\I(1,1,4)}\) are non-isomorphic.
\end{cor}
\begin{proof}
    We show that the former has a unique projective crepant resolution whilst the latter does not, so they cannot be isomorphic. Applying \cref{cor:nwgroup} to \cref{tbl:decomp}, we see that the two dimension vectors share a Namikawa Weyl group of type \(\A[3]\) (\(\Sym[4]\)), which has six hyperplanes in the associated Coxeter arrangement. We now count hyperplanes in each secondary arrangement by considering decompositions into two positive roots, starting with \(\II(\D[4])\). In the spirit of \cref{ssec:decompositions}, we see
    \begin{equation}\label{eqn:iid4decomp}
        \ubrace*{\begin{tikzcd}[scale cd=.8, cramped, sep=0]
            \&\&\&1 \\ \bal{1}=2 \\ \&\&\&1
        \end{tikzcd}}{\II(\D[4])} \quad=\quad \ubrace*{\begin{tikzcd}[scale cd=.8, cramped, sep=0]
            \&\&\&1 \\ 1=1 \\ \&\&\&0
        \end{tikzcd}}{\root*<1>} \quad+\quad \ubrace*{\begin{tikzcd}[scale cd=.8, cramped, sep=0]
            \&\&\&0 \\ 0=1 \\ \&\&\&1
        \end{tikzcd}}{\hrr(\A[2])} \quad=\quad \ubrace*{\begin{tikzcd}[scale cd=.8, cramped, sep=0]
            \&\&\&1 \\ 1=2 \\ \&\&\&0
        \end{tikzcd}}{\root*<2>} \quad+\quad \ubrace*{\begin{tikzcd}[scale cd=.8, cramped, sep=0]
            \&\&\&0 \\ 0=0 \\ \&\&\&1
        \end{tikzcd}}{\hrr(\A[1])} \quad=\quad \ubrace*{\begin{tikzcd}[scale cd=.8, cramped, sep=0]
            \&\&\&1 \\ 1=1 \\ \&\&\&1
        \end{tikzcd}}{\root*<3>} \quad+\quad \ubrace*{\begin{tikzcd}[scale cd=.8, cramped, sep=0]
            \&\&\&0 \\ 0=1 \\ \&\&\&0
        \end{tikzcd}}{\hrr(\A[1])}
    \end{equation}
    the first two decompositions each occurring in two different ways due to symmetry. There is also a single isotropic decomposition since \(\Weyl[\II(\D[4])]\) has one factor. The dimension vectors \(\root*<1>,\root*<2>,\root*<3>\) are readily checked to be roots, as they all lie in the image of the fundamental imaginary root \(\mir(\A[1])\) under simple reflections. This gives at least six hyperplanes in \(\Hpln[\II(\D[4])]\). To see that there are no other permitted decompositions of \(\II(\D[4])\) into positive roots, denote the branch vertex in some decomposition by \(i\) and the first summand by \(\root*\). Then as roots necessarily have connected support, \(i\) must have at least two neighbours with nonzero weight whenever \(\root*[i]\neq0\). In \eqref{eqn:iid4decomp}, the first equality gives the cases where \(\mir(\A[1])\leq*\root*\) and \(\root*[i]=1\). The cases where \(\mir(\A[1])\leq*\root*\) and \(\root*[i]=2\) are the second equality. The case where every neighbour has nonzero weight and \(\root*[i]=1\) is the third equality, and the isotropic decomposition is the case where \(\mir(\A[1])\not\leq*\root*\) and \(\root[i]=1\). This leaves the case where \(\mir(\A[1])\not\leq*\root*\) and \(\root[i]=2\), that is
    \begin{equation*}
        \ubrace*{\begin{tikzcd}[scale cd=.8, cramped, sep=0]
            \&\&\&1 \\ \bal{1}=2 \\ \&\&\&1
        \end{tikzcd}}{\II(\D[4])} \quad=\quad \ubrace*{\begin{tikzcd}[scale cd=.8, cramped, sep=0]
            \&\&\&1 \\ 0=1 \\ \&\&\&1
        \end{tikzcd}}{\root} \quad+\quad \ubrace*{\begin{tikzcd}[scale cd=.8, cramped, sep=0]
            \&\&\&0 \\ 1=0 \\ \&\&\&0
        \end{tikzcd}}{\hrr(\A[1])}
    \end{equation*}
    but a quick computation gives \(\p(\root*)=-1\), so \(\root*\not\in\Root\). Hence the Coxeter arrangement and secondary arrangement associated with \(\II(\D[4])\) coincide, so \(\QV{\II(\D[4])}\) has a unique projective crepant resolution by \eqref{eqn:countres}. It now suffices to prove that \(|\Hpln[\I(1,1,4)]|\geq*6\), as then \(N(\I(1,1,4))\geq*1\) thus concluding the proof. To this end, consider the decompositions
    \begin{equation*}
        \ubrace*{\begin{tikzcd}[scale cd=.8, cramped, sep=0]
            \&\bal{1}=\bal{1} \\\\ 1\&1\&1
        \end{tikzcd}}{\I(1,1,4)} \quad=\quad \ubrace*{\begin{tikzcd}[scale cd=.8, cramped, sep=0]
            \&1=0 \\\\ 1\&1\&1
        \end{tikzcd}}{\hrr(\A[4])} \quad+\quad \ubrace*{\begin{tikzcd}[scale cd=.8, cramped, sep=0]
            \&0=1 \\\\ 0\&0\&0
        \end{tikzcd}}{\hrr(\A[1])} \quad=\quad \ubrace*{\begin{tikzcd}[scale cd=.8, cramped, sep=0]
            \&1=0 \\\\ 1\&0\&0
        \end{tikzcd}}{\hrr(\A[3])} \quad+\quad \ubrace*{\begin{tikzcd}[scale cd=.8, cramped, sep=0]
            \&0=1 \\\\ 0\&1\&1
        \end{tikzcd}}{\hrr(\A[2])} \quad=\quad \ubrace*{\begin{tikzcd}[scale cd=.8, cramped, sep=0]
            \&1=1 \\\\ 1\&0\&0
        \end{tikzcd}}{\root*} \quad+\quad \ubrace*{\begin{tikzcd}[scale cd=.8, cramped, sep=0]
            \&0=0 \\\\ 0\&1\&1
        \end{tikzcd}}{\hrr(\A[1])}
    \end{equation*}
    where \(\root*=\root*<1>\in\Root\) from \eqref{eqn:iid4decomp}. Since each decomposition occurs in two different ways, if we include the isotropic decomposition from \cref{sssec:typei} then this is at least seven hyperplanes.
 \end{proof}


We can now prove our main result.

\begin{proof}[Proof of \cref{thm:quivvar}]
    Combining \cref{thm:most_quiver_varieties_are_different} and \cref{cor:two_last_different_varieties} we prove \cref{thm:quivvar}. Including \cref{thm:no_sym_products_in_table} completes the classification.
\end{proof}

\subsection{The group quotient G4}\label{ssec:gfour}

The rank-two complex reflection group \(G_4\) (classified in \cite{shephard_finite_1954}) is the binary tetrahedral group of order \(24\). It can be realised as a finite subgroup
\begin{equation*}
    G_4 = \left\{ \pm 1, \pm i, \pm j, \pm k, \tfrac{1}{2}(\pm 1\pm i\pm j \pm k)\right\}
\end{equation*}
of the group of units in the quaternions. It is generated by the elements $x=\frac{1}{2}(-1+i+j-k)$ and $y=\frac{1}{2}(-1+i-j+k)$ and has presentation as \(\gen{x, y}[x^3=y^3=(xy)^6=1]\). 

In \cite[Exm. 5.6]{bulois_towards_2017}, the 4d conical symplectic singularity \(\Comp<4>/G_4\) was realised as a Hamiltonian reduction by a semisimple group. This led the authors in \cite[\S 1.5]{bellamy_all_2024} to ask whether \(\Comp<4>/G_4\) admits a construction as a quiver variety. Applying our classification, we show that this is not the case.

\begin{proof}[Proof of \cref{cor:g4}]
    Following \cite[\S 7.2]{bellamy_hyperplane_2018}, \(\Comp<4>/G_4\) has Namikawa Weyl group \(\Sym[3]\), and admits two projective symplectic resolutions. Thus for \(\Comp<4>/G_4\) to be a quiver variety it must be isomorphic to either \(\QV{\II(\A[0],\A[2])}\) or \(\QV{\I(1,1,3)}\), as by \cref{thm:quivvar} these are the only 4d quiver varieties with Weyl group \(\Sym[3]\). However, the former has a unique projective symplectic resolution (being the direct product of Kleinian singularities), and the latter has four projective symplectic resolutions, as seen in \cref{exm:secarr}. Hence \(\Comp<4>/G_4\) is a 4d conical symplectic singularity that cannot be realised as a Nakajima quiver variety.

    Moreover, we can show that neither projective symplectic resolution $Y$ of \(X\defn\Comp<4>/G_4\) is a quiver variety. Suppose that \(\variety*\) is isomorphic to \(\QV[\stab]{\root}\) for some dimension vector \(\root\) and $\stab\in\Stab$. This isomorphism then gives rise to a map
    \(\uppi\colon Y\to \QV{\root}\). If \(\uppi\) is a resolution of singularities, then both \(X\) and \(\QV{\root}\) are isomorphic to \(\Spec\Upgamma(Y,\mathcal{O})\), but this contradicts the above argument. If \(\uppi\) is not a resolution of singularities, then by \cite[Thm. A.1]{bellamy_birational_2020}, \(\QV[\stab]{\root}\) resolves the normalisation of a leaf closure \(\cl{\leaf}\). By \cite[Thm. 1.9]{bellamy_singularities_2025}, the latter is isomorphic to a product of symmetric powers of quiver varieties. Hence \(X \iso \Spec\Upgamma(Y,\mathcal{O}) \iso \cl{\leaf}\), which means, by the argument in the previous paragraph, that 
    \(X=(\Comp<2>/\Upgamma_1\times\Comp<2>/\Upgamma_2)=\Comp<4>/(\Upgamma_1\times\Upgamma_2)\) or \(X=S^2(\Comp<2>/\Upgamma)=\Comp<4>/(S_2\wr\Upgamma)\). But \(G_4\) is an irreducible complex reflection group, so cannot take either of these forms.
\end{proof}

\subsection{Counting resolutions}\label{ssec:reso}

To conclude our work we provide computer-aided\footnote{\texttt{Magma} and \texttt{Julia} scripts used are available from the authors upon request.} numerology regarding the number of projective symplectic resolutions admitted by 4d quiver varieties. By \cite[Thm. 1.5]{bellamy_all_2024}, all admit at least one projective symplectic resolution, as every \(\root\) in \cref{tbl:dimvec} is indivisible. Using the computer algebra package \texttt{Magma}, \cite{bosma_magma_1997}, we compute the hyperplanes present in the secondary arrangement \(\Hpln[\root]\) by iterating over subroots \(\root*\leq*\root\), removing those that necessarily cannot be roots (for example \(\p(\root*)\leq*0\) or \(\supp{\root*} \) disconnected) and checking the Weyl orbits of the remaining vectors for an element of the relevant fundamental set. To compute the number of chambers in \(\Hpln[\root]\) we then pass our hyperplanes to the \texttt{Julia} package \texttt{CountingChambers}, developed in \cite{brysiewicz_computing_2023}. This process finishes in a reasonable time up to approximately \(\dim\Stab=10\), giving suggestions for the behaviour of \(N(\root)\) for many of the infinite families in our classification.

\begin{con}\label{con:typei}
    In Type \(\I\) we have the following counts of projective symplectic resolutions:
    \begin{enumerate}
        \item \(N(\I(\D[4],n))=\tfrac{n+3}{3}\tbinom{n+5}{5}\).
        \item \(N(\I(\D[m],1))=2^{m-1}\).
        \item \(N(\I(\D[m],n)') = \tfrac{2(m+n-1)}{n!}\Prod{(2m+i-2)}[i=1]<n-1>\) and \(N(\I(\D[m],0)') = 1\).
        \item \(N(\I(\E[6,7],n))\) coincides with \cite[A030648, A030649]{oeis_foundation_inc_-line_2025}, the dimensions of multiples of the \(\E[6,7]\) complex Lie algebra's minimal representation. 
        \item \(N(\I(m)) = \Prod{\tfrac{(2i)!(m+2i+1)!}{(3i+1)!(m+i)!}}[i=0]<m-1>\), which by \cite[Cor. 5.2]{stembridge_enumeration_1995} coincides with \cite[A005157]{oeis_foundation_inc_-line_2025}, the number of totally symmetric plane partitions that fit in an \(m\times m\times m\), \(m\geq4\).
    \end{enumerate}
\end{con}

Our second conjecture is motivated by observing that in each computed case, the number of chambers in the secondary arrangement of \(\root\) coincided with the size of its Namikawa Weyl group in \cref{tbl:quivvar}. This forces \(\Hpln[\root]=\Hpln*[\root]\) and a unique resolution for the quiver variety, following \cref{rmk:coxarr}.

\begin{con}\label{con:typeii}
    If \(\root\in\Sig\) is a dimension vector of Type \(\II\) then \(N(\root)=1\).
\end{con}

Clearly this holds for \(\II(\root<1>,\root<2>)\), its associated quiver variety being the product of Kleinian singularities. Similarly for \(\II(\A[m])\), which is just \(\II(\A[d-1],\A[m-d])\) for some \(d\). However we have also observed a \say{trivial} secondary arrangement for \(\II(\E[6])\), \(\II(\E[7])\), \(\II[a]\), as well as for \(\II(\D[m])\), \(\II(\D[m])'\), and \(\II(m)\) for low \(m\geq4\). We believe that more abstractly it should be possible to rule out the existence of special root decompositions (that is, outside the Coxeter arrangement) for \(\root\) of Type \(\II\).

  \appendix
    
  \section{Hasse diagrams}\label{app:hasse}
  For ease of reference, this appendix contains the Hasse diagram for each 4d quiver variety. By the theory in 
\cref{sec:symplectic}, the below tables are simply a visualisation of \cref{tbl:decomp}.

\rowcolors{1}{}{}
\begin{longtable}{|C|C|C|C|C|}\caption{Minimal degenerations of \(\QV{\root}\) for \(\root\in\Sig\) of Type \(\I\)}\label{tbl:hasseI}\\\hline
    \rowcolor{gray!50} \phantom{mm}\root\phantom{mm} & \multicolumn{4}{C|}{\phantom{mmmmmmmmmmmmmmmm}\text{Hasse diagram(s)}\phantom{mmmmmmmmmmmmmmmm}} \\\hline\endfirsthead
    \caption[]{(continued)}\\\hline\rowcolor{gray!50} \root & \multicolumn{4}{C|}{\text{Hasse diagram(s)}} \\\hline\endhead
    \hline\endfoot
    \mrow{\I(\ell,m,n) \\[3.2em] \begin{tikzcd}[scale cd=.8]
        \bal{1} \ar[r, bend left, "\cdots" description] \ar[r, "\cdots" description] \ar[r, bend right, "\cdots" description] \& \bal{1}
    \end{tikzcd} \\[3.2em] 1\leq\ell\leq m\leq n} & \multicolumn{1}{C}{\mrow{\begin{tikzcd}[column sep=2.5em]
        \& \leaf[4] \ar[dl, "{\A[\ell-1]}"'] \ar[d, "{\A[m-1]}"] \ar[dr, "{\A[n-1]}"] \\
        \leaf[21] \& \leaf[22] \& \leaf[23] \\
        \& \leaf[0] \ar[ul, "{\A[m+n-1]}"] \ar[u, "{\A[\ell+n-1]}"'] \ar[ur, "{\A[\ell+m-1]}"']
    \end{tikzcd}\\\\\ell\geq2}} & \multicolumn{1}{C}{\mrow{\begin{tikzcd}
        \& \leaf[4] \ar[dl, "{\A[m-1]}"'] \ar[dr, "{\A[n-1]}"] \\
        \leaf[21] \&\& \leaf[22] \\
        \& \leaf[0] \ar[ul, "{\A[n]}"] \ar[ur, "{\A[m]}"']
    \end{tikzcd}\\\\\ell=1,m\geq2}} & \multicolumn{1}{C}{\mrow{\begin{tikzcd}
        \leaf[4] \ar[d, "{\A[n-1]}", overlay] \\
        \leaf[2] \\
        \leaf[0] \ar[u, "{\A[1]}"', overlay]
    \end{tikzcd}\\\\m=1,n\geq2}} & \mrow{\begin{tikzcd}
        \leaf[4] \ar[dd, "{a_1}"] \\
        \phantom{\leaf[2]} \\
        \leaf[0]
    \end{tikzcd}\\\\n=1} \\
    \mrow[gray!20]{\I(\D[m],n) \\[3.2em] \begin{tikzcd}[scale cd=.8, cramped, sep=0, background color=gray!20]
    1\&\&\&\&1 \\ \&2\&\cdots\&2 \\ \bal{1} \ar[rrrr, bend right, "\cdots" description]\&\&\&\&\bal{1} \end{tikzcd} \\[3.2em] m\geq4, n\geq0} & \multicolumn{2}{C}{
        \mrow[gray!20]{\begin{tikzcd}
        \& \leaf[4] \ar[dl, "{\A[n-1]}"'] \ar[dr, "{\A[m-1]}"] \\
        \leaf[21] \&\& \leaf[22] \\
        \& \leaf[0] \ar[ul, "{\D[m]}"] \ar[ur, "{\A[m+n-3]}"']
    \end{tikzcd} \\\\ n\geq2}
    } & \multicolumn{2}{C|}{\mrow[gray!20]{\begin{tikzcd}
        \leaf[4] \ar[d, "{\A[m-1]}", overlay] \\
        \leaf[2] \\
        \leaf[0] \ar[u, "{\A[m+n-3]}"', overlay]
    \end{tikzcd}\\\\n\in\set{0,1}}} \\
    \mrow{\I(\D[m],n)' \\[3.2em] \begin{tikzcd}[scale cd=.8, cramped, sep=0]
    1\&\&\&\&\bal{1}\ar[dd, bend left = 5em, "\cdots" description, sloped] \\ \&2\&\cdots\&2 \\ 1\&\&\&\&\bal{1} \end{tikzcd} \\[3.2em] m\geq4, n\geq0} & \multicolumn{2}{C}{
        \mrow{\begin{tikzcd}
        \& \leaf[4] \ar[dl, "{\A[n-1]}"'] \ar[dr, "{\D[m-1]}"] \\
        \leaf[21] \&\& \leaf[22] \\
        \& \leaf[0] \ar[ul, "{\D[m]}"] \ar[ur, "{\A[n+1]}"']
    \end{tikzcd} \\\\ n\geq2}
    } & \multicolumn{2}{C|}{\mrow{\begin{tikzcd}
        \leaf[4] \ar[d, "{\D[m-1]}", overlay] \\
        \leaf[2] \\
        \leaf[0] \ar[u, "{\A[n+1]}"', overlay]
    \end{tikzcd}\\\\n\in\set{0,1}}} \\
    \mrow[gray!20]{\I(\E[6],n) \\[3em] \begin{tikzcd}[scale cd=.8, cramped, sep=0, background color=gray!20]
    \&\&1 \\ \&\&2 \\ \bal{1} \ar[rrrr, bend right, "\cdots" description] \&2\&3\&2\&\bal{1} \end{tikzcd} \\[3.3em] n\geq0} & \multicolumn{2}{C}{
        \mrow[gray!20]{\begin{tikzcd}
        \& \leaf[4] \ar[dl, "{\A[n-1]}"'] \ar[dr, "{\D[5]}"] \\
        \leaf[21] \&\& \leaf[22] \\
        \& \leaf[0] \ar[ul, "{\E[6]}"] \ar[ur, "{\A[n+3]}"']
    \end{tikzcd} \\\\ n\geq2}
    } & \multicolumn{2}{C|}{\mrow[gray!20]{\begin{tikzcd}
        \leaf[4] \ar[d, "{\D[5]}", overlay] \\
        \leaf[2] \\
        \leaf[0] \ar[u, "{\A[n+3]}"', overlay]
    \end{tikzcd}\\\\n\in\set{0,1}}} \\
    \mrow{\I(\E[7],n) \\[3em] \begin{tikzcd}[scale cd=.8, cramped, sep=0] \&\&\&2 \\ \bal{1}\ar[rrrrrr, bend right, "\cdots" description] \&2\&3\&4\&3\&2\&\bal{1}\end{tikzcd} \\[3.3em] n\geq0} & \multicolumn{2}{C}{
        \mrow{\begin{tikzcd}
        \& \leaf[4] \ar[dl, "{\A[n-1]}"'] \ar[dr, "{\E[6]}"] \\
        \leaf[21] \&\& \leaf[22] \\
        \& \leaf[0] \ar[ul, "{\E[7]}"] \ar[ur, "{\A[n+5]}"']
    \end{tikzcd} \\\\ n\geq2}
    } & \multicolumn{2}{C|}{\mrow{\begin{tikzcd}
        \leaf[4] \ar[d, "{\A[m-1]}", overlay] \\
        \leaf[2] \\
        \leaf[0] \ar[u, "{\A[n+5]}"', overlay]
    \end{tikzcd}\\\\n\in\set{0,1}}} \\
    \mrow[gray!20]{\I(m) \\[2em] \begin{tikzcd}[scale cd=.8, cramped, sep=0]
    \&\&\bal{1}\&\&\bal{1} \\ 1\&2\&3\&\cdots\&3\&2\&1
    \end{tikzcd} \\[2em] m\geq4} & \multicolumn{4}{C|}{\cellcolor{gray!20}\begin{tikzcd}
        \leaf[4] \ar[d, "{\A[m+1]}", overlay] \\
        \leaf[2] \\
        \leaf[0] \ar[u, "{\D[m]}"', overlay]
    \end{tikzcd}} \\
    \mrow{\I[a] \\\\ \begin{tikzcd}[scale cd=.8, cramped, sep=0]
    \bal{1}\&\&2 \\ \&3\&4\&3\&2\&1 \\ \bal{1}
    \end{tikzcd}} & \multicolumn{4}{C|}{\begin{tikzcd}
        \leaf[4] \ar[d, "{\D[6]}", overlay] \\
        \leaf[2] \\
        \leaf[0] \ar[u, "{\D[5]}"', overlay]
    \end{tikzcd}} \\
    \mrow[gray!20]{\I[b] \\\\ \begin{tikzcd}[scale cd=.8, cramped, sep=0]
    \bal{1}\&3 \\ \&\&5\&4\&3\&2\&1 \\ \bal{1}\&3
    \end{tikzcd}} & \multicolumn{4}{C|}{\cellcolor{gray!20}\begin{tikzcd}
        \leaf[4] \ar[d, "{\D[7]}", overlay] \\
        \leaf[2] \\
        \leaf[0] \ar[u, "{\E[6]}"', overlay]
    \end{tikzcd}} \\
    \mrow{\I[c] \\\\ \begin{tikzcd}[scale cd=.8, cramped, sep=0]
    \bal{1}\&\&\&\&3 \\ \&3\&4\&5\&6\&4\&2 \\ \bal{1}
    \end{tikzcd}} & \multicolumn{4}{C|}{\begin{tikzcd}
        \leaf[4] \ar[d, "{\E[7]}", overlay] \\
        \leaf[2] \\
        \leaf[0] \ar[u, "{\D[7]}"', overlay]
    \end{tikzcd}} \\
\end{longtable}

\begin{longtable}{|C|C|C|C|}\caption{Minimal degenerations of \(\QV{\root}\) for \(\root\in\Sig\) of Type \(\II\)}\label{tbl:hasseII}\\\hline
    \rowcolor{gray!50} \phantom{mm}\root\phantom{mm} & \multicolumn{3}{C|}{\phantom{mmmmmmmmmmmmmmmm}\text{Hasse diagram(s)}\phantom{mmmmmmmmmmmmmmmm}} \\\hline\endfirsthead
    \caption[]{(continued)}\\\hline\rowcolor{gray!50} \root & \multicolumn{3}{C|}{\text{Hasse diagram(s)}} \\\hline\endhead
    \hline\endfoot
    \mrow{\II(\ade<1>, \ade<2>) \\[3.4em] \begin{tikzcd}[scale cd=.8, cramped, sep=0]
        \hrr(\ade<1>)\&\bal{1}\&\hrr(\ade<2>) \end{tikzcd} \\[3.4em] \ade<i>\in\set{\A,\D,\E}} & \multicolumn{1}{C}{\mrow{\begin{tikzcd}
        \& \leaf[4] \ar[dl, "{\ade<1>}"'] \ar[dr, "{\ade<2>}"] \\
        \leaf[21] \&\& \leaf[22] \\
        \& \leaf[0] \ar[ul, "{\ade<2>}"] \ar[ur, "{\ade<1>}"']
    \end{tikzcd}\\\\\ade<1>\neq\A[0]\neq\ade<2>}} & \multicolumn{1}{C}{\mrow{\begin{tikzcd}
        \leaf[4] \ar[d, "{\ade<2>}", overlay] \\
        \leaf[2] \\
        \leaf[0] \ar[u, "{\A[0]}"', overlay]
    \end{tikzcd}\\\\\ade<1>=\A[0]\neq\ade<2>}} & \multicolumn{1}{C|}{\mrow{\begin{tikzcd}
        \leaf[4] \ar[dd] \\
        \phantom{\leaf[2]} \\
        \leaf[0]
    \end{tikzcd}\\\\\ade<1>=\A[0]=\ade<2>}} \\
    \mrow[gray!20]{\II(m) \\[2em] \begin{tikzcd}[scale cd=.8, cramped, sep=0]
        2\&\&\&\bal{1} \\ \&4\&\cdots\&4\&3\&2\&1 \\ 2 \end{tikzcd} \\[2em] m\geq4} & \multicolumn{3}{C|}{\cellcolor{gray!20}\begin{tikzcd}
        \leaf[4] \ar[d, "{\D[m+2]}", overlay] \\
        \leaf[2] \\
        \leaf[0] \ar[u, "{\D[m]}"', overlay]
    \end{tikzcd}} \\
    \mrow{\II[a]\\\\\begin{tikzcd}[scale cd=.8, cramped, sep=0]
    \&2\&\&3 \\ \bal{1}\&4\&5\&6\&4\&2 \end{tikzcd}} & \multicolumn{3}{C|}{\begin{tikzcd}
        \leaf[4] \ar[d, "{\E[7]}", overlay] \\
        \leaf[2] \\
        \leaf[0] \ar[u, "{\D[6]}"', overlay]
    \end{tikzcd}} \\
    \mrow[gray!20]{\II[b]\\\\\begin{tikzcd}[scale cd=.8, cramped, sep=0]
        \&\&\&5 \\ \bal{1}\&4\&7\&10\&8\&6\&4\&2
    \end{tikzcd}} & \multicolumn{3}{C|}
    {\cellcolor{gray!20}\begin{tikzcd}
        \leaf[4] \ar[d, "{\E[8]}", overlay] \\
        \leaf[2] \\
        \leaf[0] \ar[u, "{\E[7]}"', overlay]
    \end{tikzcd}}
\end{longtable}

\begin{longtable}{|C|C|C|C|C|}\caption{Minimal degenerations of \(\QV{\root}\) for \(\root\in\Sig\) of Type \(\III\)}\label{tbl:hasseIII} \\\hline
    \rowcolor{gray!50} \root & \multicolumn{4}{C|}{\phantom{mmmmmmmmmmmmmmm}\text{Hasse diagram(s)}\phantom{mmmmmmmmmmmmmmm}} \\\hline\endfirsthead
    \caption[]{(continued)}\\\hline\rowcolor{gray!50} \root & \multicolumn{4}{C|}{\text{Hasse diagram(s)}} \\\hline\endhead
    \hline\endfoot
    \mrow{\III(\D[m],i) \\[3.2em] \begin{tikzcd}[scale cd=.8, cramped, sep=0]
    1\&\&\&1\&\&\&1 \\ \&2\&\cdots\&\bal{2}\&\cdots\&2 \\ 1\&\&\&\&\&\&1 \end{tikzcd} \\[3.2em] 1\leq i\leq\ceil{\tfrac{m-3}{2}}} & \multicolumn{1}{C}{\mrow{\begin{tikzcd}[column sep=1.5em]
        \& \leaf[4] \ar[d, "{\A[1]}"] \ar[dl, "{\D[i+1]}"'] \ar[dr, "{\D[m-i-1]}"] \\
        \leaf[21] \& \leaf[22] \& \leaf[23] \\
        \& \leaf[0] \ar[u, "{\D[m]}"'] \ar[ul, "{\D[m-i+1]}"] \ar[ur, "{\D[i+3]}"']
    \end{tikzcd}\\\\i\neq1}} & \multicolumn{2}{C}{\mrow{\begin{tikzcd}[column sep=.5em]
        \&\& \leaf[4] \ar[dll, "{\A[1]}"'] \ar[dl, "{\A[1]}"] \ar[dr, "{\A[1]}"'] \ar[drr, "{\D[m-2]}"] \\
        \leaf[21] \& \leaf[22] \&\& \leaf[23] \& \leaf[24] \\
        \&\& \leaf[0] \ar[ull, "{\D[m]}"] \ar[ul, "{\D[m]}"'] \ar[ur, "{\D[m]}"] \ar[urr, "{\D[4]}"']
    \end{tikzcd}\\\\i=1}} & \mrow{\begin{tikzcd}
        \leaf[4] \ar[d] \ar[d, shift left = 1em, "{\A[1]}"] \ar[d, shift left = .5em] \ar[d, shift right = .5em] \ar[d, shift right =1em] \\
        \leaf[21]\cdots\leaf[25] \\
        \leaf[0] \ar[u] \ar[u, shift left = 1em] \ar[u, shift left = .5em] \ar[u, shift right = .5em] \ar[u, shift right =1em, "{\D[4]}"']
    \end{tikzcd}\\\\m=4} \\
    \mrow[gray!20]{\III(\E[7])\\\\\begin{tikzcd}[scale cd=.8, cramped, sep=0]
    \&\&\&1 \\ \&\&\&\bal{2} \\ 1\&2\&3\&4\&3\&2\&1 \end{tikzcd}} & \multicolumn{4}{C|}{\cellcolor{gray!20}\begin{tikzcd}
        \& \leaf[4] \ar[dl, "{\A[7]}"'] \ar[dr, "{\A[1]}"] \\
        \leaf[21] \&\& \leaf[22] \\
        \& \leaf[0] \ar[ul, "{\E[6]}"] \ar[ur, "{\E[7]}"']
    \end{tikzcd}} \\
    \mrow{\III(\E[8])\\\\\begin{tikzcd}[scale cd=.8, cramped, sep=0]
    \&\&\&\&\&3 \\ 1\&2\&3\&4\&5\&6\&4\&\bal{2}\&1 \end{tikzcd}} & \multicolumn{4}{C|}{\begin{tikzcd}
        \& \leaf[4] \ar[dl, "{\D[8]}"'] \ar[dr, "{\A[1]}"] \\
        \leaf[21] \&\& \leaf[22] \\
        \& \leaf[0] \ar[ul, "{\E[7]}"] \ar[ur, "{\E[8]}"']
    \end{tikzcd}} \\
    \mrow[gray!20]{\III(m,6) \\[3em] \begin{tikzcd}[scale cd=.8, cramped, sep=0]
    \&\&\&\&1 \\ 1\&\&\&\&2 \\ \&2\&\cdots\&\bal{2}\&3\&2\&1 \\ 1 \end{tikzcd} \\[3.4em] m\geq4} & \multicolumn{2}{C}{\mrow[gray!20]{\begin{tikzcd}
        \& \leaf[4] \ar[dl, "{\A[5]}"'] \ar[dr, "{\D[m-2]}"] \\
        \leaf[21] \&\& \leaf[22] \\
        \& \leaf[0] \ar[ul, "{\D[m+1]}"] \ar[ur, "{\E[6]}"']
    \end{tikzcd}\\\\m\neq4}} & \multicolumn{2}{C|}{\mrow[gray!20]{\begin{tikzcd}
        \& \leaf[4] \ar[dl, "{\A[5]}"'] \ar[d, "{\A[1]}"] \ar[dr, "{\A[1]}"] \\
        \leaf[21] \& \leaf[22] \& \leaf[23] \\
        \& \leaf[0] \ar[ul, "{\D[5]}"] \ar[u, "{\E[6]}"'] \ar[ur, "{\E[6]}"']
    \end{tikzcd}\\\\m=4}} \\
    \mrow{\III(m,7) \\[3.4em] \begin{tikzcd}[scale cd=.8, cramped, sep=0]
    1\&\&\&\&\&2 \\ \&2\&\cdots\&\bal{2}\&3\&4\&3\&2\&1 \\ 1 \end{tikzcd} \\[3.4em] m\geq4} & \multicolumn{2}{C}{\mrow{\begin{tikzcd}
        \& \leaf[4] \ar[dl, "{\D[6]}"'] \ar[dr, "{\D[m-2]}"] \\
        \leaf[21] \&\& \leaf[22] \\
        \& \leaf[0] \ar[ul, "{\D[m+2]}"] \ar[ur, "{\E[7]}"']
    \end{tikzcd}\\\\m\neq4}} & \multicolumn{2}{C|}{\mrow{\begin{tikzcd}
        \& \leaf[4] \ar[dl, "{\D[6]}"'] \ar[d, "{\A[1]}"] \ar[dr, "{\A[1]}"] \\
        \leaf[21] \& \leaf[22] \& \leaf[23] \\
        \& \leaf[0] \ar[ul, "{\D[6]}"] \ar[u, "{\E[7]}"'] \ar[ur, "{\E[7]}"']
    \end{tikzcd}\\\\m=4}} \\
    \mrow[gray!20]{\III(m,8) \\[3.4em] \begin{tikzcd}[scale cd=.8, cramped, sep=0]
    1\&\&\&\&\&\&\&3 \\ \&2\&\cdots\&\bal{2}\&3\&4\&5\&6\&4\&2 \\ 1
    \end{tikzcd} \\[3.4em] m\geq4} & \multicolumn{2}{C}{\mrow[gray!20]{\begin{tikzcd}
        \& \leaf[4] \ar[dl, "{\E[7]}"'] \ar[dr, "{\D[m-2]}"] \\
        \leaf[21] \&\& \leaf[22] \\
        \& \leaf[0] \ar[ul, "{\D[m+4]}"] \ar[ur, "{\E[8]}"']
    \end{tikzcd} \\\\ m\neq4}} & \multicolumn{2}{C|}{\mrow[gray!20]{\begin{tikzcd}
        \& \leaf[4] \ar[dl, "{\E[7]}"'] \ar[d, "{\A[1]}"] \ar[dr, "{\A[1]}"] \\
        \leaf[21] \& \leaf[22] \& \leaf[23] \\
        \& \leaf[0] \ar[ul, "{\D[8]}"] \ar[u, "{\E[8]}"'] \ar[ur, "{\E[8]}"']
    \end{tikzcd} \\\\ m=4}} \\
    \mrow{\III[a] \\\\ \begin{tikzcd}[scale cd=.8, cramped, sep=0]
    \&\&\&\&\bal{2} \\ 1\&2\&3\&4\&5\&4\&3\&2\&1
    \end{tikzcd}} & \multicolumn{4}{C|}{\begin{tikzcd}
        \leaf[4] \ar[d, "{\A[9]}", overlay] \\
        \leaf[2] \\
        \leaf[0] \ar[u, "{\E[7]}"', overlay]
    \end{tikzcd}} \\
    \mrow[gray!20]{\III[b] \\\\ \begin{tikzcd}[scale cd=.8, cramped, sep=0]
    \&\&4 \\ \bal{2}\&5\&8\&7\&6\&5\&4\&3\&2\&1
    \end{tikzcd}} & \multicolumn{4}{C|}{\cellcolor{gray!20}\begin{tikzcd}
        \leaf[4] \ar[d, "{\D[10]}", overlay] \\
        \leaf[2] \\
        \leaf[0] \ar[u, "{\E[8]}"', overlay]
    \end{tikzcd}} \\
\end{longtable}

  \bibliographystyle{plain}
  \renewcommand{\addcontentsline}[3]{}
  \bibliography{QuivVar}
\end{document}